\documentclass[twoside,11pt]{article}

\usepackage{url}
\usepackage{amsmath}
\usepackage[mathscr]{eucal}
\usepackage{subfig}
\usepackage{bm}
\usepackage{ntheorem}
\usepackage[all]{xy}
\usepackage{amssymb}
\usepackage{graphicx}

    \newcommand{\R}{\mathbb{R}}
    \newcommand{\F}{\ensuremath{\mathscr{F}}}
    \newcommand{\G}{\ensuremath{\mathscr{G}}}
    \renewcommand{\S}{\ensuremath{\mathscr{S}}}
    \renewcommand{\b}{\mathbf}
    \newtheorem{note}{Note}
    \newtheorem*{cor}{Corollary}
    \newtheorem{theorem}{Theorem}
    
    \newtheorem{lemma}{Lemma}
    \newtheorem*{proof2}{Proof}

\title{Undercomplete Blind Subspace Deconvolution}
\author{Zolt\'an Szab\'o, Barnab\'as P\'oczos, and Andr\'as L{\H
o}rincz\\\\
        \footnotesize Department of Information Systems, E\"{o}tv\"{o}s Lor{\'a}nd University\\
        \footnotesize P{\'a}zm{\'a}ny P. s{\'e}t{\'a}ny 1/C, Budapest H-1117, Hungary\\
        \footnotesize WWW home page: \url{http://nipg.info}\\
        \footnotesize \texttt{szzoli@cs.elte.hu, pbarn@cs.elte.hu, lorincz@inf.elte.hu}}

\begin{document}
\date{}
\maketitle
\begin{abstract}
We introduce the blind subspace deconvolution (BSSD) problem,
which is the extension of both the blind source deconvolution
(BSD) and the independent subspace analysis (ISA) tasks. We
examine the case of the undercomplete BSSD (uBSSD). Applying
temporal concatenation we reduce this problem to ISA. The
associated `high dimensional' ISA problem can be handled by a
recent technique called joint \mbox{f-decorrelation} (JFD).
Similar decorrelation methods have been used previously for kernel
independent component analysis (kernel-ICA). More precisely, the
kernel canonical correlation (KCCA) technique is a member of this
family, and, as is shown in this paper, the kernel generalized
variance (KGV) method can also be seen as a decorrelation method
in the feature space. These kernel based algorithms will be
adapted to the ISA task. In the numerical examples, we (i) examine
how efficiently the emerging higher dimensional ISA tasks can be
tackled, and (ii) explore the working and advantages of the
derived kernel-ISA methods.
\end{abstract}

\section{Introduction}\label{sec:introduction}
Independent component analysis  (ICA)
\cite{jutten91blind,comon94independent} aims to recover linearly
or non-linearly mixed independent and hidden sources. There is a
broad range of applications for ICA, such as blind source
separation, feature extraction and denoising. Particular
applications include the analysis of financial and neurobiological
data, fMRI, EEG, and MEG. For recent review concerning ICA the
reader is referred to the literature
\cite{hyvarinen01independent,cichocki02adaptive}.

Traditional ICA algorithms are  \mbox{one-dimensional} in the
sense that all sources are assumed to be independent real valued
random variables. Nonetheless, applications in which only certain
groups of sources are independent may be highly relevant in
practice. In this case, the independent sources can be
multidimensional. For instance, consider the generalization of the
cocktail-party problem, where \emph{independent groups of people}
are talking about  independent topics or more than one \emph{group
of musicians} is playing at the party. The separation task
requires an extension of ICA, which can be called independent
subspace analysis (ISA) \cite{hyvarinen00emergence},
multidimensional independent component analysis (MICA)
\cite{cardoso98multidimensional}, and group ICA
\cite{theis05blind}. We will use the first of these abbreviations
throughout this paper.

Strenuous efforts have been made to develop ISA algorithms
\cite{cardoso98multidimensional,akaho99MICA,hyvarinen00emergence,hyvarinnen06fastisa,vollgraf01multi,bach03beyond,strogbauer04least,poczos05independent1,poczos05independent2,theis05blind,theis05multidimensional,theis06towards,szabo06real,nolte06identifying}.
For the most part, ISA-related theoretical problems concern the
estimation of entropy or of mutual information. For this, the
$k$-nearest neighbors \cite{poczos05independent1} and the geodesic
spanning tree methods \cite{poczos05independent2} can be applied.
Other recent approaches seek independent subspaces via kernel
methods \cite{bach03beyond} and joint block diagonalization
\cite{theis05blind,theis06towards}.

Another extension of the original ICA task is the blind source
deconvolution (BSD) problem. Such a problem emerges, for example,
at a cocktail-party being held in an \emph{echoic room}.  Several
BSD algorithms were developed in the past. See, for example, the
review of \cite{pedersen07survey}. Like ICA, BSD has several
applications: (i) remote sensing applications; passive radar/sonar
processing \cite{macdonald05derivation,hedgepeth99expectation},
(ii) image-deblurring, image restoration \cite{vural06blind},
(iii) speech enhancement using microphone arrays, acoustics
\cite{douglas05natural,mitianoudis03audio,roan03blind,araki03fundamental},
(iv) multi-antenna wireless communications, sensor networks
\cite{akyldiz02wireless,deligianni06source}, (v) biomedical
signal---EEG, ECG, MEG, fMRI---analysis
\cite{jung00independent,glover99deconvolution,dyrholm06model},
(vi) optics \cite{kotzer98generalized}, (vii) seismic exploration
\cite{karsli06further}.

The simultaneous assumption of the two extensions, that is, ISA
combined with BSD, seems to be a more realistic model than either
of the two models alone. For example, at the cocktail-party,
groups of people or groups of musicians may form independent
source groups and echoes could be present. This task will be
called blind subspace deconvolution (BSSD). We treat the
undercomplete case (uBSSD) here. In terms of the cocktail-party
problem, it is assumed that there are more microphones than
acoustic sources. Here we note that the complete, and in
particular the overcomplete, BSSD task is challenging and as of
yet no general solution is known. We can show that temporal
concatenation turns the uBSSD task into an ISA problem. One of the
most stringent applications of BSSD could be the analysis of EEG
or fMRI signals. The ICA assumptions could be highly problematic
here, because some sources may depend one another, so an ISA model
seems better. Furthermore, the passing of information from one
area to another and the related delayed and transformed activities
may be modeled as echoes. Thus, one can argue that BSSD may fit
this important problem domain better than ICA or even ISA.

In principle, the ISA problem can be treated with the methods
listed above. However, the dimension of the ISA problem derived
from an uBSSD task is not amenable to state-of-the-art ISA
methods. According to a recent decomposition principle, the ISA
Separation Theorem \cite{szabo06separation}, the ISA task can be
divided into two consecutive steps under certain conditions: after
the application of the ICA algorithm, the ICA elements need to be
grouped.\footnote{The possibility of such a decomposition
principle was suspected by \cite{cardoso98multidimensional}, who
based his conjecture on numerical experiments. To the best of our
knowledge, a proof encompassing sufficient conditions for this
intriguing hypothesis was first published by \cite{szabo06cross}.}
The importance of this direction stems from the fact that ICA
methods can deal with problems in high dimensions. The derived ISA
task will be solved with the use of the decomposition principle
augmented by the joint \mbox{f-decorrelation} (JFD) technique
\cite{szabo06real}.

We show other ISA approaches beyond the JFD method: We adapted the
kernel canonical correlation analysis (KCCA) and the kernel
generalized variance (KGV) methods \cite{bach02kernel} to measure
the mutual dependency of multidimensional variables. One can show
that similarly to the JFD and the KCCA methods, the KGV technique
deals with nonlinear decorrelation in function spaces. We found
that they can be more precise but are limited to smaller problems.

The paper is structured as follows:
Section~\ref{sec:BSSD+ISA-model} formulates the problem domain.
Section~\ref{sec:reduction-steps} shows how to transform the uBSSD
task into an ISA task. The JFD method, which we use to solve the
derived ISA task, is the subject of Section~\ref{sec:ISA-methods}.
This section also addresses how to tailor the KCCA and KGV
kernel-ICA methods to solve the ISA problem.
Section~\ref{sec:illustrations} contains the numerical
illustrations and conclusions are drawn in
Section~\ref{sec:conclusions}.

\section{The BSSD and the ISA Model}\label{sec:BSSD+ISA-model}
The BSSD task and its special case, the ISA model, are defined in
Section~\ref{sec:BSSD-model}. Section~\ref{sec:ISA-ambiguities}
details the ambiguities of the ISA task.
Section~\ref{sec:ISA-cost} introduces some possible ISA cost
functions.

\subsection{The BSSD Equations}\label{sec:BSSD-model}
Here, we define the BSSD task. Assume that we have $M$ hidden, independent,
multidimensional \emph{components} (random variables). Suppose also that only their
casual FIR filtered mixture is available for observation:\footnote{Causal: $l\ge 0$
in $\sum_l$. FIR: the number of terms in the sum is finite.}
\begin{equation}
\b{x}(t)=\sum_{l=0}^L\b{H}_l\b{s}(t-l),\label{eq:obs}\\
\end{equation}
where
$\b{s}(t)=\left[\b{s}^1(t);\ldots;\b{s}^M(t)\right]\in\R^{Md}$ is
a vector concatenated of components $\b{s}^m(t)\in\R^{d}$. For a
given $m$, $\b{s}^m(t)$ is i.i.d. (independent and identically
distributed) in time $t$, $\b{s}^m$s are non-Gaussian, and
$I(\b{s}^1,\ldots,\b{s}^M)=0$, where $I$ stands for the mutual
information of the arguments. The total dimension of the
components is \mbox{$D_s:=Md$}, the dimension of the observation
$\b{x}$ is $D_x$. Matrices $\b{H}_l\in \R^{D_x\times D_s}$
$(l=0,\ldots,L)$ describe the mixing, these are the \emph{mixing
matrices}. Without any loss of generality it may be assumed that
$E[\b{s}]=\b{0}$, where $E$ denotes the expectation value. Then
$E[\b{x}]=\b{0}$ holds, as well. The goal of the BSSD problem is
to estimate the original source $\b{s}(t)$ by using observations
$\b{x}(t)$ only.

The case $L=0$ corresponds to the ISA task, and if $d=1$ also
holds then the ICA task is recovered. In the BSD task $d=1$ and
$L$ is a non-negative integer. $D_x>D_s$ is the
\emph{undercomplete}, $D_x=D_s$ is the \emph{complete}, and
$D_x<D_s$ is the \emph{overcomplete} task. Here, we treat the
undercomplete BSSD (uBSSD) problem. We will transform the uBSSD
task to undercomplete ISA (uISA) or to complete ISA. From now on
they both will be called ISA.

\begin{note}
Mixing matrices $\b{H}_l$ ($0 \le l \le L$) have a one-to-one
mapping to polynomial matrix\footnote{$\b{H}[z]$ is also known as
\emph{channel matrix} or \emph{transfer function} in the
literature.}
\linebreak[4] \mbox{$\b{H}[z]:=\sum_{l=0}^L\b{H}_lz^{-l}\in\R[z]^{D_x\times
D_s}$}, where $z$ is the time-shift operation, that is
\mbox{$(z^{-1}\b{u})(t)=\b{u}(t-1)$}. $\b{H}[z]$ may be regarded
as an operation that maps $D_s$-dimensional series to
$D_x$-dimensional series. Equation~\eqref{eq:obs} can be written
as $\b{x}=\b{H}[z]\b{s}$.
\end{note}

\begin{note}\label{note:FIR-inverse-exists}
It can be shown \cite{rajagopal03multivariate} that in the uBSSD
task $\b{H}[z]$ has a polynomial matrix left inverse
$\b{W}[z]\in\R[z]^{D_x\times D_s}$ with probability 1, under mild
conditions. In other words, for these polynomial matrices
$\b{W}[z]$ and $\b{H}[z]$, $\b{W}[z]\b{H}[z]$ is the identity
mapping. The mild condition is as follows: Coefficients of
polynomial matrix $\b{H}[z]$, that is, random matrix
$[\b{H}_0;\ldots;\b{H}_L]$ is drawn from a continuous
distribution. Under this condition, hidden source $\b{s}(t)$ can
be estimated by a \emph{suitable} causal FIR filtered form of
observation $\b{x}(t)$.
\end{note}

For the uBSSD task it is assumed that $\b{H}[z]$ has a polynomial matrix left
inverse. For the uISA and ISA tasks it is supposed that \emph{mixing matrix}
$\b{H}_0\in\R^{D_x\times D_s}$ has full column rank, that is its rank is $D_s$.

\subsection{Ambiguities of the ISA Model}\label{sec:ISA-ambiguities}
Because the uBSSD task will be reduced to ISA, it is important to see the ambiguities of the ISA task. First, the
complete ISA problem ($L=0, D_x=D_s$) is presented, the undercomplete ISA will be treated later.

The identification of the ISA model is ambiguous. However, the
ambiguities are simple \cite{theis04uniqueness1}: hidden
multidimensional components can be determined up to permutation
and up to invertible transformation within the subspaces.
Ambiguities within the subspaces can be weakened. Namely, because
of the invertibility of mixing matrix
\mbox{$\b{H}[z]=\b{H}_0\in\R^{D_s\times D_s}$}, it can be assumed
without any loss of generality that both the sources and the
observation are \emph{white}, that is,
\begin{align*}
E[\b{s}]&=\b{0},cov\left[\b{s}\right]=\b{I}_{D_s},\\
E[\b{x}]&=\b{0},cov\left[\b{x}\right]=\b{I}_{D_x},
\end{align*}
where $\b{I}_{D_s}$ is the \mbox{$D_s$-dimensional} identity
matrix and $cov$ is the covariance matrix. It then follows that
the mixing matrix $\b{H}_0$ and thus the \emph{demixing matrix}
$\b{W}=\b{H}_0^{-1}$ are orthogonal:
\[
    \b{I}_{D_s} =cov[\b{x}]=E\left[\b{x}\b{x}^*\right]=\b{H}_0E\left[\b{s}\b{s}^*\right]\b{H}_0^*=\b{H}_0\b{I}_{D_s}\b{H}_0^*=\b{H}_0\b{H}_0^*,
\]
where $*$ denotes transposition. In sum, $\b{H}_0,\b{W}\in\mathscr{O}^{D_s}$, where
$\mathscr{O}^{D_s}$ denotes the set of \mbox{$D_s$-dimensional} orthogonal matrices.
Now, $\b{s}^m$ sources are determined up to permutation and orthogonal
transformation.

In order to transform the undercomplete ISA task into a complete ISA task with white observations let
\mbox{$\b{C}:=cov[\b{x}]=E[\b{x}\b{x}^*]=\b{H}_0\b{H}_0^*\in\R^{D_x\times D_x}$} denote the covariance matrix of the
observation. Rank of $\b{C}$ is $D_s$, since the rank of matrix $\b{H}_0$ is $D_s$ according to our assumptions. Matrix
$\b{C}$ is symmetric ($\b{C}=\b{C}^*$), thus it can be decomposed as follows: $\b{C}=\b{U}\b{D}\b{U}^*$, where
$\b{U}\in\R^{D_x\times D_s}$, and the columns of matrix $\b{U}$ are orthogonal, that is, $\b{U}^*\b{U}=\b{I}_{D_s}$.
Furthermore, the rank of diagonal matrix $\b{D}\in\R^{D_s\times D_s}$ is $D_s$. The principal component analysis can
provide a decomposition in the desired form. Let $\b{Q}:=\b{D}^{-1/2}\b{U}^*\in\R^{D_s\times D_x}$. Then the original
observation $\b{x}$ can be modified to $\b{x}':=\b{Q}\b{x}=\b{Q}\b{H}_0\b{s}\in\R^{D_s}$. The resulting $\b{x}'$ is
white and can be regarded as the observation of a \emph{complete} ISA task having mixing matrix
$\b{Q}\b{H}_0\in\mathscr{O}^{D_s}$.

\subsection{ISA Cost Functions}\label{sec:ISA-cost}
After the whitening procedure (Section~\ref{sec:ISA-ambiguities}), the ISA task can
be viewed as the minimization of the mutual information between the estimated
components on the orthogonal group:
\begin{equation}
    J_{I}(\b{W}):=I\left(\b{y}^1,\ldots,\b{y}^M\right),\label{eq:ISA-cost-I}
\end{equation}
where \mbox{$\b{y}=\b{W}\b{x}$}, \mbox{$\b{y}=\left[\b{y}^1;\ldots;\b{y}^M\right]$}, $\b{y}^m \in \R^d$, and
$\b{W}\in\mathcal{O}^D$. This formulation of the ISA task serves us in Section~\ref{sec:studied-ISA-costs}, where we
estimate the dependencies of the multidimensional variables.

The ISA task can be rewritten into the minimization of the sum of
Shannon's multidimensional differential entropies
\cite{poczos05independent1}:
\begin{equation}
J_H(\b{W}):=\sum_{m=1}^MH\left(\b{y}^m\right),
\label{eq:ISA-cost-sumH}
\end{equation}
where \mbox{$\b{y}=\b{W}\b{x}$}, \mbox{$\b{y}=\left[\b{y}^1;\ldots;\b{y}^M\right]$}, $\b{y}^m \in \R^d$,
$\b{W}\in\mathcal{O}^D$.

\begin{note}
Until now, we formulated the ISA task by means of the entropy or
the mutual information of multidimensional random variables, see
Equations~\eqref{eq:ISA-cost-I} and \eqref{eq:ISA-cost-sumH}.
However, any algorithm that treats mutual information between
\emph{1-dimensional} random variables can also be sufficient. This
statement is based on the considerations below. Well-known
identities of mutual information and entropy expressions
\cite{cover91elements} show that the minimization of cost function
\begin{equation*}
J_{H,I}(\b{W}):= \sum_{m=1}^M \sum_{i=1}^d H(y^m_i)-\sum_{m=1}^M I(y^m_1,\ldots,y^m_d),
\end{equation*}
or that of
\begin{equation*}
    J_{I,I}(\b{W}):=I\left(y_1^1,\ldots,y_d^M\right)-\sum_{m=1}^M I\left(y_1^m,\ldots,y_d^m\right)
\end{equation*}
can also solve the ISA task. Here, \mbox{$\b{y}=\b{W}\b{x}$} is
the estimated ISA source, where $\b{x}\in\R^D$ is the whitened
observation in the ISA model. $\b{W}\in\mathscr{O}^D$ is the
estimated ISA demixing matrix, and in
\mbox{$\b{y}=\left[\b{y}^1;\ldots;\b{y}^M\right]\in\R^D$} the
$\b{y}^m\in\R^d$, $m=1, \ldots, M$, represent the estimated
components with coordinates $y^m_i\in\R$. The first term of both
cost functions $J_{H,I}$ and $J_{I,I}$ is an ICA cost function.
Thus, these first terms can be fixed by means of ICA
preprocessing.\footnote{From the algorithmic point of view,
\emph{any} ICA algorithm that minimizes cost function
$I(y^1_1,\ldots,y^M_d)$ suits the ICA preprocessing step.} In this
case, if the Separation Theorem holds (for details see
Section~\ref{sec:ISA2ICA}), then term $\sum_{m=1}^M
I(y_1^m,\ldots,y_d^m)$ implies that the maximization of the sum of
mutual information between \emph{\mbox{1-dimensional}} random
variables \emph{within} the subspaces is sufficient for solving
the ISA task.
\end{note}

\section{Reduction Steps}\label{sec:reduction-steps}
Here we show that the direct search for inverse FIR filter can be
circumvented (Note~\ref{note:FIR-inverse-exists}). Namely,
temporal concatenation reduces the uBSSD task to an (u)ISA problem
(Section~\ref{sec:uBSSD2(u)ISA}). Our earlier results will allow
further simplifications. We will reduce the ISA task to an ICA
task plus a search for optimal permutation of the ICA coordinates.
This decomposition principle will be elaborated in
Section~\ref{sec:ISA2ICA} by means of the Separation Theorem.

\subsection{Reduction of uBSSD to (u)ISA}\label{sec:uBSSD2(u)ISA}
We reduce the uBSSD task to an ISA problem. The BSD literature
provides the basis for our reduction; \cite{fevotte03unified} use
temporal concatenation in their work. This method can be extended
to multidimensional $\b{s}^m$ components in a natural fashion:

Let $L'$ be such that
\begin{equation}\label{eq:L'-ineq-implicit}
D_xL'\ge D_s(L+L')
\end{equation}
is fulfilled. Such $L'$ exists due to the undercomplete assumption $D_x>D_s$:
\begin{equation}\label{eq:L'-ineq-explicit}
L'\ge\left\lceil\frac{D_sL}{D_x-D_s}\right\rceil.
\end{equation}
This choice of $L'$ guarantees that the reduction gives rise to an (under)complete ISA task: let $x_m(t)$ denote the
$m^{th}$ coordinate of observation $\b{x}(t)$ and let the matrix \mbox{$\b{H}_l\in\R^{D_x\times Md}$} be decomposed
into  $1\times d$ sized blocks. That is, $\b{H}_l={[\b H^{ij}_l]}_{i=1..D_x,j=1..M}$ $(\b H^{ij}_l\in\R^{1\times d})$,
where $i$ and $j$ denote row and column indices, respectively. Using notations
\begin{align*}
\b{S}^m(t)&:=[\b{s}^m(t); \b{s}^m(t-1);\ldots; \b{s}^m(t-(L+L')+1)]\in\R^{d(L+L')},\\
\b{X}^m(t)&:=[x_m(t); x_m(t-1);\ldots; x_m(t-L'+1)]\in\R^{L'},\\
\b{S}(t)&:=[\b{S}^1(t); \ldots; \b{S}^M(t)]\in\R^{Md(L+L')=D_s(L+L')},\\
\b{X}(t)&:=[\b{X}^1(t); \ldots; \b{X}^{D_x}(t)]\in\R^{D_xL'},\\
\b{A}^{ij}&:= \left[\begin{array}{cccccc}\b{H}^{ij}_0&\ldots&\b{H}^{ij}_L& \b{0} &\ldots&\b{0}\\&\ddots&&\ddots&&\\&&\ddots&&\ddots\\\b{0}&\ldots&\b{0}&\b{H}^{ij}_0&\ldots&\b{H}^{ij}_L\end{array}\right]\in\R^{L'\times d(L+L')},\\
\b{A}&:=[\b{A}^{ij}]_{i=1..D_x,j=1..M}\in\R^{D_xL'\times Md(L+L')=D_xL'\times D_s(L+L')},
\end{align*}
model
\begin{equation}\label{eq:uBSSD-reduced-model}
\b{X}(t)=\b{A}\b{S}(t)
\end{equation}
can be obtained. Here, $\b{s}^m(t)$s are i.i.d.\ in time $t$, they are independent
for different $m$ values, and Equation~\eqref{eq:L'-ineq-implicit} holds for $L'$. Thus,
\eqref{eq:uBSSD-reduced-model} is either an undercomplete or a complete ISA task,
depending on the relation of the l.h.s and the r.h.s of \eqref{eq:L'-ineq-implicit}:
the task is complete if the two sides are equal. The number of the components and
the dimension of the components in task \eqref{eq:uBSSD-reduced-model} are
\mbox{$M(L+L')$} and $d$, respectively.

If we end up with an undercomplete ISA problem in
\eqref{eq:uBSSD-reduced-model} then it can be reduced to a
complete one, as was shown in Section~\ref{sec:ISA-ambiguities}.
Thus, choosing the minimal value for $L'$ in
\eqref{eq:L'-ineq-explicit}, the dimension of the obtained ISA
task is
\begin{equation}
D_{\text{ISA}}:=D_s(L+L')=D_s\left(L+\left\lceil\frac{D_sL}{D_x-D_s}\right\rceil\right).\label{eq:D-ISA}
\end{equation}
Taking into account the ambiguities of the ISA task
(Section~\ref{sec:ISA-ambiguities}), the original $\b{s}^m$
components will occur $L+L^{'}$ times and up to orthogonal
transformations. As a result, in the ideal case, our estimations
are as follows
\[
\hat{\b{s}}^m_k:=\b{F}^m_k\b{s}^m\in\R^d,
\]
where $k=1,\ldots,L+L'$, $\b{F}^m_k\in\mathscr{O}^d$.

\subsection{Reduction of ISA to ICA}\label{sec:ISA2ICA}
The Separation Theorem \cite{szabo06separation} conjectured by
\cite{cardoso98multidimensional} allows one to decompose the
solution of the ISA problem, under certain conditions, into 2
steps: In the first step, ICA estimation is executed by minimizing
$I(y^1_1,\ldots,y^M_d)$. In the second step, the ICA elements are
grouped by finding an optimal permutation. This principle will be
formalized in Section~\ref{sec:ISA-sep-T}.
Section~\ref{sec:suff-conds-of-sep-T} provides sufficient
conditions for the theorem.

\subsubsection{The ISA Separation Theorem}\label{sec:ISA-sep-T}
We state the ISA Separation Theorem for components having possibly different $d_m$ dimensions:
\begin{theorem}[Separation Theorem for ISA]
Let $\b{y}=\left[y_1;\ldots;y_D\right]=\b{W}\b{x}\in\R^D$, where \mbox{$\b{W}\in \mathcal{O}^D$}, $\b{x}\in\R^D$ is the
whitened observation of the ISA model, and \mbox{$D=\sum_{m=1}^Md_m$}. Let $\S^{d_m}$ denote the surface of the
\mbox{$d_m$-dimensional} unit sphere, that is \mbox{$\S^{d_m}:=\{\b{w}\in\R^{d_m}:\sum_{i=1}^{d_m}w_i^2=1\}$}.

Presume that the $\b{u}:=\b{s}^m\in\R^{d_m}$ sources
$(m=1,\ldots,M)$ of the ISA model satisfy condition
\begin{equation}
H\left(\sum_{i=1}^{d_m}
w_iu_i\right)\ge\sum_{i=1}^{d_m}w_i^2H\left(u_i\right),
\forall\b{w}\in \S^{d_m},\label{eq:suff}
\end{equation}
and that the ICA cost function
$J_{\text{ICA}}(\b{W})=\sum_{i=1}^DH(y_i)$ has minimum over the
orthogonal matrices in $\b{W}_{\mathrm{ICA}}$. Then it is
sufficient to search for the solution to the ISA task as a
permutation of the solution of the ICA task. Using the concept of
demixing matrices, it is sufficient to explore forms
\[
    \b{W}_{\mathrm{ISA}}=\b{P}\b{W}_{\mathrm{ICA}},\label{eq:Wform}
\]
where $\b{P}\in\R^{D\times D}$ is a permutation matrix to be
determined and $\b{W}_\mathrm{ISA}$ is the ISA demixing matrix.
\end{theorem}
The proof of the theorem is presented in
Appendix~\ref{sec:RISA-sep-theorem}. It is intriguing that if
\eqref{eq:suff} is satisfied then the simple decomposition
principle provides the \emph{global} minimum of
\eqref{eq:ISA-cost-I}. In the literature on joint block
diagonalization (JBD) \cite{meraim04algorithms} have put forth a
similar \emph{conjecture} recently. According to this conjecture,
for quadratic cost function, if Jacobi optimization is applied,
the block-diagonalization of the matrices can be found by the
optimization of permutations following the joint diagonalization
of the matrices. ISA solutions formulated within the JBD framework
\cite{theis05blind,theis05multidimensional,theis06towards,szabo06real}
make efficient use of this idea in practice. \cite{theis06towards}
could justify this approach for \emph{local} minimum points.

\subsubsection{Sufficient Conditions of the ISA Separation Theorem}\label{sec:suff-conds-of-sep-T}
The question of which types of sources satisfy the Separation
Theorem is open. Equation~\eqref{eq:suff} provides only a
sufficient condition. Below, we list sources $\b{s}^m$ that
satisfy \eqref{eq:suff}. Details and the extension of the
Separation Theorem for complex variables can be found in a
technical report of \cite{szabo06separation}.

\begin{enumerate}
    \item
    Assume that variables $\b{u}=\b{s}^m$ satisfy the so-called w-EPI condition (EPI is shorthand for the \emph{entropy power inequality} \cite{cover91elements}), that is,
    \begin{equation}
            e^{2H\left(\sum_{i=1}^dw_iu_i\right)}\ge \sum_{i=1}^d e^{2H(w_iu_i)}, \forall\b{w}\in \S^d.\label{eq:w-EPI-first}
    \end{equation}
    Then inequality \eqref{eq:suff} holds for these variables too. The proof can be found in Lemma~\ref{lem:suff} of Appendix~\ref{sec:RISA-sep-theorem}.
    \item
        The \eqref{eq:w-EPI-first} w-EPI condition is valid
        \begin{enumerate}
            \item
                for spherically symmetric or shortly spherical variables
                \cite{fang90symmetric}. The distribution of such variables is invariant for orthogonal transformations.\footnote{In the ISA
                task the non-degenerate affine transformations of spherical variables, the so called elliptical
                variables, do not provide valuable generalizations due to the ambiguities of the ISA task.} Sketch of the proof ($\b{u}=\b{s}^m$):
                the w-EPI condition concerns projections to unit vectors. For spherical
                variables, the distribution and thus the entropy of these projections are independent of $\b{w}\in\S^d$.
                Because  $e^{2H(w_iu_i)}=e^{2H(u_i)}w_i^2$ and $\b{w}\in \S^d$,
                the \mbox{w-EPI} is satisfied with equality $\forall \b{w}\in \S^d$.  $\square$
            \item
                for \mbox{2-dimensional} variables invariant to $90^{\circ}$
                rotation. Under this condition, density function $h$ of component $\b{s}^m$ is subject to the
                following invariance
                    \begin{equation*}
                        h(u_1,u_2)=h(-u_2,u_1)=h(-u_1,-u_2)=h(u_2,-u_1)\quad\left(\forall
                        \b{u}\in\R^2\right).
                    \end{equation*}

                Sketch of the proof ($\b{u}=\b{s}^m$): Assume that function
                $f:\S^2\ni\b{w}\mapsto H\left(\sum_{i=1}^dw_iu_i\right)$ has global minimum on set
                 $\S^2\cap\{\b{w}\ge \b{0}\}$.\footnote{Relation $\b{w}\ge \b{0}$ concerns each coordinate.}
                Let this minimum be at $\b{w}_m\in\R^2$. Then, the $90^{\circ}$ invariance warrants that function $f$
                take its global minimum also on  $\b{w}_m^{\bot}\in\R^2$, which is perpendicular to  $\b{w}_m$. Let
                $(\b{C}^m)^*=[\b{w}_m,\b{w}_m^{\bot}]\in\mathcal{O}^2$. Now, we can estimate variables $\b{C}^m\b{s}^m$. This
                is sufficient because the ISA solution is ambiguous up to orthogonal transformations within
                each subspace.  $\square$

                A special case of this requirement is invariance
                to permutation and sign changes
                \begin{equation*}
                    h(\pm u_1,\pm u_2)=h(\pm u_2,\pm u_1).
                \end{equation*}
                In other words, there exists a function
                $g:\R^2\rightarrow\R$,
                which is symmetric in its variables and
                \begin{equation*}
                    h(\b{u})=g(|u_1|,|u_2|).
                \end{equation*}
                Special cases within this family are distributions
                \begin{equation*}
                    h(\b{u})=g\left(\sum_i|u_i|^p\right)\quad(p>0),
                \end{equation*}
                which are constant over the spheres of $L^p$-space. They are called\ $L^p$ spherical variables which, for $p=2$,
                corresponds to spherical variables.
            \item
                for certain weakly dependent variables: \cite{takano95inequalities} has determined sufficient conditions when EPI
                holds.\footnote{The constraint of $d=2$ may be generalized to higher
                dimensions. We are not aware of such generalizations.} If the EPI property is satisfied on unit sphere
                $\S^d$, then the ISA Separation Theorem holds (Lemma~\ref{lem:suff}).
        \end{enumerate}
\end{enumerate}
These results are summarized schematically in Table~\ref{tab:suffcond-summary}.
\begin{table}
  \[\xymatrixcolsep{-1.7cm}
  \xymatrix{
  &\txt{invariance to $90^{\circ}$ rotation ($d=2$)}\ar@{=>}[ddd]\ar[dr]|-{\text{specially}}&\\
  &&\text{invariance to sign and permutation}\ar[d]|-{\text{specially}}\\
  &&\text{$L^p$ spherical ($p>0$)}\\
  \txt{Takano's dependency\\($d=2$)} \ar@{=>}[r] & \text{w-EPI}\ar@{=>}[d] & \txt{spherical symmetry}\ar@{=>}[l]\ar[u]|-{\text{generalization for $d=2$}}\\
  &\txt{Equation~\eqref{eq:suff}: sufficient\\ for the ISA Separation Theorem}& }
  \]
  \caption{Sufficient conditions for the ISA Separation Theorem.}\label{tab:suffcond-summary}
\end{table}

\section{ISA Methods}\label{sec:ISA-methods}
We showed how to convert the uBSSD task to an ISA task in
Section~\ref{sec:uBSSD2(u)ISA}. In the following we will present
methods that can solve the ISA task. In
Section~\ref{sec:studied-ISA-costs} we treat estimations of the
mutual information of the ISA cost functions in
Section~\ref{sec:ISA-cost}. Methods that can optimize these cost
functions are elaborated in Section~\ref{sec:Opt-of-ISA-costs}. We
also present here the pseudocode of the procedures studied. In
Section~\ref{sec:diff-or-unknown-dims} we review methods that can
treat non-equal or unknown component dimensions. In what follows,
and in accordance with \eqref{eq:obs}, let $\b{x}\in\R^D$ denote
the whitened observation, while
\mbox{$\b{y}=[\b{y}^1;\ldots;\b{y}^M]=\b{W}\b{x}\in\R^D$}
($\b{W}\in\mathscr{O}^D$) and $\b{y}^m\in\R^d$
\mbox{$(m=1,\ldots,M)$}  stand for the estimated source and its
components in the ISA task, respectively.

\subsection{Dependency Estimations}\label{sec:studied-ISA-costs}
Here we introduce two dependency estimators. First, in
Section~\ref{sec:JFD} we describe a decorrelation method that uses
a set of functions jointly. This method is called joint
\mbox{f-decorrelation} (JFD) method \cite{szabo06real}. Our second
technique (Section~\ref{sec:kernelISA-methods}) generalizes
earlier kernel-ICA methods for the ISA task. The motivation for
this latter method is the efficiency and precision of kernel-ICA
methods in finding independent components \cite{bach02kernel}. Our
experiences are similar with kernel-ISA methods, see
Section~\ref{sec:simulations}. We found that kernel-ISA methods
need more computations, but can provide more precise solutions
than the JFD technique.

\subsubsection{The JFD Method}\label{sec:JFD} The JFD method
estimates the hidden $\b{s}^m$ components through the
decorrelation over a function set $\F(\ni\b{f})$
\cite{szabo06real}. Formally, let the empirical $\b{f}$-covariance
matrix of $\b{y}(t)$ and $\b{y}^m(t)$ for function
$\b{f}=[\b{f}^1;\ldots;\b{f}^M]\in\F$ over $t=1,\ldots,T$ be
denoted by
\begin{align*}
     {\Sigma}(\b{f},T,\b{W})&=\widehat{cov}\left[\b{f}\left(\b{y}\right),\b{f}\left(\b{y}\right)\right]=\\
     &=
     \frac{1}{T}\sum_{t=1}^T\left\{\b{f}[\b{y}(t)]-\frac{1}{T}\sum_{k=1}^T\b{f}[\b{y}(k)]\right\}\left\{\b{f}[\b{y}(t)]-\frac{1}{T}\sum_{k=1}^T\b{f}[\b{y}(k)]\right\}^*,\nonumber\\
     {\Sigma}^{i,j}(\b{f},T,\b{W})&=\widehat{cov}\left[\b{f}^i\left(\b{y}^i\right),\b{f}^j\left(\b{y}^j\right)\right]=\\
     &=\frac{1}{T}\sum_{t=1}^T\left\{\b{f}^i[\b{y}^i(t)]-\frac{1}{T}\sum_{k=1}^T\b{f}^i[\b{y}^i(k)]\right\}\left\{\b{f}^j[\b{y}^j(t)]-\frac{1}{T}\sum_{k=1}^T\b{f}^j[\b{y}^j(k)]\right\}^*.\nonumber
\end{align*}
Then, the joint decorrelation on $\F$ can be formulated as the minimization of cost function
\begin{equation}
    J_{\text{JFD}}(\b{W}):=\sum_{\b{f}\in\F}\left\|\b{N}\circ{\Sigma}(\b{f},T,\b{W})\right\|^2_F.\label{eq:JFD-cost}
\end{equation}
Here: (i) $\b{W}\in\mathscr{O}^D$, (ii) $\F$ denotes a set of $\R^D\rightarrow \R^D$ functions, and each function acts
on each coordinate separately, (iii) $\circ$ denotes the point-wise multiplication, called the Hadamard-product, (iv)
$\b{N}$ masks according to the subspaces, $\b{N}:=\b{E}_D-\b{I}_M\otimes \b{E}_d$, where all elements of matrix
\mbox{$\b{E}_D\in\R^{D\times D}$} and \mbox{$\b{E}_d\in\R^{d\times d}$} are equal to 1, $\otimes$ is the
Kronecker-product, (v) $\left\|\cdot\right\|^2_F$ denotes the square of the Frobenius norm, that is, the sum of the
squares of the elements.

Cost function \eqref{eq:JFD-cost} can be interpreted as follows: for \emph{any} function $\b{f}^m:\R^d\rightarrow\R^d$
that acts on independent variables $\b{y}^m$ $(m=1,\ldots,M)$ the variables $\b{f}^m(\b{y}^m)$ remain independent.
Thus, covariance matrix ${\Sigma}(\b{f},T,\b{W})$ of variable
$\b{f}(\b{y})=\left[\b{f}^1\left(\b{y}^1\right);\ldots;\b{f}^M\left(\b{y}^M\right)\right]$ is block-diagonal.
Independence of estimated sources $\b{y}^m$ is gauged by the uncorrelatedness on the function set $\F$. Thus, the
non-block-diagonal portions (${\Sigma}^{i,j}(\b{f},T,\b{W})$, $i\ne j$) of covariance matrices
${\Sigma}(\b{f},T,\b{W})$ are punished. This principle is expressed by the term
$\left\|\b{N}\circ{\Sigma}(\b{f},T,\b{W})\right\|^2_F$.

\subsubsection{Kernel-ISA Methods}\label{sec:kernelISA-methods} Two
alternatives for the ISA cost function of \eqref{eq:JFD-cost} are
presented. They estimate the mutual information based ISA cost
defined in \eqref{eq:ISA-cost-I} via kernels: the KCCA and KGV
kernel-ICA methods of \cite{bach02kernel} are extended to the ISA
task. The original methods estimate pair-wise independence between
\emph{1-dimensional} random variables.\footnote{We note that if
our observations are generated by an ISA model then---unlike in
the ICA task when $d=1$---pairwise independence is \emph{not}
equivalent to mutual independence
\cite{comon94independent,poczos05independent2}. Nonetheless,
according to our numerical experiences it is an efficient
approximation in many situations.} The extension to the
multidimensional case is straightforward, the arguments of the
kernels can be modified to multidimensional variables and the
derivation of \cite{bach02kernel} can be followed. The main steps
are provided below for the sake of completeness. The resulting
expressions can be used for the estimation of dependence between
multidimensional random variables. The performance of these simple
extensions on the related ISA applications is shown in
Section~\ref{sec:kernelISA-vs-JFD}.

\paragraph{The KCCA Method}\label{sec:KCCA} First, the
2-variable-case is treated and then it will be generalized to many
variables. \subparagraph{2-variable-case} Assume that the mutual
dependence of two random variables $\b{u}\in\R^{d_1}$ and
$\b{v}\in\R^{d_2}$ has to be measured. Let positive semi-definite
kernels \mbox{$k^{\b{u}}(\cdot,\cdot): \R^{d_1}\times \R^{d_1} \to
\R $}, and $k^{\b{v}}(\cdot,\cdot): \R^{d_2}\times \R^{d_2} \to \R
$ be chosen in the respective spaces. Let $\F^{\b{u}}$ and
$\F^{\b{v}}$ denote the reproducing kernel Hilbert spaces (RKHS)
\cite{aronszajn50theory,wahba99support,scholkopf99advances}
associated with the kernels. Here, $\F^{\b{u}}$ and $\F^{\b{v}}$
are function spaces having elements that perform mappings
$\R^{d_1}\rightarrow\R$ and $\R^{d_2}\rightarrow\R$, respectively.
Then the mutual dependence between $\b{u}$ and $\b{v}$ can be
measured, for instance, by the following expression:

\begin{align*}
J^{*}_{\text{KCCA}}(\b{u},\b{v},\F^{\b{u}},\F^{\b{v}})&:= \sup_{g \in \F^{\b{u}}, h \in \F^{\b{v}}}
\textrm{corr}[g(\b{u}),h(\b{v})],
\end{align*}
where $corr$ denotes correlation.

The value of $J^*_{\text{KCCA}}$ can be estimated empirically: assume that we have
$T$ samples both from $\b{u}$ and from $\b{v}$. These samples are
$\b{u}_1,\ldots,\b{u}_T\in\R^{d_1}$ and $\b{v}_1,\ldots,\b{v}_T\in\R^{d_2}$. Then,
using notations \mbox{$\bar{g}:=\frac{1}{T}\sum \limits_{k=1}^T g(\b{u}_k)$},
$\bar{h}:=\frac{1}{T}\sum \limits_{k=1}^T g(\b{v}_k)$, the empirical estimation of
$J^*_{\text{KCCA}}$ could be the following:
\[
J^{*,emp}_{\text{KCCA}}(\b{u},\b{v},\F^{\b{u}},\F^{\b{v}}):= \sup_{g \in \F^{\b{u}}, h \in \F^{\b{v}}} \frac{\frac{1}{T}\sum_{t=1}^T
[g(\b{u}_t)-\bar{g}][h(\b{v}_t)-\bar{h}]} {\sqrt{\frac{1}{T}\sum_{t=1}^T[g(\b{u}_t)-\bar{g}]^2} \sqrt{\frac{1}{T}\sum_{t=1}^T[
h(\b{v}_t)-\bar{h}]^2} }.
\]
However, it is worth including some regularization for
$J^{*}_{\text{KCCA}}$ \cite{fukumizu07statistical}, therefore
$J^{*}_{\text{KCCA}}$ is modified to
\begin{equation}
J_{\text{KCCA}}(\b{u},\b{v},\F^{\b{u}},\F^{\b{v}}):=\sup_{g \in \F^{\b{u}}, h \in \F^{\b{v}}}
\frac{\textrm{cov}[g(\b{u}),h(\b{v})]}{\sqrt{\textrm{var}\left[g(\b{u})\right]+\kappa \left\|g\right\|^2_{\F^{\b{u}}}}
\sqrt{\textrm{var}\left[h(\b{v})\right]+\kappa \left\|h\right\|^2_{\F^{\b{v}}}}},\label{eq:J-KCCA-reg}
\end{equation}
where  expression `$\textrm{var}$' stands for variance, $\kappa>0$
is the regularization parameter, $\left\|.\right\|^2_{\F^{\b{u}}}$
and $\left\|.\right\|^2_{\F^{\b{v}}}$ denote the RKHS norm of
their arguments in $\F^{\b{u}}$ and $\F^{\b{v}}$, respectively.
Now, expanding the denominator up to second order in $\kappa$,
setting the expectation value of the samples to zero in the
respective RKHSs, and using the notation $\kappa_2:=\frac{\kappa
T}{2}$ \cite{bach02kernel}, the empirical estimation of
\eqref{eq:J-KCCA-reg} is

\begin{equation}
\hat{J}^{emp}_{\text{KCCA}}(\b{u},\b{v},\F^{\b{u}},\F^{\b{v}}) = \sup_{\b{c}_1 \in \R^T, \b{c}_2 \in \R^T}
\frac{ \b{c}_1^* \widetilde{\b{K}}^{\b{u}}\widetilde{\b{K}}^{\b{v}} \b{c}_2 } {\sqrt{\b{c}_1^*
\left(\widetilde{\b{K}}^{\b{u}}+\kappa_2\b{I}_T\right)^2\b{c}_1} \sqrt{\b{c}_2^*\left(
\widetilde{\b{K}}^{\b{v}}+\kappa_2\b{I}_T\right)^2 \b{c}_2}}, \label{eq:J-KCCA-reg-dual}
\end{equation}
where $\widetilde{\b{K}}^{\b{u}}$,
$\widetilde{\b{K}}^{\b{v}}\in\R^{T\times T}$ are the so-called
centered kernel matrices: These matrices are derived from  kernel
matrices
$\b{K}^{\b{u}}=[k(\b{u}_i,\b{u}_j)]_{i,j=1,\ldots,T},\b{K}^{\b{v}}=[k(\b{v}_i,\b{v}_j)]_{i,j=1,\ldots,T}\in\R^{T\times
T}$, as is described below. Let $\b{1}_T\in\R^T$ denote a vector
whose all elements are equal to $1$ and let
\mbox{$\b{H}:=\b{I}_T-\frac{1}{T}\b{1}_T\b{1}_T^*\in\R^{T\times
T}$} denote the so-called $T$-dimensional centering matrix. Then
$\widetilde{\b{K}}^{\b{u}}:= \b{H} \b{K}^{\b{u}}\b{H}$,
$\widetilde{\b{K}}^{\b{v}}:= \b{H} \b{K}^{\b{v}}\b{H}$.

Computing the stationary points of $\hat{J}^{\,emp}_{\text{KCCA}}$ in \eqref{eq:J-KCCA-reg-dual}, that is, setting
$\b{0}=\frac{\partial \hat{J}^{\,emp}_{\text{KCCA}}}{\partial \b{c}}$, the resulting task is to solve a
\emph{generalized eigenvalue problem} of the form $\b{C}\bm{\xi}=\mu\b{D}\bm{\xi}$:
$$
\begin{pmatrix}
  (\widetilde{\b{K}}^{\b{u}}+\kappa_2\b{I}_T)^2 & \widetilde{\b{K}}^{\b{u}} \widetilde{\b{K}}^{\b{v}} \\
  \widetilde{\b{K}}^{\b{v}} \widetilde{\b{K}}^{\b{u}} & (\widetilde{\b{K}}^{\b{v}}+\kappa_2\b{I}_T)^2 \\
\end{pmatrix}
\begin{pmatrix}
  \b{c}_1 \\
  \b{c}_2 \\
\end{pmatrix}=(1+\gamma)\begin{pmatrix}
                    (\widetilde{\b{K}}^{\b{u}}+\kappa_2\b{I}_T)^2 & \b{0} \\
                    \b{0} &  (\widetilde{\b{K}}^{\b{v}}+\kappa_2\b{I}_T)^2\\
                  \end{pmatrix}
                  \begin{pmatrix}
  \b{c}_1 \\
  \b{c}_2 \\
\end{pmatrix},
$$
where the objective is to maximize $\gamma:=\b{c}_1^* \widetilde{\b{K}}^{\b{u}}\widetilde{\b{K}}^{\b{v}} \b{c}_2$. Our
task is to estimate $\hat{J}^{\,emp}_{\text{KCCA}}$, the maximum of $\gamma$.

\subparagraph{Generalization for many variables} The KCCA method
can be generalized for more than two random variables and can be
used to measure pair-wise dependence: Let us introduce the
following notations: Let $\b{y}^1 \in \R^{d_1} ,\ldots,\b{y}^M \in
\R^{d_M}$ be random variables. We want to measure the dependence
between these variables. Let positive semi-definite kernels
$k^m(\cdot,\cdot): \R^{d_m}\times \R^{d_m} \to \R$
$(m=1,\ldots,M)$ be chosen in the respective spaces. Let $\F^m$
denote the RKHS associated with kernel $k^m(\cdot,\cdot)$. Having
$T$ samples $\b{y}^m_1,\ldots,\b{y}^m_T$ for all random variables
$\b{y}^m$ $(m=1,\ldots,M)$, matrices $\b{K}^m:=
[k^m(\b{y}^m_i,\b{y}^m_j)]_{i,j=1,\ldots,T}\in \R^{T \times T}$
and $\widetilde{\b{K}}^m:=\b{H}\b{K}^m\b{H}\in \R^{T \times T}$
can be created. Let the regularization parameter be chosen as
$\kappa>0$ and let $\kappa_2$ denote the auxiliary variable
$\kappa_2:=\frac{\kappa T}{2}$. It can be proven that the
computation of $\hat{J}^{\,emp}_{\text{KCCA}}$ involves the
solution of the following generalized eigenvalue problem:

\begin{align}
\lefteqn{\begin{pmatrix}
  (\widetilde{\b{K}}^1+\kappa_2\b{I}_T)^2 & \widetilde{\b{K}}^1\widetilde{\b{K}}^2 & \cdots & \widetilde{\b{K}}^1\widetilde{\b{K}}^M \\
  \widetilde{\b{K}}^2\widetilde{\b{K}}^1 &   (\widetilde{\b{K}}^2+\kappa_2\b{I}_T)^2 & \cdots & \widetilde{\b{K}}^2\widetilde{\b{K}}^M \\
  \vdots &           \vdots               &    & \vdots   \\
  \widetilde{\b{K}}^M\widetilde{\b{K}}^1 & \widetilde{\b{K}}^M\widetilde{\b{K}}^2 & \cdots & (\widetilde{\b{K}}^M+\kappa_2\b{I}_T)^2
\end{pmatrix}
\begin{pmatrix}
  \b{c}_1 \\
  \b{c}_2 \\
  \vdots \\
  \b{c}_M
\end{pmatrix}=}\label{eq:KCCA-gen-eig-problem-multi}\\
&&=\lambda
\begin{pmatrix}
  (\widetilde{\b{K}}^1+\kappa_2\b{I}_T)^2 & \b{0} & \cdots & \b{0}\\
  \b{0} & (\widetilde{\b{K}}^2+\kappa_2\b{I})^2 &  \cdots & \b{0}  \\
  \vdots &\vdots & &\vdots \\
  \b{0} & \b{0} & \cdots & (\widetilde{\b{K}}^M+\kappa_2\b{I})^2 \\
\end{pmatrix}
\begin{pmatrix}
  \b{c}_1 \\
  \b{c}_2 \\
  \vdots\\
  \b{c}_M
\end{pmatrix}. \nonumber
\end{align}
Analogously to the two-variable-case, the largest eigenvalue of this task is a
measure of the value of the pair-wise dependence of the random variables.

\paragraph{The KGV Method} \label{sec:KGV}
Equation~\eqref{eq:ISA-cost-I} in Section~\ref{sec:ISA-cost}
indicates that the ISA task can be seen as the minimization of the
mutual information. The basic idea of the KGV technique is
that---even for non-Gaussian variables---it estimates the mutual
information by the Gaussian approximation \cite{bach02kernel}.
Namely, let $\b{y}=[\b{y}^1;\ldots,\b{y}^M]$ be multidimensional
normal random variable with covariance matrix $\b{C}$. Let
$\b{C}^{i,j}\in\R^{d_m\times d_m}$ denote the cross-covariance
between components of $\b{y}^m\in\R^{d_m}$. The mutual information
between components $\b{y}^1,\ldots,\b{y}^M$ is
\cite{cover91elements}:
\[
I(\b{y}^1,\ldots,\b{y}^M)=-\frac{1}{2} \log \left(\frac{\det \b{C}}{\prod_{m=1}^{M}\det\b{C}^{m,m}}\right).
\]
The quotient $\frac{\det\b{C}}{\prod_{m=1}^M\det\b{C}^{m,m}}$ is
called the \emph{generalized variance}. If $\b{y}$ is \emph{not
normal}---this is the typical situation in the ISA task---then let
us transform the individual components $\b{y}^m$ using feature
mapping $\bm{\varphi}$ associated with the reproducing kernel and
assume that the image is a normal variable. Thus, the cost
function
\begin{equation}\label{eq:KGV}
    J_{\text{KGV}}(\b{W}):=-\frac{1}{2}\log\left[\frac{\det(\bm{\mathcal{K}})}{\prod_{m=1}^M\det(\bm{\mathcal{K}}^{m,m})}\right]
\end{equation}
is associated with the ISA task. In Equation~\eqref{eq:KGV}
$\bm{\phi}(\b{y}):=[\bm{\varphi}(\b{y}^1);\ldots;\bm{\varphi}(\b{y}^M)]$,
$\bm{\mathcal{K}}:=cov[\bm{\phi}(\b{y})]$, and the sub-matrices
are $\bm{\mathcal{K}}^{i,j}=cov[\bm{\varphi}(\b{y}^i),
\bm{\varphi}(\b{y}^j)]$. Expression
$\frac{\det(\bm{\mathcal{K}})}{\prod_{m=1}^M\det(\bm{\mathcal{K}}^{m,m})}$
is called the \emph{kernel generalized variance} (KGV).

The next theorem shows that the KGV technique can be interpreted as a decorrelation based method:

\begin{theorem}
Let ${\Sigma}\in\R^{D\times D}$ be a positive semi-definite matrix, let ${\Sigma}^{m,m}\in\R^{d_m\times d_m}$
denote the $m^{th}$ block in the diagonal of matrix ${\Sigma}$, and let $D=\sum_{m=1}^Md_m$. Then the function
\[
    0\le Q({\Sigma}):=-\frac{1}{2}\log\left[\frac{\det({\Sigma})}{\prod_{m=1}^M\det({\Sigma}^{m,m})}\right]
\]
is 0 iff ${\Sigma}=blockdiag({\Sigma}^{1,1},\ldots,{\Sigma}^{M,M})$.
\end{theorem}

This theorem can be proven for $d_m\ge 1$, as in the case of
$d_m=1$ \cite{cover91elements}, see the work of
\cite{szabo06real}. The theorem implies the following:
\begin{cor}\label{conseq:KGV-JFD:equiv}
Setting ${\Sigma}:=\bm{\mathcal{K}}$, the KGV technique is a
decorrelation technique according to feature mapping
$\bm{\varphi}$. The KGV technique aims at minimizing of
cross-covariances
\mbox{$\bm{\mathcal{K}}^{i,j}=cov[\bm{\varphi}(\b{y}^i),\bm{\varphi}(\b{y}^j)]$}
to $\b{0}$.
\end{cor}

We note that the kernel covariance (KC) ICA method
\cite{gretton05kernel}---similarly to the KCCA method---can be
extended to measure the mutual dependence of multidimensional
random variables and thus to solve the ISA task. Again, the
computation of the cost function can be converted to the solution
of a generalized eigenvalue problem. This eigenvalue problem is
provided in Appendix~\ref{sec:KC-deduction} for the sake of
completeness.

\begin{note}
The KCCA, KGV and KC methods can estimate only \emph{pair-wise} dependence.
Nonetheless, the joint mutual information can be estimated by recursive methods
computing pair-wise mutual information: for the mutual information of random
variables $\b{y}^m\in\R^{d_m}$ $(m=1,\ldots,M)$ it can be shown that the recursive
relation
\begin{equation}
I(\b{y}^1,\ldots,\b{y}^M)=\sum_{m=1}^MI\left(\b{y}^m,\left[\b{y}^{m+1},\ldots,\b{y}^M\right]\right)\label{eq:ISA-cost-I-rec}
\end{equation}
holds \cite{cover91elements}. Thus, for example, the KCCA
eigenvalue problem of \eqref{eq:KCCA-gen-eig-problem-multi} can be
replaced by pair-wise estimation of mutual information. We note
that the tree-dependent component analysis model
\cite{bach03beyond} estimates the joint mutual information from
the pair-wise mutual information.
\end{note}

\subsection{Optimization of ISA Costs}\label{sec:Opt-of-ISA-costs}
There are several possibilities to optimize ISA cost functions:
\begin{enumerate}
    \item
        Without ICA preprocessing, optimization problems concern either the \emph{Stiefel manifold}
        \cite{edelman98geometry,lippert98nonlinear,plumbley04lie,quinquis06efficient}
        or the \emph{flag manifold} \cite{nishimori06riemannian}. According to our experiences, these gradient based optimization
        methods may be stuck in poor local minima.
    \item
        According to the ISA Separation Theorem, it may be sufficient to search for
        optimal permutation of the ICA components provided by ICA preprocessing. We applied greedy permutation search:
        two coordinates of different subspaces are exchanged provided that this change decreases cost
        function $J$. Here, $J$ denotes, for example, $J_{\text{JFD}}$,
        $J_{\text{KCCA}}$, or $J_{\text{KGV}}$ depending on the ISA
        technique applied. The variable of $J$ is the permutation matrix
        $\b{P}$ using the parametrization
        $\b{W}_{\text{ISA}}=\b{P}\b{W}_{\text{ICA}}$. Pseudocode is provided in
        Table~\ref{tab:ISA-pseudocode}. Our experiences show that greedy permutation search is often sufficient for the
        estimation of the ISA subspaces. However, it is easy to generate examples in which this is not true \cite{poczos05independent2}.
        In such cases, global permutation search method of higher computational burden may become necessary \cite{szabo06cross}.
\end{enumerate}

\begin{table}
\centering
  \begin{minipage}{10.5cm}
  \centering
  \begin{tabular}{|l|}
        \hline
        \textbf{Input of the algorithm}\\
        \verb|   |ISA observation: $\{\mathbf{x}(t)\}_{t=1,\ldots,T}$\\
        \textbf{Optimization}\footnotemark \\
        \verb|   |\textbf{ICA}: on the whitened observation $\b{x}$ $\Rightarrow$ $\hat{\b{s}}_{\text{ICA}}$ estimation\\
        \verb|   |\textbf{Permutation search}\\
        \verb|       |$\b{P}:=\b{I}_D$\\
        \verb|       |repeat\\
        \verb|           |sequentially for $\forall p\in\G^{m_1},q\in\G^{m_2}(m_1\ne m_2):$\\
        \verb|               |if $J(\b{P}_{pq}\b{P})<J(\b{P})$\\
        \verb|                   |$\b{P}:=\b{P}_{pq}\b{P}$\\
        \verb|               |end\\
        \verb|       |until $J(\cdot)$ decreases in the \emph{sweep} above\\
        \textbf{Estimation}\\
        \verb|    |$\hat{\b{s}}_{\text{ISA}}=\b{P}\hat{\b{s}}_{\text{ICA}}$\\
        \hline
  \end{tabular}
  \end{minipage}
  \caption{Pseudocode of the ISA Algorithm. Cost $J$ stands for the ISA cost function of JFD, KCCA, or KGV methods.
  The permutation matrix of the ISA Separation Theorem is the variable of $J$.}
  \label{tab:ISA-pseudocode}
\end{table}
\footnotetext{Let $\G^1, \ldots, \G^M$ denote the indices of the $1^{st}, \ldots , M^{th}$ subspaces, that is, $\G^m:=\{(m-1)d+1,\ldots,md\}$, and permutation matrix $\b{P}_{pq}$ exchanges coordinates $p$ and $q$.}

\subsection{Different and Unknown Component Dimensions}\label{sec:diff-or-unknown-dims}

Here we give a quick overview how one can handle situations when
the dimensions of the subspaces are unequal, unknown, or both.
Note that the introduced uBSSD-ISA reduction, the ISA ambiguities
\cite{theis06towards} and the Separation Theorem remain the same
for subspaces of different dimensions, and thus it is sufficient
to consider the ISA problem.
\begin{enumerate}
    \item
        If the dimensions of the subspaces are different but known, the ISA task can be solved
        \begin{enumerate}
            \item
                the mask $\b{N}$ of the JFD method should be modified (see Equation~\ref{eq:JFD-cost}).
            \item
                the kernel-ISA methods include this situation; they were introduced for different subspace dimensions.
        \end{enumerate}
    \item
        If the dimension of the hidden source $\b{s}$ is known, but the individual dimensions of components
        $\b{s}^m$ are not, then clustering can exploit the dependencies between the coordinates of the estimated
        sources.
        For example:
        \begin{enumerate}
            \item
                If we assume that the hidden $\b{s}$ source has block diagonal AR dynamics of the form $\b{s}(t+1)=\b{F}\b{s}(t)+\b{e}(t)$---$\b{F}=blockdiag\left(\b{F}^1,\ldots,\b{F}^M\right)$---then connectivity of the estimated matrix
                $\b{\hat{F}}$ may
                help \cite{poczos06noncombinatorial}. One may assume that the $i$ and the $j$ coordinates are `\mbox{$\b{\hat{F}}$-connected}' if
                value \mbox{$\max\{|\hat{F}_{ij}|, |\hat{F}_{ji}|\}$} is above a certain threshold.
            \item
                Similar considerations can be applied in the ISA problem, for example, by using
                cumulant based matrices \cite{theis06towards}.
            \item
                Weaknesses of the threshold based method include (i) the uncertainty in choosing the threshold, and
                (ii) the fact that the methods are sensitive to the threshold. More robust solutions can be designed if the dependencies, for example, the mutual information
                 amongst the coordinates, are used to construct an adjacency matrix and apply a clustering method for this matrix.
                One might use, for example, hierarchical \cite{strogbauer04least} or tree-structured clustering methods \cite{bach03beyond}.
        \end{enumerate}
\end{enumerate}

\section{Illustrations}\label{sec:illustrations}
The efficiency of the algorithms of Section~\ref{sec:ISA-methods}
are illustrated. Test cases are introduced in
Section~\ref{sec:databases}. The quality of the solutions will be
measured by the normalized Amari-error, the Amari-index
(Section~\ref{sec:amaridist}). Numerical results are presented in
Section~\ref{sec:simulations}.

\subsection{Databases}\label{sec:databases}
We define five databases to study our identification algorithms.
We do not know whether they satisfy \eqref{eq:suff} or not.
According to our experiences, the ISA Separation Theorem works on
these examples.

\subsubsection{The 3D-geom, the Celebrities and the ABC Database}

The first 3 databases are illustrated in
Figure~\ref{fig:database:3D-geom,celebs,ABC}. In the
\emph{3D-geom} test $\b{s}^m$s were random variables uniformly
distributed on 3-dimensional geometric forms ($d=3$). We chose 6
different components ($M=6$) and, as a result, the dimension of
the hidden source $\b{s}$ is $D_s=18$. The \emph{celebrities} test
has 2-dimensional source components generated from cartoons of
celebrities \mbox{($d=2$)}.\footnote{See http://www.smileyworld.com.}
Sources $\b{s}^m$ were generated by sampling 2-dimensional
coordinates proportional to the corresponding pixel intensities.
In other words, 2-dimensional images of celebrities were
considered as density functions. $M=10$ was chosen. In the
\emph{ABC} database, hidden sources $\b{s}^m$ were uniform
distributions defined by 2-dimensional images ($d=2$) of the
English alphabet. The number of components varied as
$M=2,5,10,15$, and thus the dimension of the source $D_s$ was
$4,10,20,30$, respectively.

\captionsetup[subfloat]{labelformat=empty} 
\begin{figure}%
\centering%
\hfill\subfloat[][(a)]{\includegraphics[width=1.5cm]{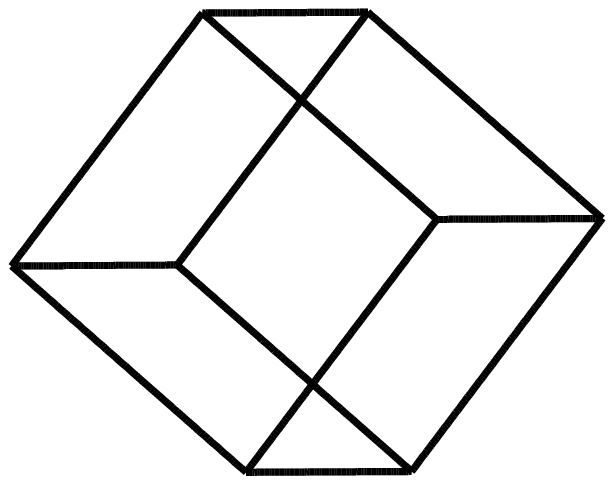}\hfill%
\includegraphics[width=1.5cm]{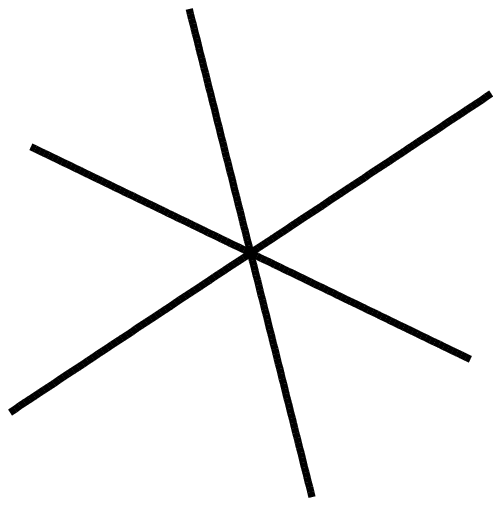}\hfill%
\includegraphics[width=1.5cm]{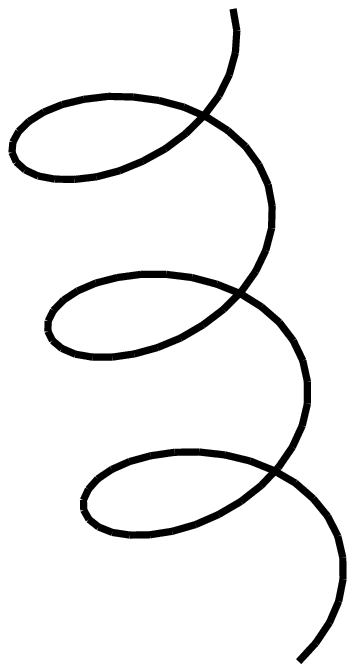}\hfill%
\includegraphics[width=1.5cm]{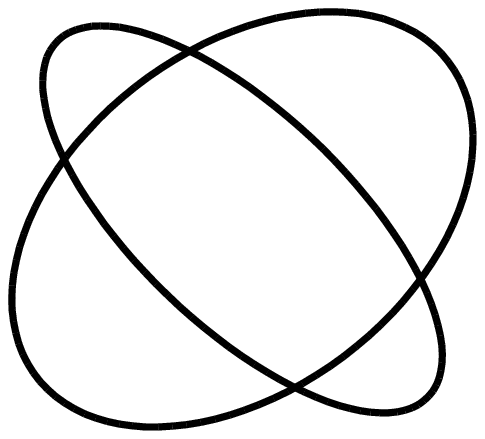}\hfill%
\includegraphics[width=1.5cm]{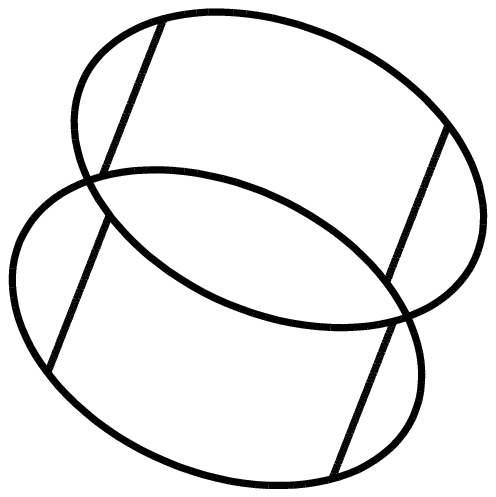}\hfill%
\includegraphics[width=1.5cm]{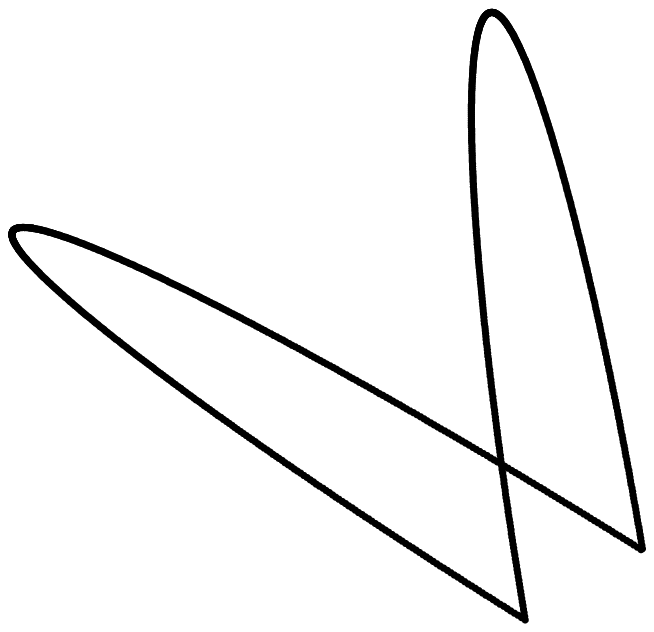}}
\hfill
\subfloat[][(c)]{\includegraphics[width=5.2cm]{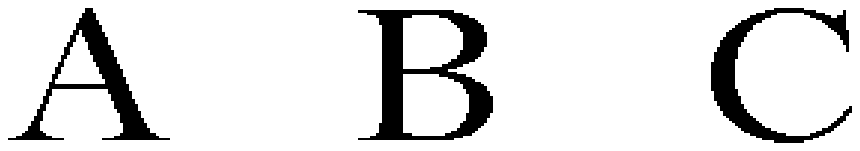}}\hfill\hfill\hfill\\
\subfloat[][(b)]{\includegraphics[width=1.4cm]{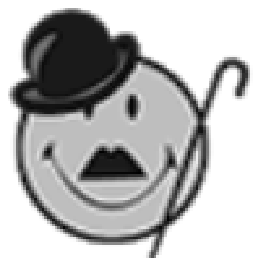}
\includegraphics[width=1.4cm]{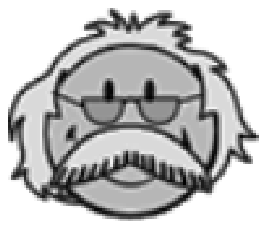}
\includegraphics[width=1.4cm]{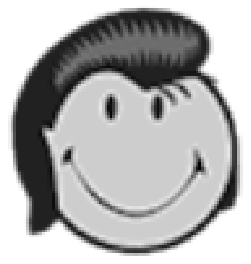}
\includegraphics[width=1.4cm]{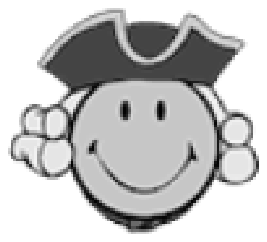}
\includegraphics[width=1.4cm]{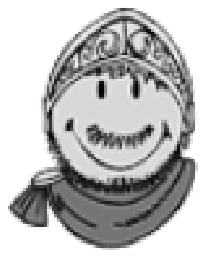}
\includegraphics[width=1.4cm]{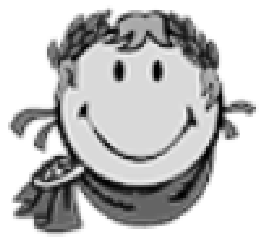}
\includegraphics[width=1.4cm]{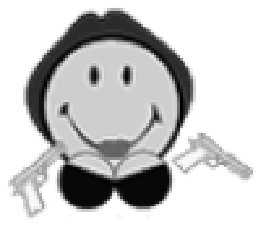}
\includegraphics[width=1.4cm]{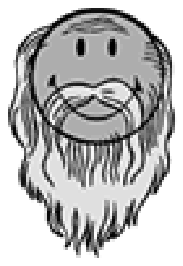}
\includegraphics[width=1.4cm]{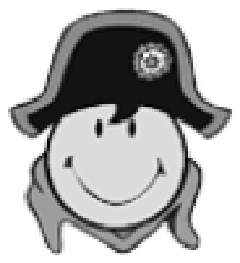}
\includegraphics[width=1.4cm]{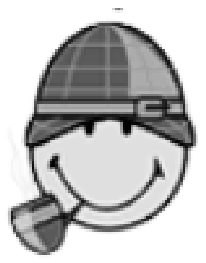}}
\caption[]{Illustration of the  \emph{3D-geom}, \emph{celebrities}
and  \emph{ABC} databases. (a): database \emph{3D-geom}, 6
3-dimensional components ($M=6$, $d=3$). Hidden sources are
uniformly distributed variables on 3-dimensional geometric
objects. (b): database $celebrities$. Density functions of the
hidden sources are proportional to the pixel intensities of the
2-dimensional images ($d=2$). Number of hidden components: $M=10$.
(c): database $ABC$. Here, the hidden sources $\b{s}^m$ are
uniformly
distributed on images ($d=2$) of letters. Number of components $M$ varies between 2 (A and B) and 15 (A-O).}%
\label{fig:database:3D-geom,celebs,ABC}
\end{figure}\captionsetup[subfloat]{labelformat=parens}

\subsubsection{The all-k-independent Database}
The $d$-dimensional hidden components $\b{u}:=\b{s}^m$ were
created as follows: coordinates $u_i(t)$ $(i=1,\ldots, k)$ were
uniform random variables on the set \{0,\ldots,k-1\}, whereas
$u_{k+1}$ was set to $mod(u_1 + \ldots + u_k,k)$. In this
construction, every $k$-element subset of $\{u_1,\ldots,u_{k+1}\}$
is made of independent variables. This database is called the
\emph{all-k-independent} problem
\cite{poczos05independent2,szabo06cross}. In our simulations
$d=k+1$ was set to $3$ or $4$ and we used $2$ components ($M=2$).
Thus, source dimension $D_s$ was either $6$ or $8$.

\subsubsection{The Beatles Database}

Our \emph{Beatles} test is a non-i.i.d.\ example. Here, hidden sources are stereo Beatles
songs.\footnote{See http://rock.mididb.com/beatles/.} $8$ kHz sampled portions of two songs (A Hard Day's Night, Can't Buy
Me Love) made the hidden $\b{s}^m$s. Thus, the dimension of the components $d$ was $2$, the number of the components
$M$ was $2$, and the dimension of the problem $D_s$ was $4$.

\subsection{Normalized Amari-error, the Amari-index}\label{sec:amaridist}

We have shown in Section~\ref{sec:uBSSD2(u)ISA} that the uBSSD
task can be reduced to an ISA task. Consequently, we use ISA
performance measure to evaluate our algorithms. Assume that there
are $M$ pieces of $d$-dimensional hidden components in the ISA
task, $\b{A}$ is the mixing matrix, and $\b{W}$ is the estimated
demixing matrix. Then optimal estimation provides matrix
$\b{G}:=\b{W}\b{A}$, a block-permutation matrix made of $d\times
d$ sized blocks. Let matrix $\b{G}\in\R^{D\times D}$ be decomposed
into $d\times d$ blocks:
$\b{G}=\left[\b{G}^{ij}\right]_{i,j=1,\ldots,M}$. Let $g^{i,j}$
denote the sum of the absolute values of the elements of matrix
$\b{G}^{i,j}\in\R^{d\times d}$. We used the normalized version
\cite{szabo06cross} of the Amari-error \cite{amari96new} adapted
to the ISA task \cite{theis05blind,theis05multidimensional}
defined as:
\[
    r(\b{G}):=\frac{1}{2M(M-1)}\left[\sum_{i=1}^M\left(\frac{
\sum_{j=1}^Mg^{ij}}{\max_jg^{ij}}-1\right)+ \sum_{j=1}^M\left(\frac{ \sum_{i=1}^Mg^{ij}}{\max_ig^{ij}}-1\right)\right].
\]
We refer to the normalized Amari-error as the Amari-index. One can
see that $0\le r(\b{G})\le 1$ for any matrix $\b{G}$, and
$r(\b{G})=0$ if and only if $\b{G}$ is a block-permutation matrix
with $d\times d$ sized blocks. Normalization is advantageous; we
can compare the precision of ISA procedures and procedures, which
are reduced to ISA tasks.

\subsection{Simulations}\label{sec:simulations}

Results on databases \emph{3D-geom}, \emph{celebrities}, \emph{ABC}, \emph{all-$k$-independent} and \emph{Beatles} are
provided here. These experimental studies have two main parts:
\begin{enumerate}
    \item
        The efficiency of the JFD method on the uBSSD task is demonstrated in Section~\ref{sec:JFD-on-uBSSD}.
    \item
        The derived KCCA, KGV kernel-ISA methods were tested on ISA tasks. We show examples in which these methods are favorable
        over the JFD method in Section~\ref{sec:kernelISA-vs-JFD}.
\end{enumerate}
In both cases the tasks are either ISA tasks or can be reduced to ISA
(Section~\ref{sec:uBSSD2(u)ISA}). Thus, we used the Amari-index
(Section~\ref{sec:amaridist}) to measure and compare the performance of the
different methods. For each individual parameter, $50$ random runs were averaged.
Our parameters are: $T$, the sample number  of observations $\b{x}(t)$, $L$, the
parameter of the length of the convolution (the length of the convolution is $L+1$),
$M$, the number of the components, and $d$, the dimension of the components,
depending on the test. Random run means random choice of quantities $\b{H}[z]$ and
$\b{s}$.

\subsubsection{JFD on uBSSD}\label{sec:JFD-on-uBSSD}
Our results concerning the uBSSD task are delineated. As we showed
in Section~\ref{sec:uBSSD2(u)ISA}, the temporal concatenation can
turn the uBSSD task into an ISA problem. These ISA tasks
associated with simple uBSSD problems can easily become more than
\mbox{100-dimensional}. Earlier ISA methods cannot deal with such
`high dimensional' problems. This is why we resorted to the recent
JFD method (Section~\ref{sec:JFD}), which seemed to be efficient
in solving such large problems under the following circumstances:
Equation~\eqref{eq:D-ISA} implies that the dimension
$D_{\text{ISA}}$ of the derived ISA task with fixed $L$ and $D_s$
decreases provided that the difference \mbox{$D_x-D_s\ge 1$}
increases. This coincides with our experiences: the higher this
difference is, the smaller number of samples can reach the same
precision. Below, studies for $D_x-D_s=D_s$ ($D_x=2D_s$) are
presented. This choice was amenable to the JFD method.\footnote{We
note that the hardest $D_x-D_s=1$ task is also feasible. However,
the sample number necessary to find the solution grows
considerably, as can be expected from
\eqref{eq:L'-ineq-explicit}.} In this case the dimension of the
ISA task in \eqref{eq:D-ISA} simplifies to the form
\[
    D_{\text{ISA}}=2D_sL.
\]

The JFD technique works with the pseudocode given in
Table~\ref{tab:ISA-pseudocode}: it reduces the \mbox{uBSSD} task
to the ISA task, where the fastICA algorithm
\cite{hyvarinen97fast} was chosen to perform the ICA computation.
In the JFD cost, we chose manifold $\F$ as
\mbox{$\F:=\{\b{u}\rightarrow \cos(\b{u}), \b{u}\rightarrow
\cos(2\b{u})\}$}, where the functions operated on the coordinates
separately \cite{szabo06real}. For the `observations', the
elements of mixing matrices $\b{H}_l$ in Equation~\eqref{eq:obs}
were generated independently from standard normal
distributions.\footnote{Uniform distribution on $[0.1]$, instead
of normal distribution, showed similar performance.}
\begin{figure}%
\centering%
\subfloat[][]{\includegraphics[width=7.9cm]{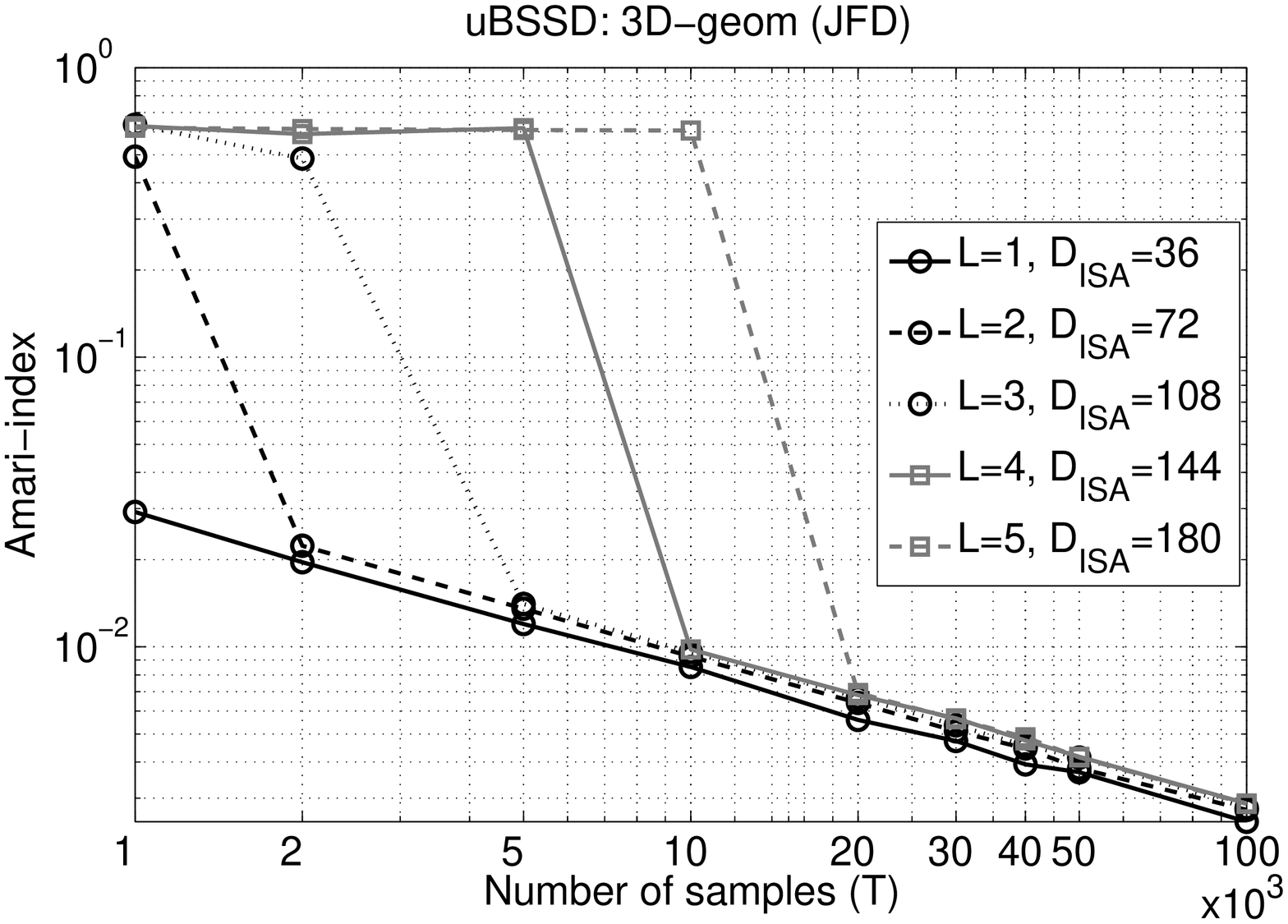}}%
\subfloat[][]{\includegraphics[width=7.9cm]{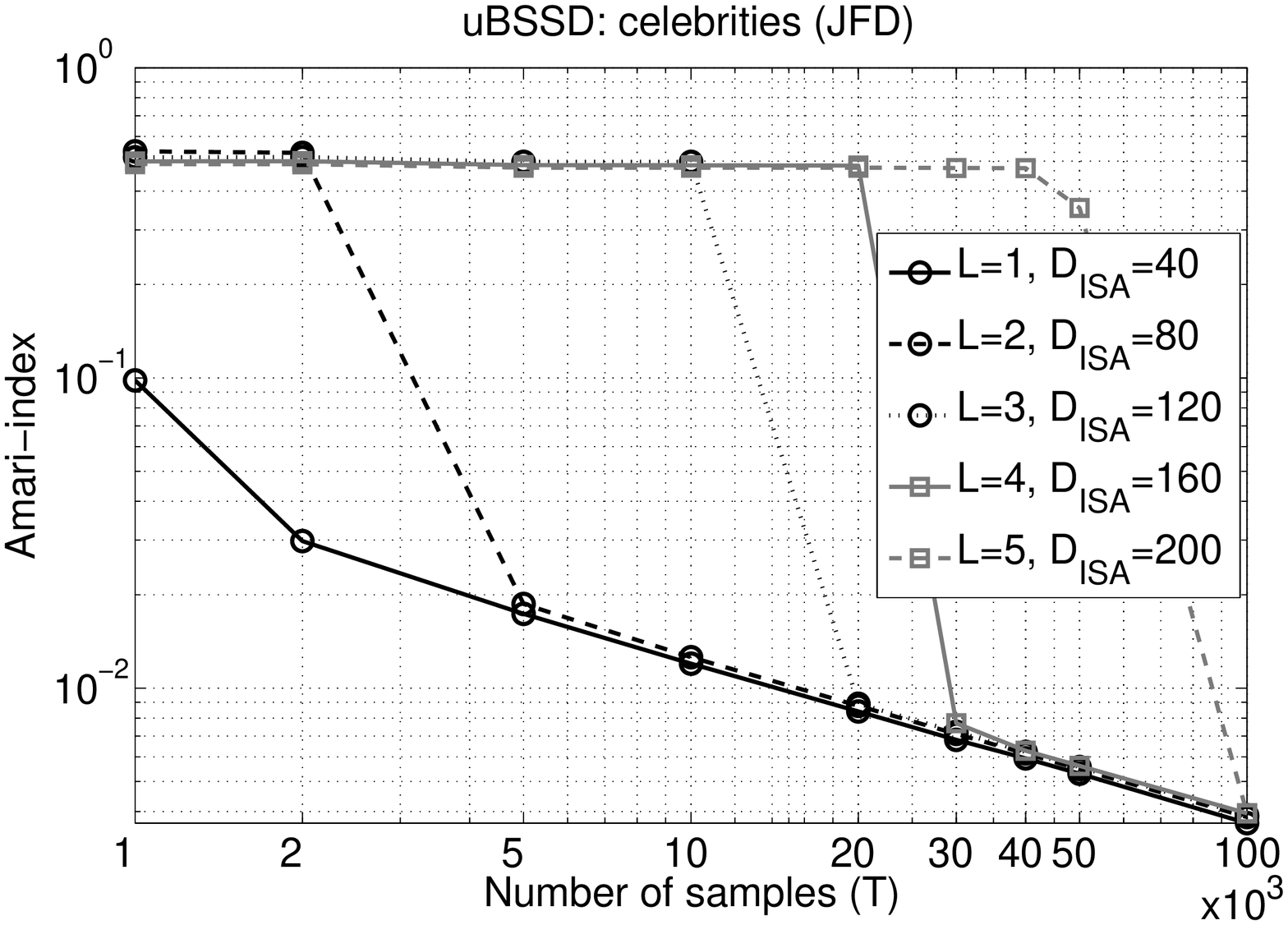}}%
\caption[]{Estimation error of the JFD method on the \emph{3D-geom} and
\emph{celebrities} databases: Amari-index as a function of sample number on log-log
scale for different convolution lengths. (a): \emph{3D-geom}, (b):
\emph{celebrities} database. Dimension of the ISA task: $D_{\text{ISA}}$.
For further information, see Table~\ref{tab:JFD-amari-dists-3D-geom-celebrities}.}%
\label{fig:JFD-amari-3D-geom-celebs}%
\end{figure}

\begin{table}
    \centering
    \begin{tabular}{|@{\hspace{3pt}}c@{\hspace{3pt}}|@{\hspace{3pt}}c@{\hspace{3pt}}|@{\hspace{3pt}}c@{\hspace{3pt}}|@{\hspace{3pt}}c@{\hspace{3pt}}|@{\hspace{3pt}}c@{\hspace{3pt}}|@{\hspace{3pt}}c@{\hspace{3pt}}|}
    \hline
        &$L=1$ & $L=2$ & $L=3$ & $L=4$ & $L=5$\\
    \hline\hline
        3D-geom&$0.25\%$ $(\pm 0.01)$ & $0.27\%$ $(\pm 0.03)$ & $0.28\%$ $(\pm 0.02)$ & $0.29\%$ $(\pm 0.03)$ & $0.29\%$ $(\pm 0.01)$\\
    \hline
        celebrities&$0.37\%$ $(\pm 0.01)$ & $0.38\%$ $(\pm 0.01)$ & $0.39\%$ $(\pm 0.01)$ & $0.39\%$ $(\pm 0.01)$ & $0.40\%$ $(\pm 0.01)$\\
    \hline
    \end{tabular}
    \caption{The Amari-index of the JFD method for database \emph{3D-geom} and $celebrities$, for different convolution lengths: average $\pm$ deviation.
    Number of samples: $T=100,000$. For other sample numbers between $1,000\le T < 100,000$ see Figure~\ref{fig:JFD-amari-3D-geom-celebs}.}
    \label{tab:JFD-amari-dists-3D-geom-celebrities}
\end{table}

\begin{figure}%
\centering%
\subfloat[][]{\includegraphics[width=7.9cm]{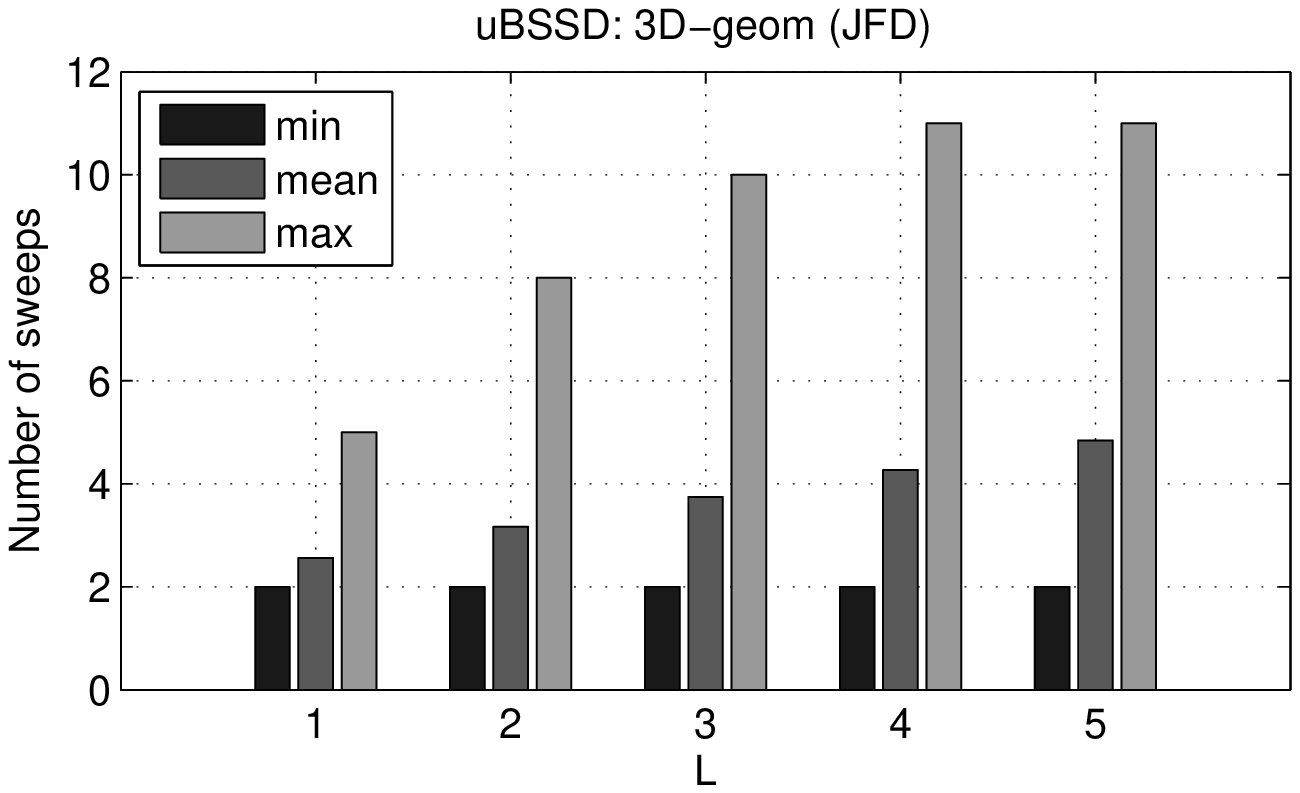}}%
\subfloat[][]{\includegraphics[width=7.9cm]{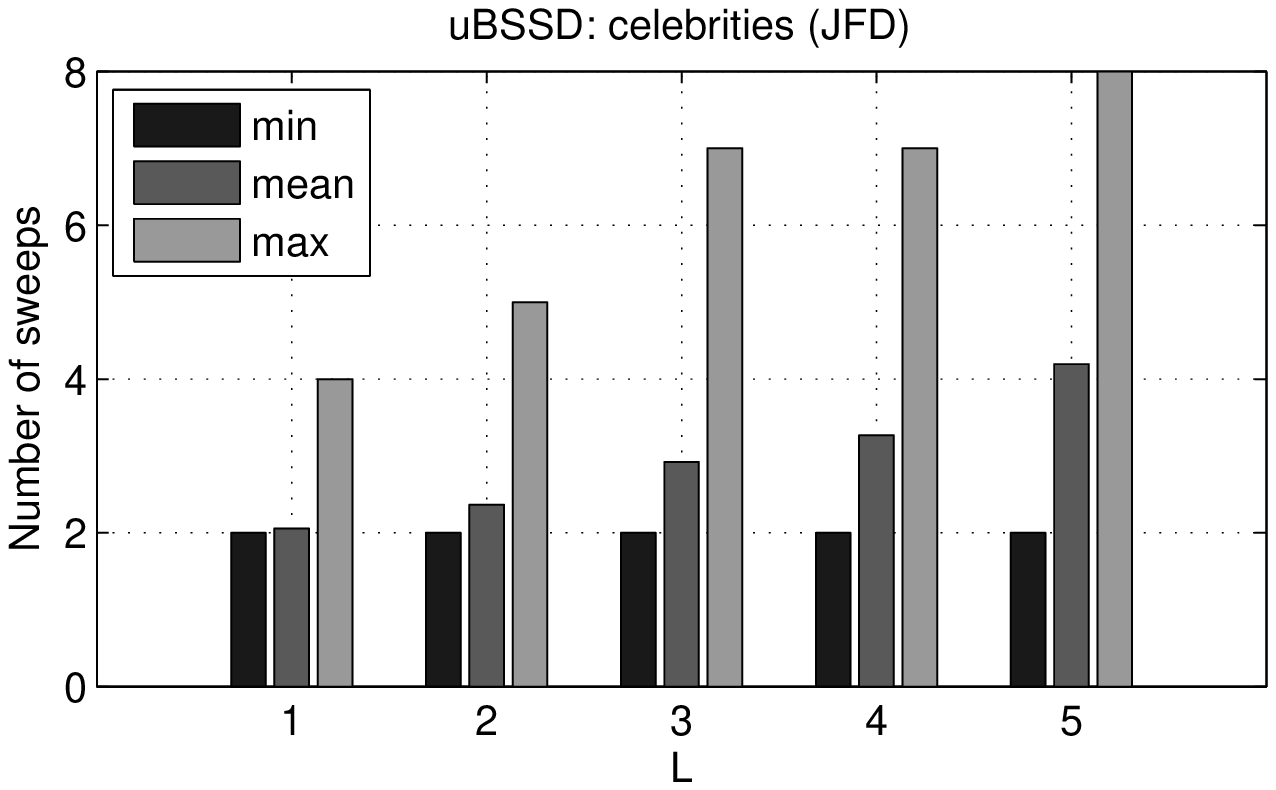}}%
\caption[]{Number of sweeps in permutation search needed for the JFD method as a function of the convolution
length. (a): \emph{3D-geom}, (b): celebrities database. Black: minimum, gray: average, light gray: maximum.}%
\label{fig:JFD-sweep-3D-geom-celebs}%
\end{figure}
\begin{figure}%
\centering%
\hfill\subfloat[][]{\includegraphics[width=3.5cm]{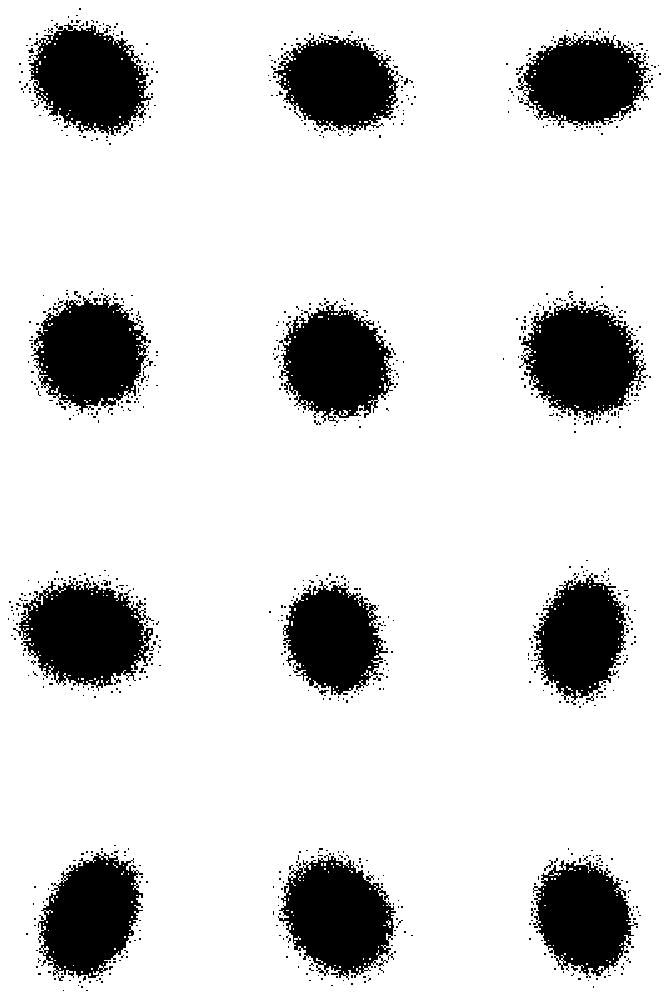}}\hfill%
\subfloat[][]{\includegraphics[width=4.85cm]{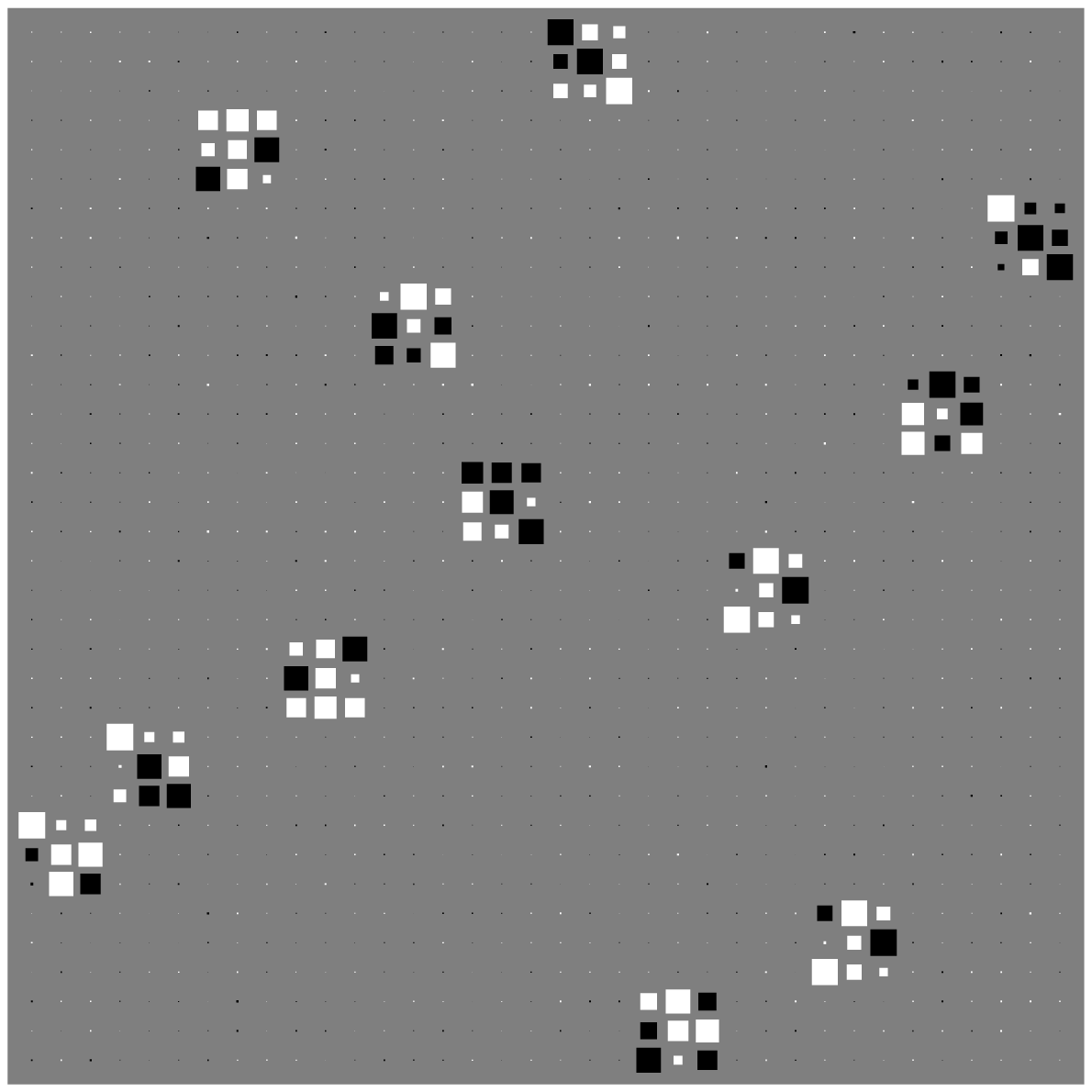}}\hfill%
\subfloat[][]{\includegraphics[width=3.4cm]{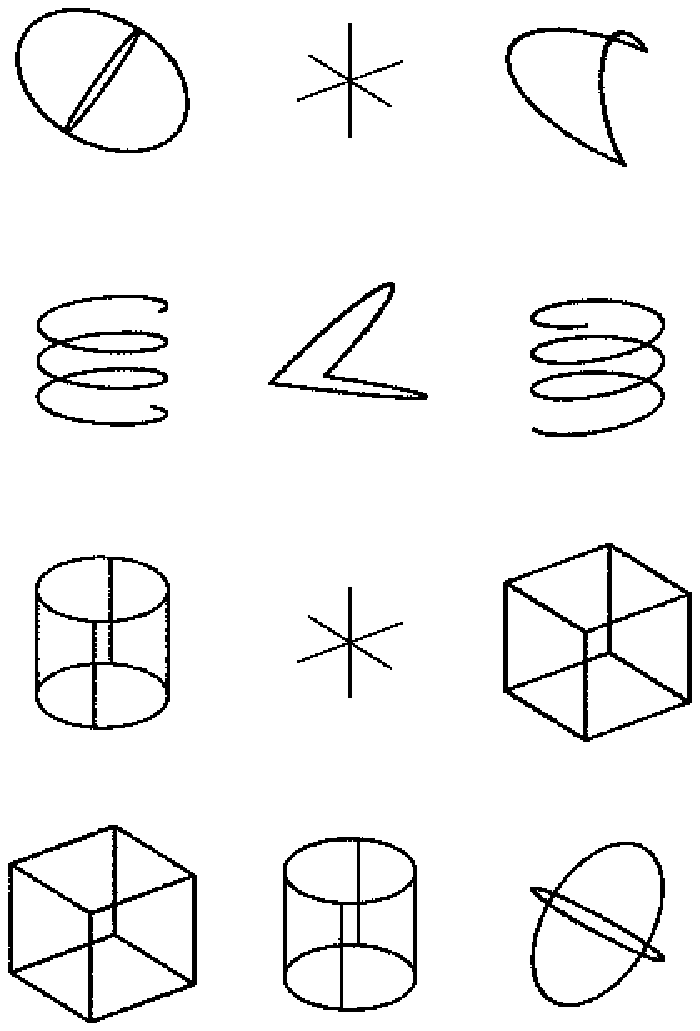}}\hfill\hfill%
\caption[]{Illustration of the JFD method on the uBSSD task for
the \emph{3D-geom} database. Sample number $T=100,000$,
convolution length $L=1$, Amari-index: $0.25\%$. (a): observed
convolved signals $\b{x}(t)$. (b) Hinton-diagram: the product of
the mixing matrix of the derived ISA task and the estimated
demixing matrix (= approximately block-permutation matrix with
$3\times 3$ blocks). (c): estimated components.
Note: hidden components are recovered $L+L'=2$ times, up to permutation and orthogonal transformation.}%
\label{fig:JFD-3D-geom-demo}%
\end{figure}
\begin{figure}%
\centering%
\hfill\subfloat[][]{\includegraphics[width=4.3cm]{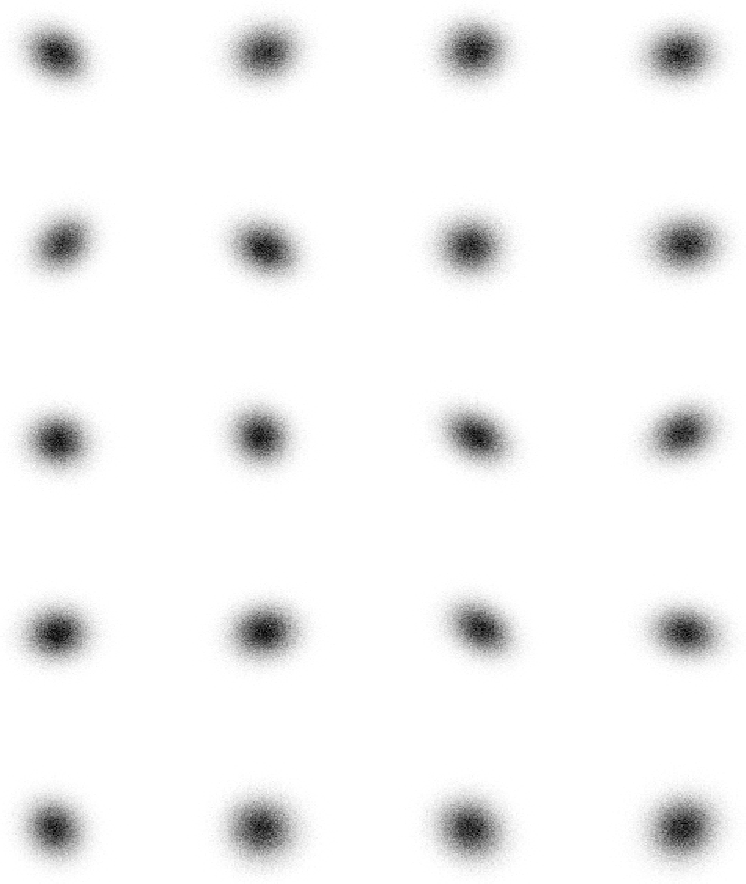}}\hfill%
\subfloat[][]{\includegraphics[width=4.65cm]{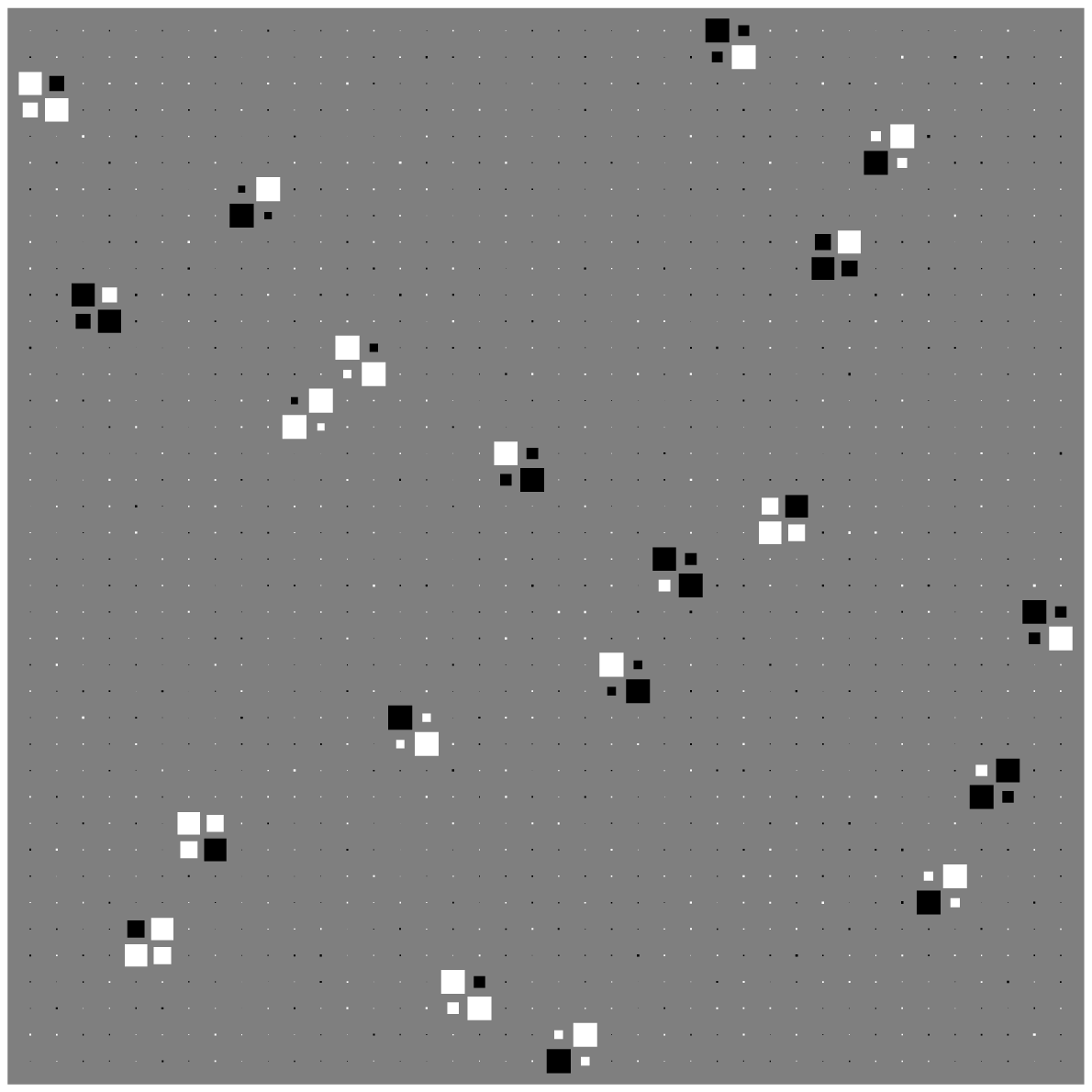}}\hfill%
\subfloat[][]{\includegraphics[width=4.3cm]{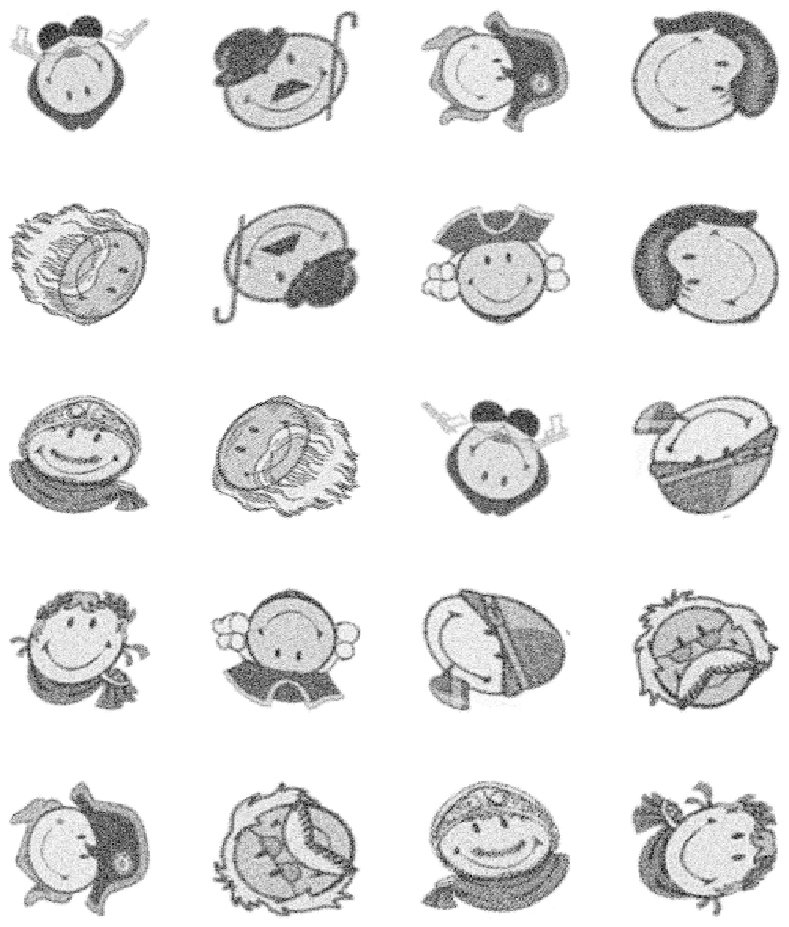}}\hfill\hfill%
\caption[]{Illustration of the JFD method on the uBSSD task for
the \emph{celebrities} database. Sample number $T=100,000$,
convolution length $L=1$, Amari-index: $0.37\%$. (a): observed
convolved signals $\b{x}(t)$. (b) Hinton-diagram: the product of
the mixing matrix of the derived ISA task and the estimated
demixing matrix (= approximately block-permutation matrix with
$2\times 2$ blocks).  (c): estimated components.
Note: hidden components are recovered $L+L'=2$ times, up to permutation and orthogonal transformation.}%
\label{fig:JFD-celebs-demo}%
\end{figure}

We studied the dependence of the precision versus the sample number on databases \mbox{\emph{3D-geom}} and
\emph{celebrities}. The dimension and the number of the components were $d=3$ and $M=6$ for the \emph{3D-geom} database
and $d=2$ and $M=10$ for the \emph{celebrities} database, respectively. In both cases the sample number $T$ varied
between $1,000$ and $100,000$. The parameter of the length of the convolution took $L=1,\ldots,5$ values. Thus, the
length of the convolution changed between $2$ and $6$. Our results are summarized in
Figure~\ref{fig:JFD-amari-3D-geom-celebs}. The values of the errors are given in
Table~\ref{tab:JFD-amari-dists-3D-geom-celebrities}. The number of sweeps needed to optimize the permutations after
performing ICA is provided in Figure~\ref{fig:JFD-sweep-3D-geom-celebs}. Figures~\ref{fig:JFD-3D-geom-demo} and
\ref{fig:JFD-celebs-demo} illustrate the estimations of the JFD technique on the \emph{3D-geom} and the
\emph{celebrities} databases, respectively.

Figure~\ref{fig:JFD-amari-3D-geom-celebs} demonstrates that the
JFD algorithm was able to uncover the hidden components with high
precision. The precision of the estimations shows similar
characteristics on the \mbox{\emph{3D-geom}} and the
\emph{celebrities} databases. The Amari-index is approximately
constant for small sample numbers. For each curve, above a certain
threshold the Amari-index decreases suddenly and after the sudden
decrease the precision follows a power law $r(T)\propto T^{-c}$
$(c>0)$. The power law decline is manifested by straight line on
log-log scale. The slopes of these straight lines are very close
to one another. The number of sweeps was between $2$ and $11$ ($2$
and $8$) for the \emph{3D-geom} (\emph{celebrities}) tests over
all sample numbers, for $1\le L\le 5$ and for $50$ random
initializations. According to
Table~\ref{tab:JFD-amari-dists-3D-geom-celebrities}, the
Amari-index for sample number $T=100,000$ is below $1\%$
($0.25-0.40\%$) with small standard deviations ($0.01-0.03$).

\begin{figure}%
\centering%
\subfloat[][]{\includegraphics[width=7.9cm]{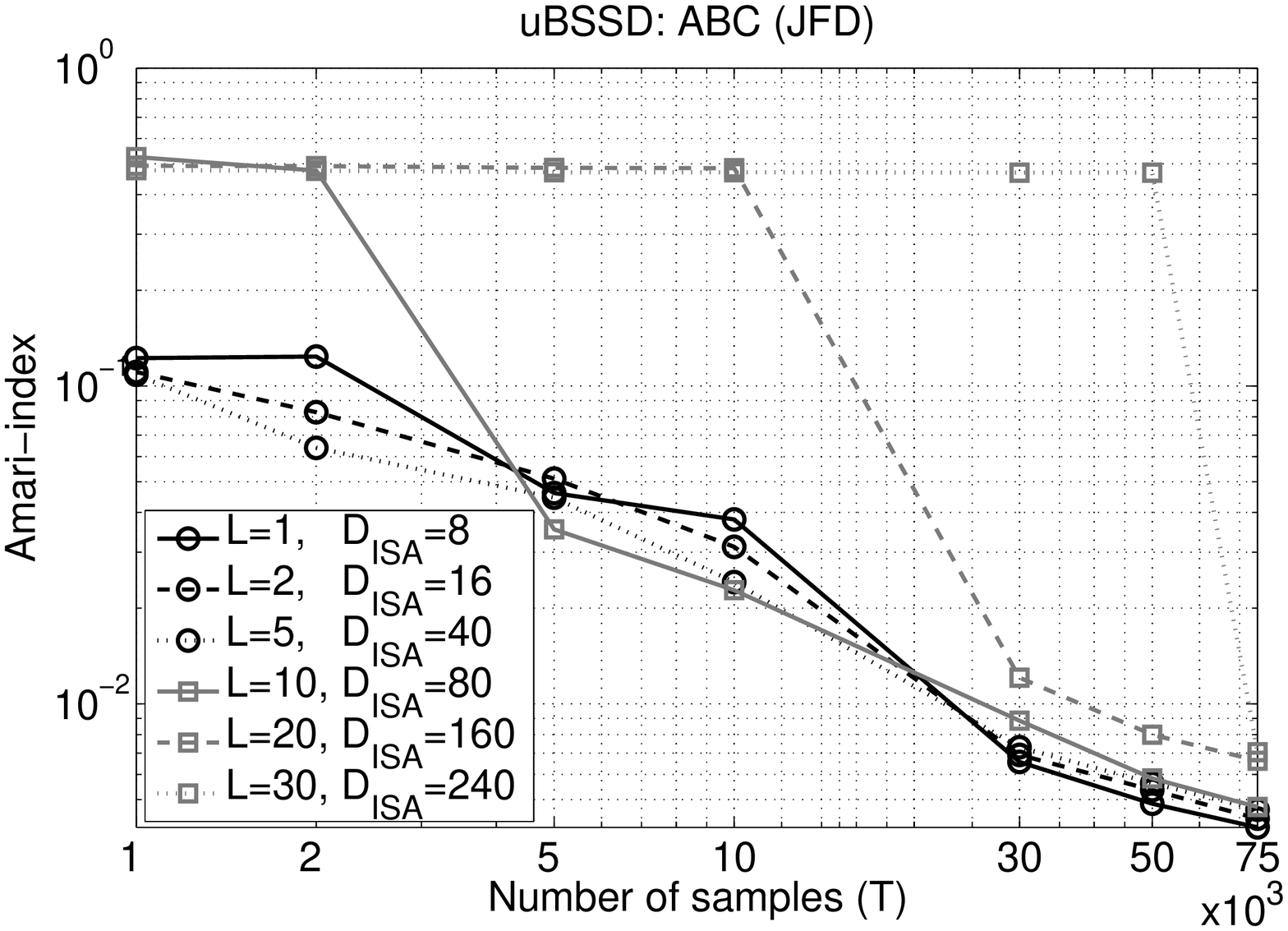}}%
\subfloat[][]{\includegraphics[width=7.9cm]{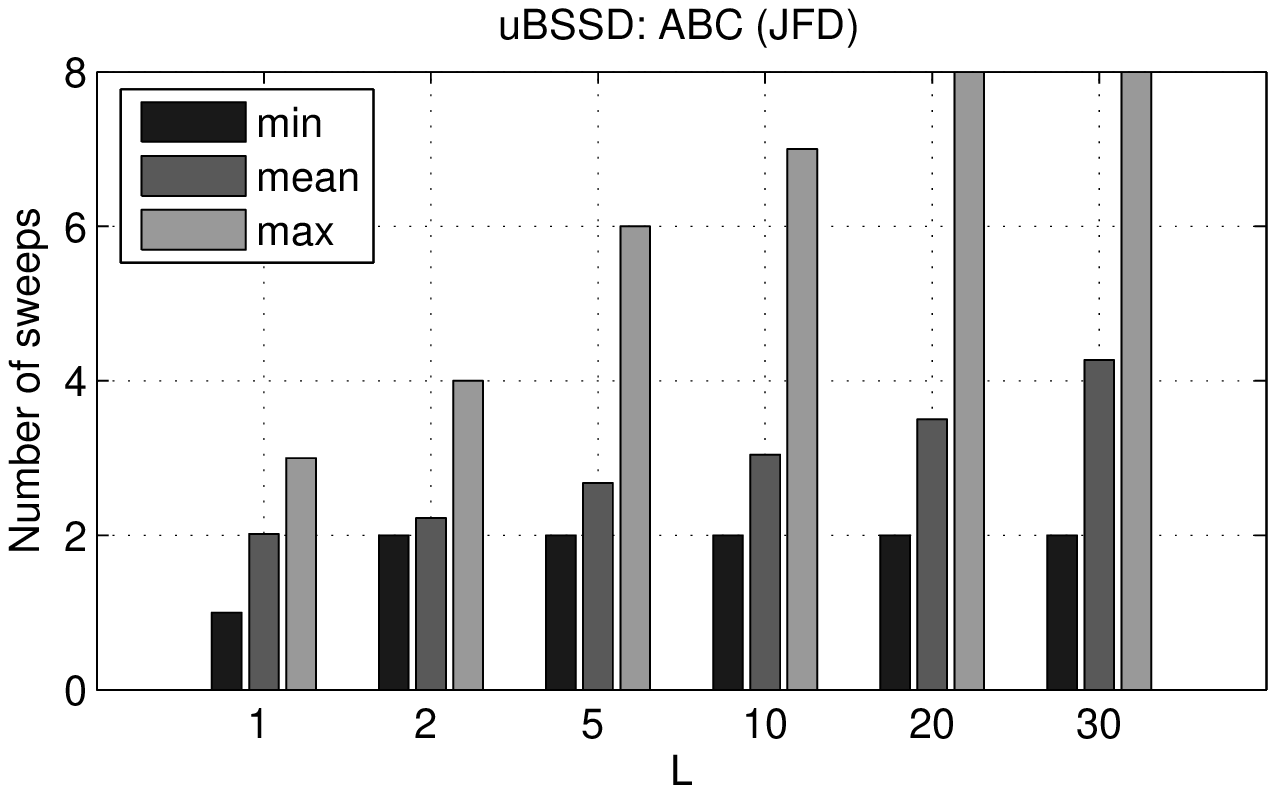}}%
\caption[]{(a): Amari-index of the JFD method on the \emph{ABC} database as a function of sample number and for
different convolution lengths on log-log scale. (b): Number of sweeps of permutation optimization on the derived ISA
task as a function of convolution length. Dimension of the ISA task: $D_{\text{ISA}}$. Black: minimum, gray: average,
light gray: maximum. For further information, see Table~\ref{tab:JFD-amari-dists-ABC}.}%
\label{fig:JFD-amari-sweep-ABC}%
\end{figure}

In another test the \emph{ABC} database was used. The number and the dimensions of the components were minimal ($d=2$,
$M=2$) and the dependence on the convolution length was tested. Parameter $L$ took values on $1,2,5,10,20,30$. The
number of observations varied between \mbox{$1,000\le T\le 75,000$}. The Amari-index and the sweep number of the
optimization are illustrated in Figure~\ref{fig:JFD-amari-sweep-ABC}. Precise values of the Amari-index are provided in
Table~\ref{tab:JFD-amari-dists-ABC}.

\begin{table}
    \centering
    \begin{tabular}{|@{\hspace{2pt}}c@{\hspace{1pt}}|@{\hspace{2pt}}c@{\hspace{1pt}}|@{\hspace{2pt}}c@{\hspace{1pt}}|@{\hspace{2pt}}c@{\hspace{1pt}}|@{\hspace{2pt}}c@{\hspace{1pt}}|@{\hspace{2pt}}c@{\hspace{1pt}}|}
    \hline
        $L=1$ & $L=2$ & $L=5$ & $L=10$ & $L=20$& $L=30$\\
    \hline\hline
        $0.41\%$ $(\pm 0.06)$ & $0.44\%$ $(\pm 0.05)$ & $0.46\%$ $(\pm 0.05)$ & $0.47\%$ $(\pm 0.03)$ & $0.66\%$ $(\pm 0.13)$& $0.70\%$ $(\pm 0.11)$\\
    \hline
    \end{tabular}
    \caption{Amari-index of the JFD method for \emph{ABC} database for different convolution lengths: average $\pm$ deviation. Number of samples:
    $T=75,000$.
    For other sample numbers between $1,000\le T <  75,000$ see Figure~\ref{fig:JFD-amari-sweep-ABC}(a).}
    \label{tab:JFD-amari-dists-ABC}
\end{table}

According to Figure~\ref{fig:JFD-amari-sweep-ABC}, the JFD method
found the hidden components. The `power law' decline of the
Amari-index, which was apparent in the \emph{3D-geom} and the
\emph{celebrities} databases, appears for the \emph{ABC} test,
too. The figure indicates that for $75,000$ samples and for $L=30$
(convolution length is $31$) the problem is still amenable to the
JFD technique. The number of sweeps required for the optimization
of the permutations was between $1$ and $8$ for all sample numbers
\mbox{$1,000\le T \le 75,000$}, parameters $1\le L\le 30$ and for
all $50$ random initializations. According to
Table~\ref{tab:JFD-amari-dists-ABC}, for sample number $T=75,000$
the Amari-index stays below $1\%$ on average ($0.41-0.7\%$) and
has a small ($0.03-0.11$) standard deviation.

\begin{figure}%
\centering%
\subfloat[][]{\includegraphics[width=7.9cm]{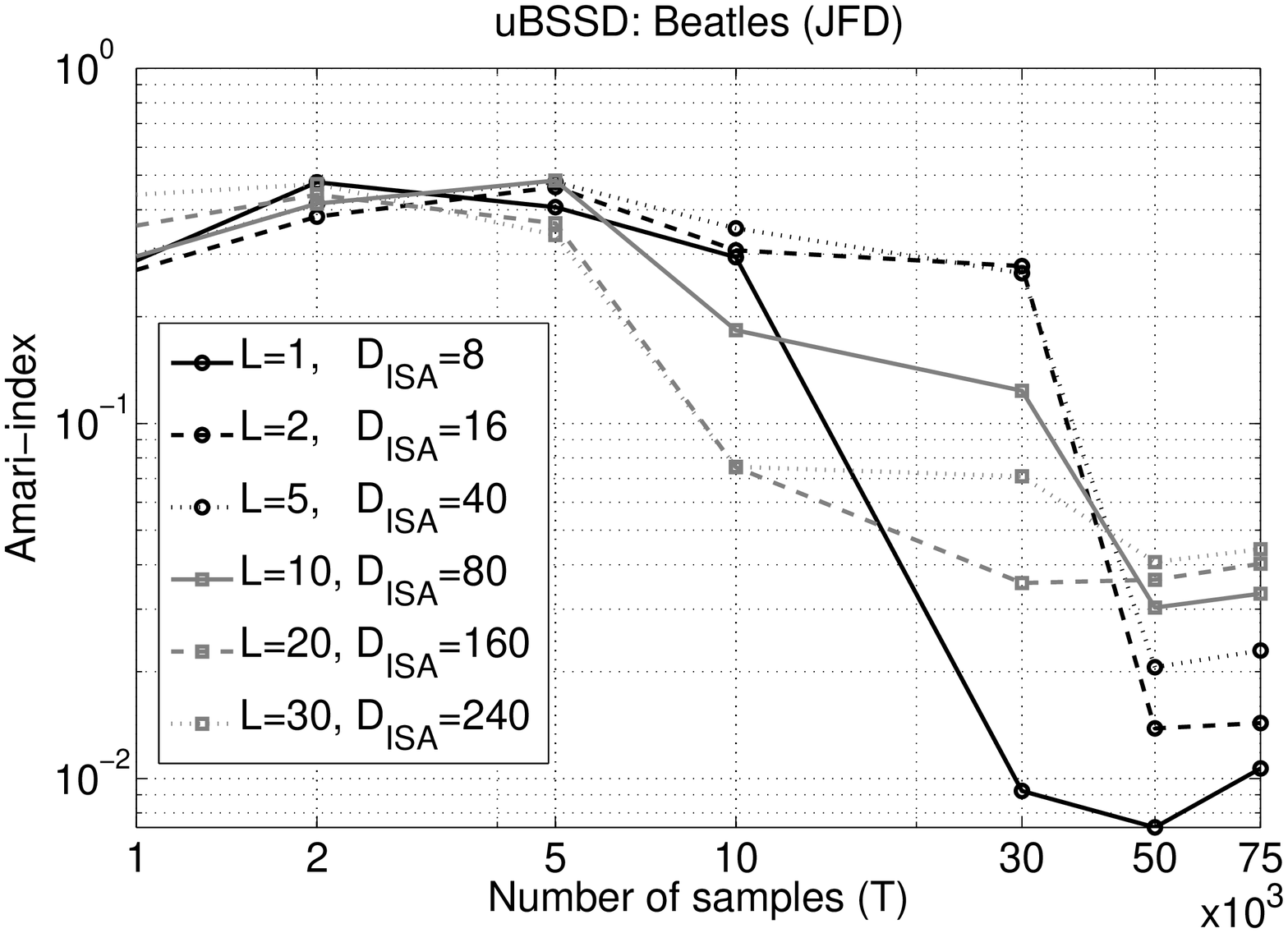}}%
\subfloat[][]{\includegraphics[width=7.9cm]{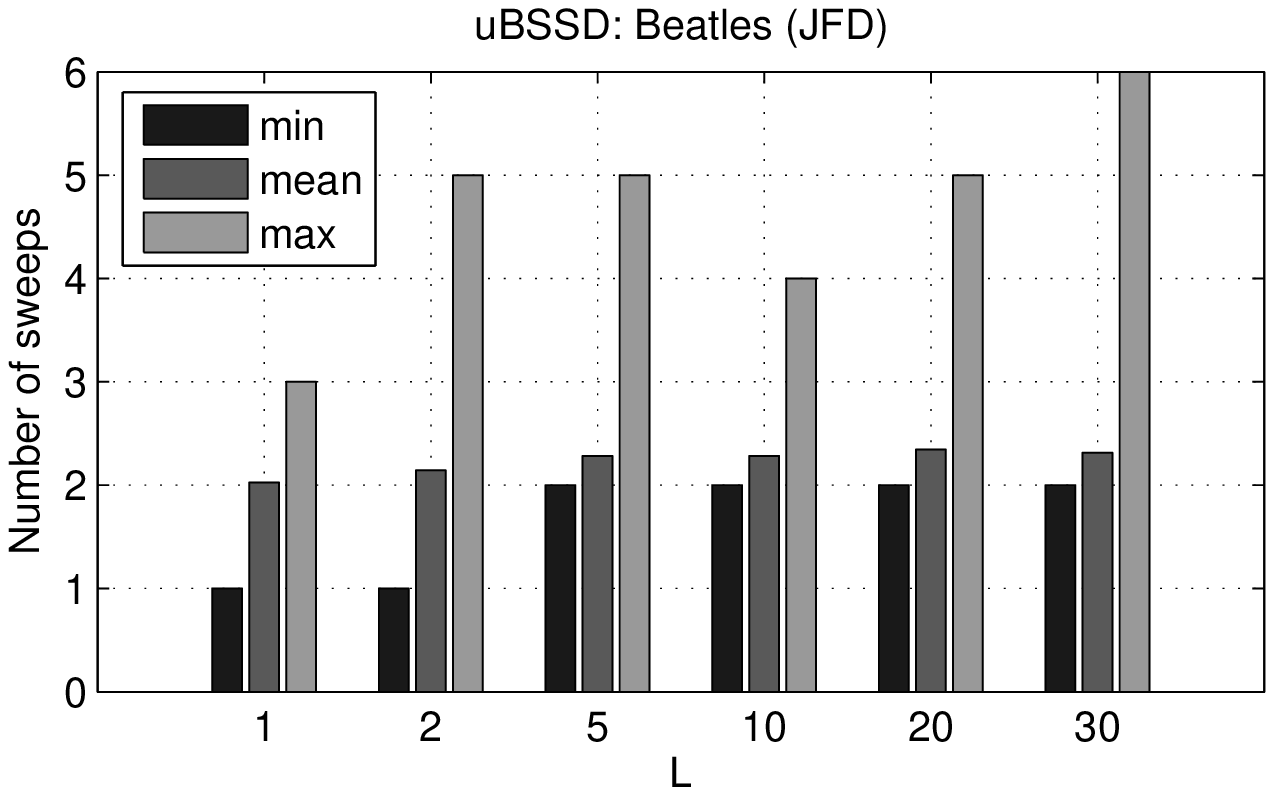}}%
\caption[]{(a): Amari-index of the JFD method on the
\emph{Beatles} database as a function of sample number and for
different convolution lengths on log-log scale. (b): Number of
sweeps of permutation optimization on the derived ISA task as a
function of convolution length. Dimension of the ISA task:
$D_{\text{ISA}}$. Black: minimum, gray: average,
light gray: maximum. For further information, see Table~\ref{tab:JFD-amari-dists-Beatles}.}%
\label{fig:JFD-amari-sweep-Beatles}%
\end{figure}

\begin{figure}%
\centering%
\subfloat[][]{\includegraphics[width=3.2cm]{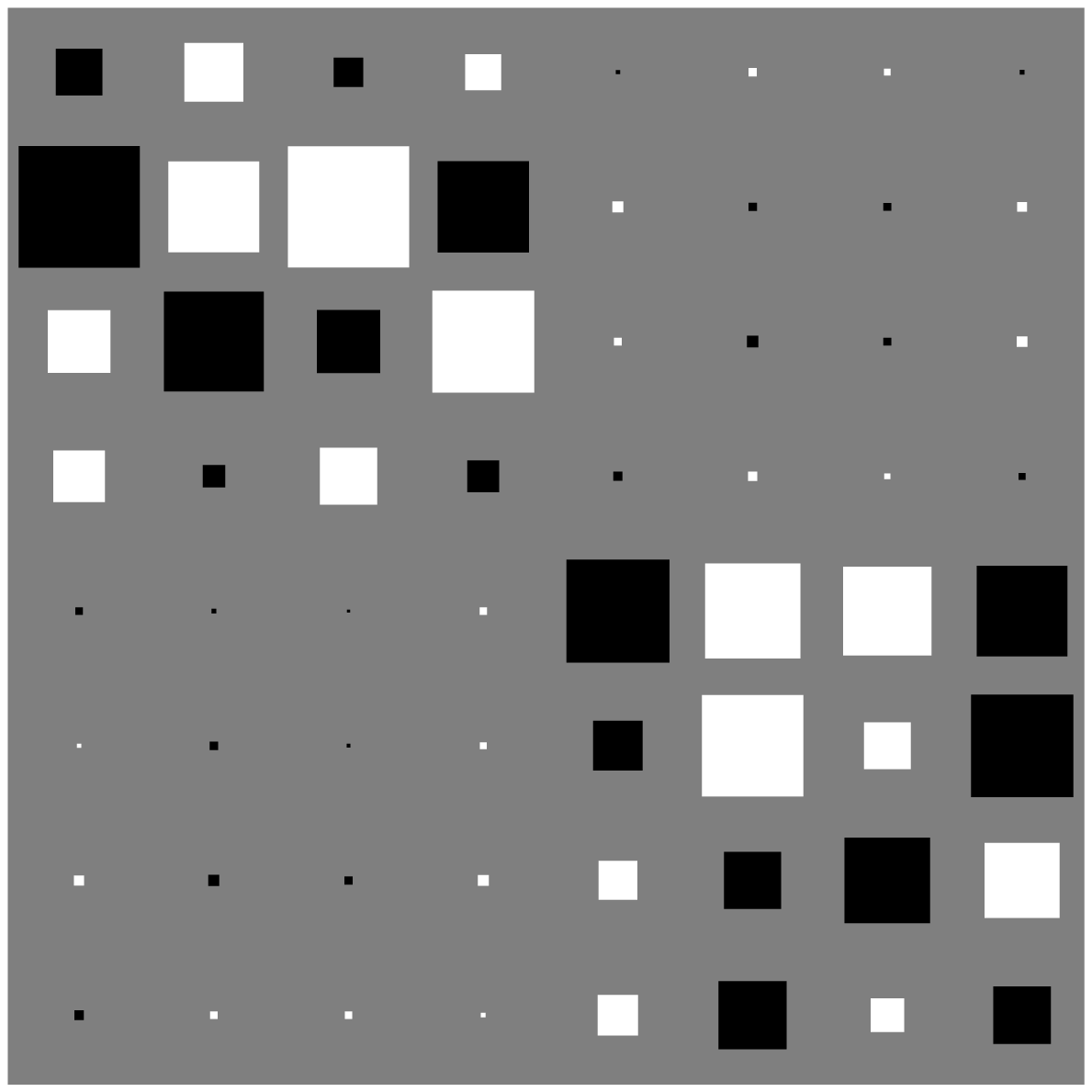}}\hspace*{2cm}%
\subfloat[][]{\includegraphics[width=3.2cm]{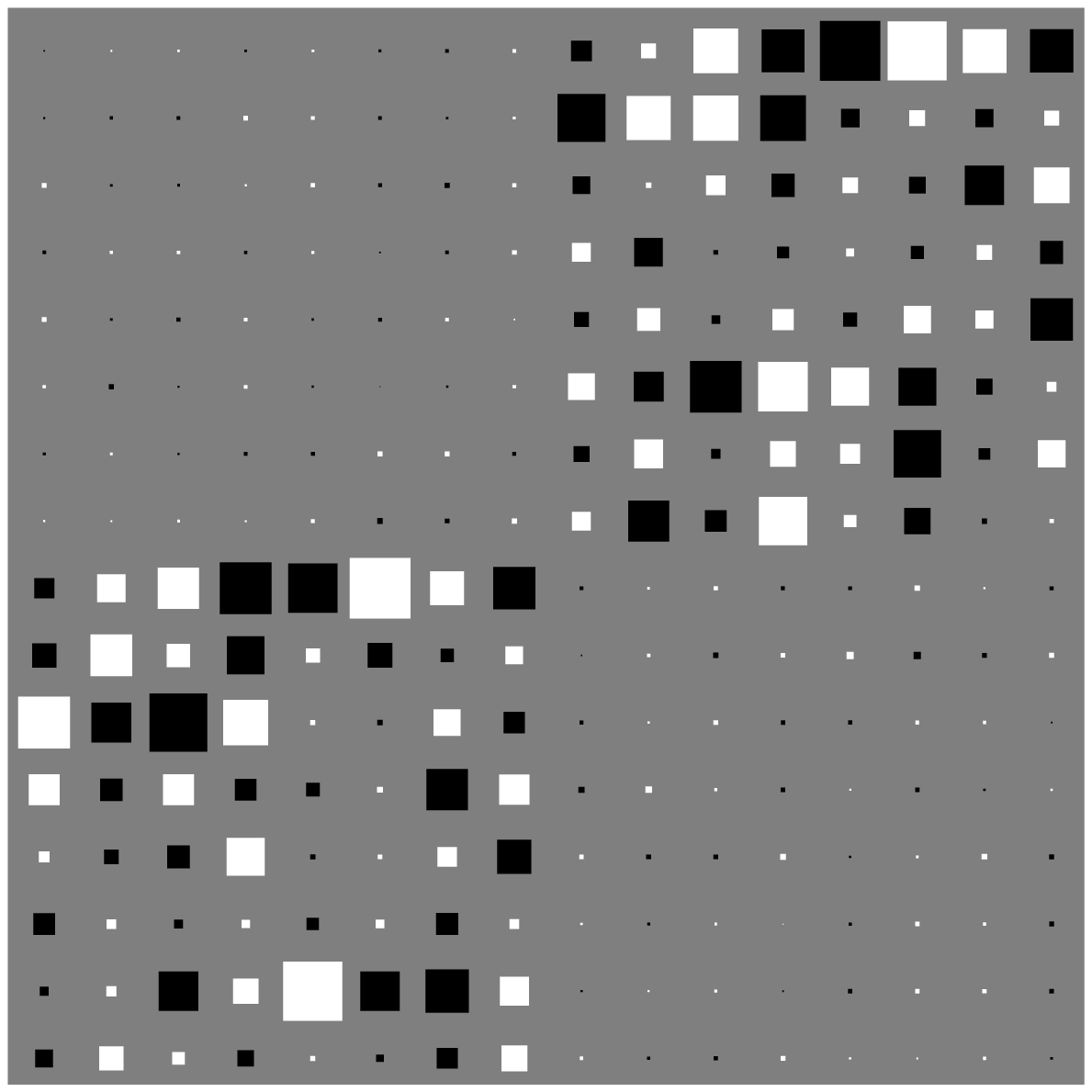}}\hspace*{2cm}%
\subfloat[][]{\includegraphics[width=3.2cm]{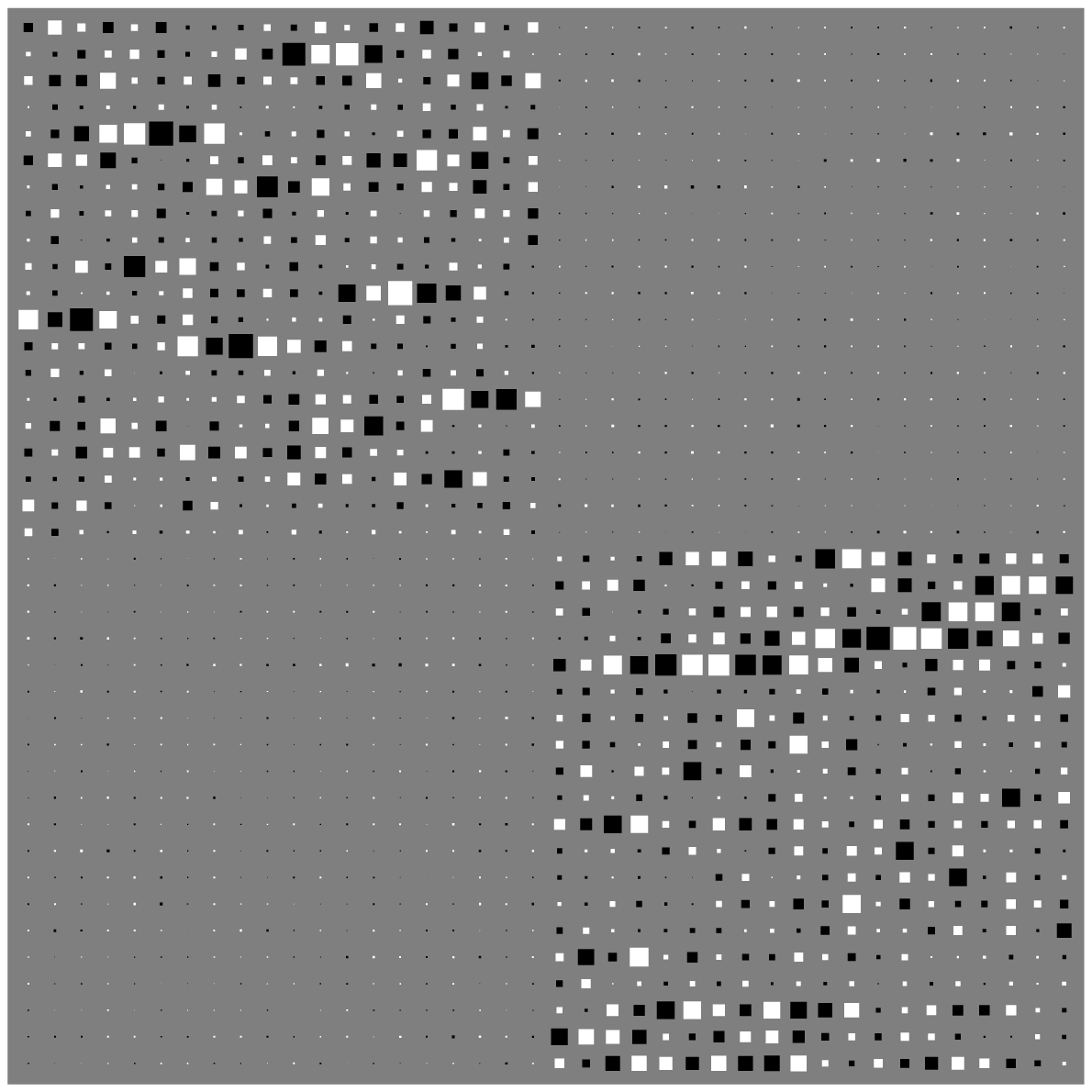}}\hfill%
\caption[]{Hinton-diagrams of the JFD methods on the \emph{Beatles} database for different convolution parameters. (a):
$L=1$ $(D_{\text{ISA}}=8)$, (b): $L=2$ $(D_{\text{ISA}}=16)$, (c): $L=5$ $(D_{\text{ISA}}=40)$. In the ideal case, $2$
pieces of $D_{\text{ISA}}/2$-dimensional blocks are formed by multiplying the mixing matrix of the derived ISA task and
the estimated demixing matrix.}
\label{fig:JFD-hinton-Beatles}%
\end{figure}

In the case of \emph{Beatles} database, test parameters were
similar to those of the \emph{ABC} database: the number and the
dimensions of the components were minimal ($d=2$, $M=2$) and the
dependence on the convolution length was tested. Parameter $L$
took values on $1,2,5,10,20,30$. The number of observations varied
between \mbox{$1,000\le T\le 75,000$}. The Amari-index and the
sweep number of the optimization are illustrated in
Figure~\ref{fig:JFD-amari-sweep-Beatles}. Precise values of the
Amari-index are provided in
Table~\ref{tab:JFD-amari-dists-Beatles}. The Hinton-diagrams are
in Figure~\ref{fig:JFD-hinton-Beatles}.

The Beatles songs are non-i.i.d sources and subsequent samples
$\b{s}^m$(t) and $\b{s}^m(t-1)$ have dependencies. Thus, in
$\b{S}(t)$ the components of the \eqref{eq:uBSSD-reduced-model}
model that belong to the same song can not be distinguished. In
the ideal case, however, the songs can be separated, that is, the
separation of the $2$ pieces of
\mbox{$(D_{\text{ISA}}/2)$-dimensional} song-subspaces is
possible.  We measure the appearance of the 2 of
\mbox{$(D_{\text{ISA}}/2)$-dimensional} blocks for the
\emph{Beatles} songs by means of the Amari-index. The results,
which demonstrate the success of our method, are shown in
Figure~\ref{fig:JFD-amari-sweep-Beatles}. It can be seen that the
JFD method found the hidden components, for $50$ and $75$ thousand
samples for $L=30$ (convolution length is $31$) too. The number of
sweeps required for the optimization of the permutations was
between $1$ and $6$ for all sample numbers \mbox{$1,000\le T \le
75,000$}, parameters $1\le L\le 30$ and for all $50$ random
initializations. According to
Table~\ref{tab:JFD-amari-dists-Beatles}, for sample number
$T=75,000$ the Amari-index is between $1.07\%$ $(L=1)$ and
$4.43\%$ $(L=30)$ on the average and has a small ($0.04-0.09$)
standard deviation. Separation of the $2$ of $D_{\text{ISA}}/2$
dimensional subspaces are illustrated in
Figure~\ref{fig:JFD-hinton-Beatles} through the Hinton-diagrams.

\begin{table}
    \centering
    \begin{tabular}{|@{\hspace{2pt}}c@{\hspace{1pt}}|@{\hspace{2pt}}c@{\hspace{1pt}}|@{\hspace{2pt}}c@{\hspace{1pt}}|@{\hspace{2pt}}c@{\hspace{1pt}}|@{\hspace{2pt}}c@{\hspace{1pt}}|@{\hspace{2pt}}c@{\hspace{1pt}}|}
    \hline
        $L=1$ & $L=2$ & $L=5$ & $L=10$ & $L=20$& $L=30$\\
    \hline\hline
        $1.07\%$ $(\pm 0.04)$ & $1.43\%$ $(\pm 0.09)$ & $2.29\%$ $(\pm 0.07)$ & $3.31\%$ $(\pm 0.06)$ & $4.03\%$ $(\pm 0.06)$& $4.43\%$ $(\pm 0.04)$\\
    \hline
    \end{tabular}
    \caption{Amari-index of the JFD method for \emph{Beatles} database for different convolution lengths: average $\pm$ deviation. Number of samples:
    $T=75,000$.
    For other sample numbers between $1,000\le T <  75,000$ see Figure~\ref{fig:JFD-amari-sweep-Beatles}(a).}
    \label{tab:JFD-amari-dists-Beatles}
\end{table}

\subsubsection{Kernel-ISA Techniques}\label{sec:kernelISA-vs-JFD}
We study the efficiency of the KCCA and KGV kernel-ISA methods of
Section~\ref{sec:kernelISA-methods}.  The kernel-ISA techniques
was found to have a higher computational burden, but they also
have advantages compared with the JFD technique for ISA tasks.

For the KCCA and KGV methods we also applied the pseudocode of
Table~\ref{tab:ISA-pseudocode}. ICA was executed by the fastICA
algorithm \cite{hyvarinen97fast}. The Gaussian kernel
$k(\b{u},\b{v})\propto
\exp\left(\frac{-\left\|\b{u}-\b{v}\right\|^2}{2\sigma^2}\right)$
was chosen for both the KCCA and KGV methods and parameter
$\sigma$ was set to $5$. In the KCCA method, regularization
parameter $\kappa=10^{-4}$ was applied. In the experiments,
parameters $(\sigma,\kappa)$ were proved to be reasonably robust.
Mixing matrix $\b{A}$ was generated randomly from the orthogonal
group and the sample number was chosen from the interval $100 \le
T\le 5,000$.

\begin{figure}%
\centering%
\subfloat[][]{\includegraphics[width=7.9cm]{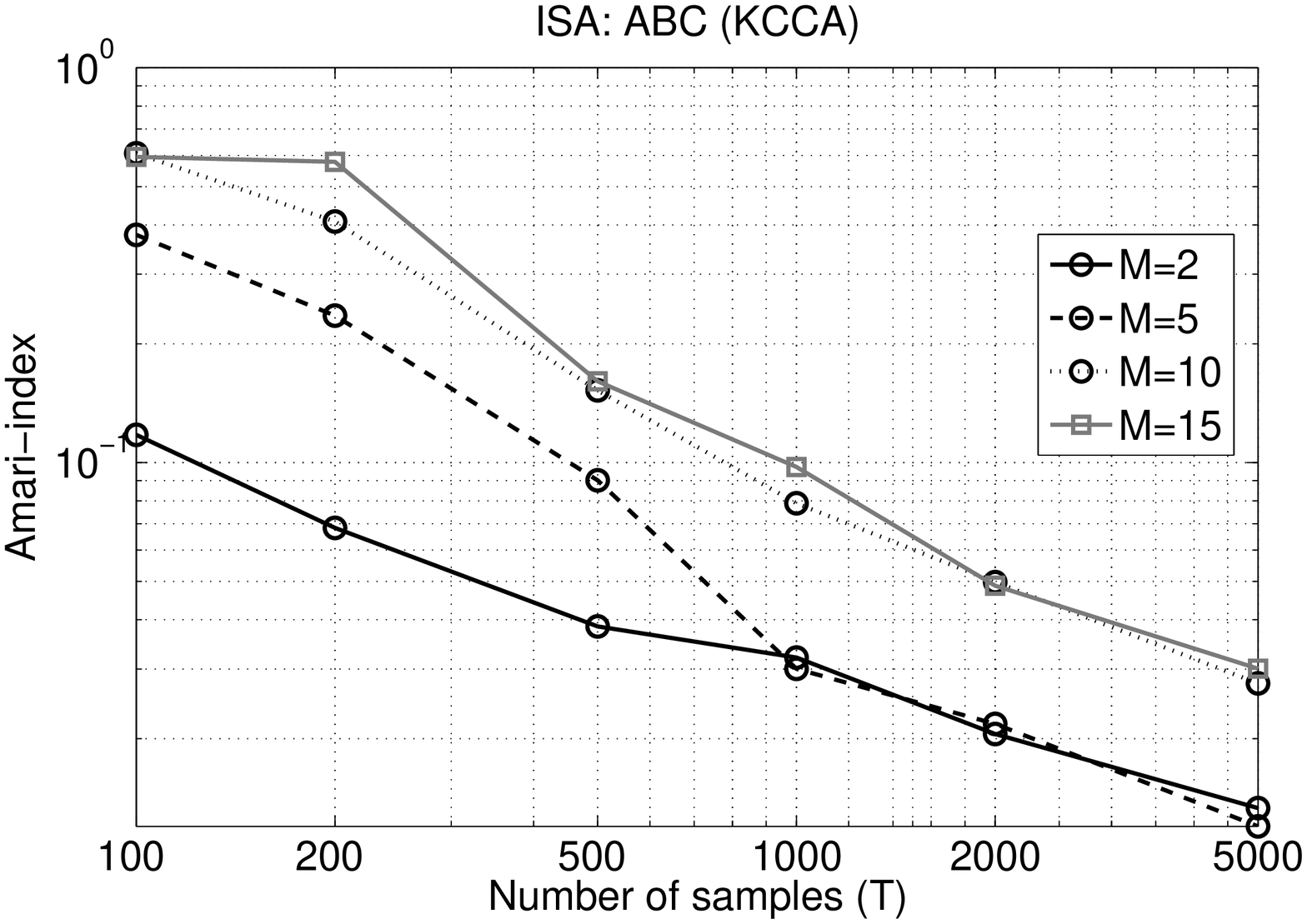}}%
\subfloat[][]{\includegraphics[width=7.9cm]{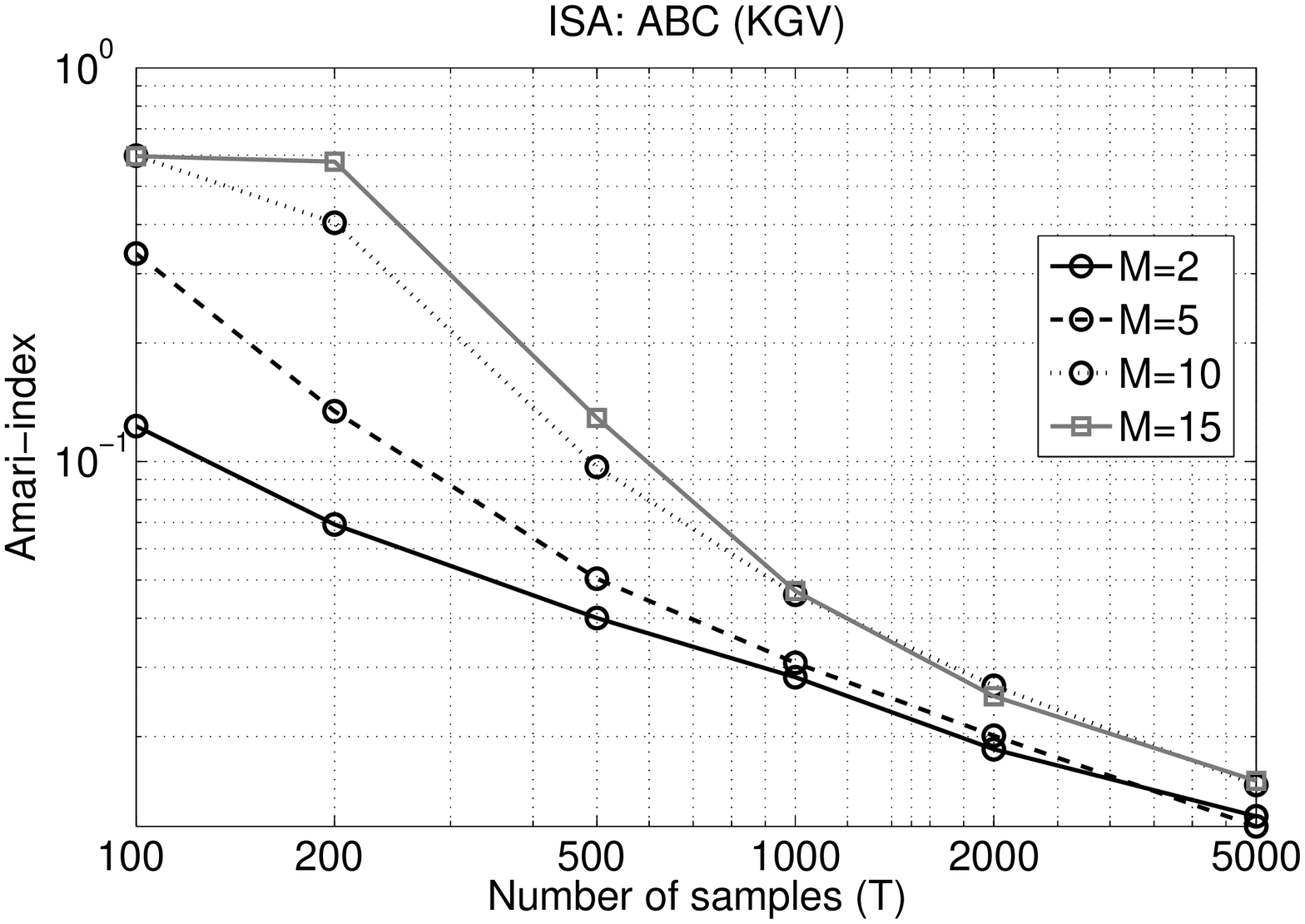}}\\%
\subfloat[][]{\includegraphics[width=7.9cm]{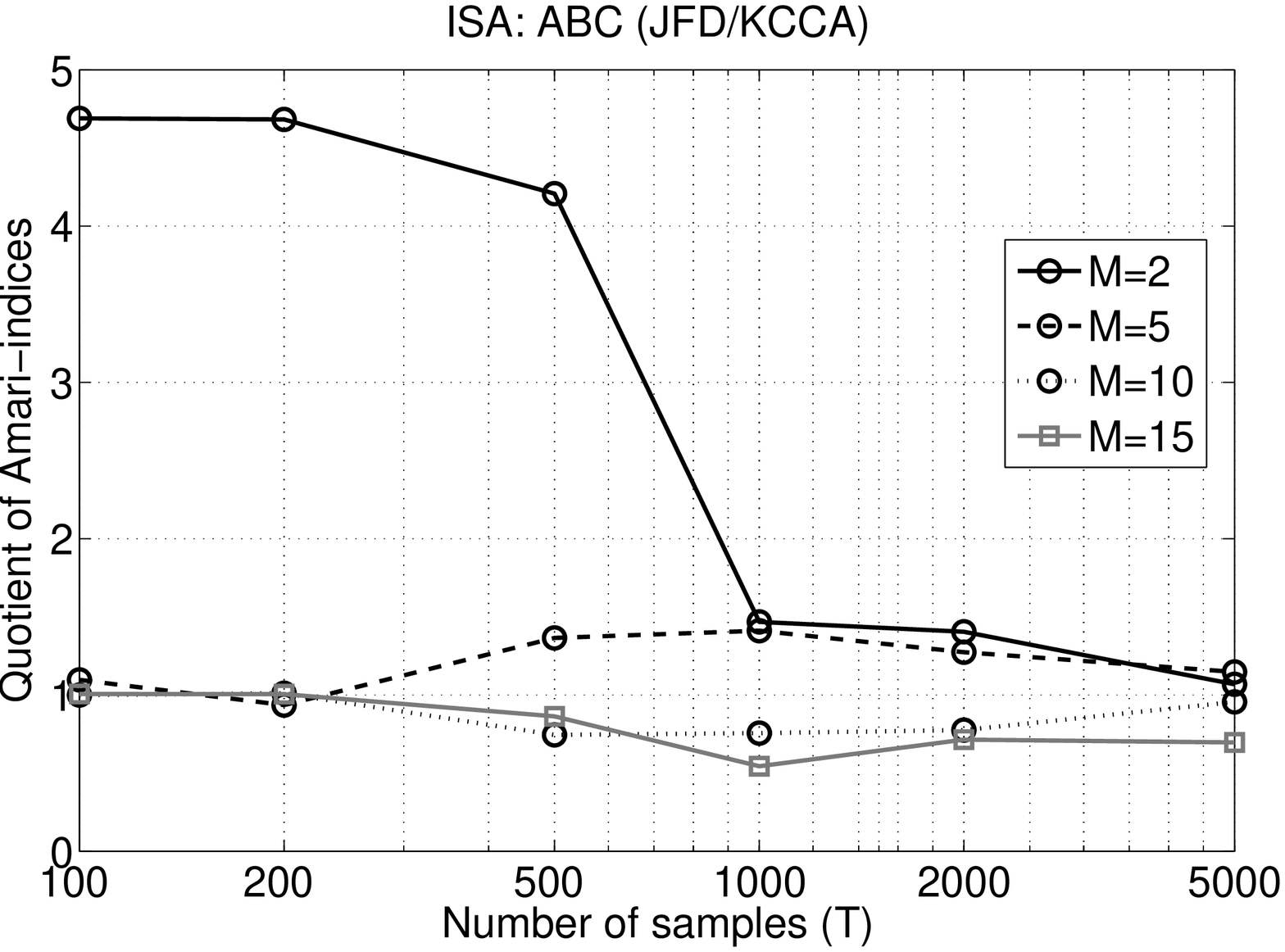}}%
\subfloat[][]{\includegraphics[width=7.9cm]{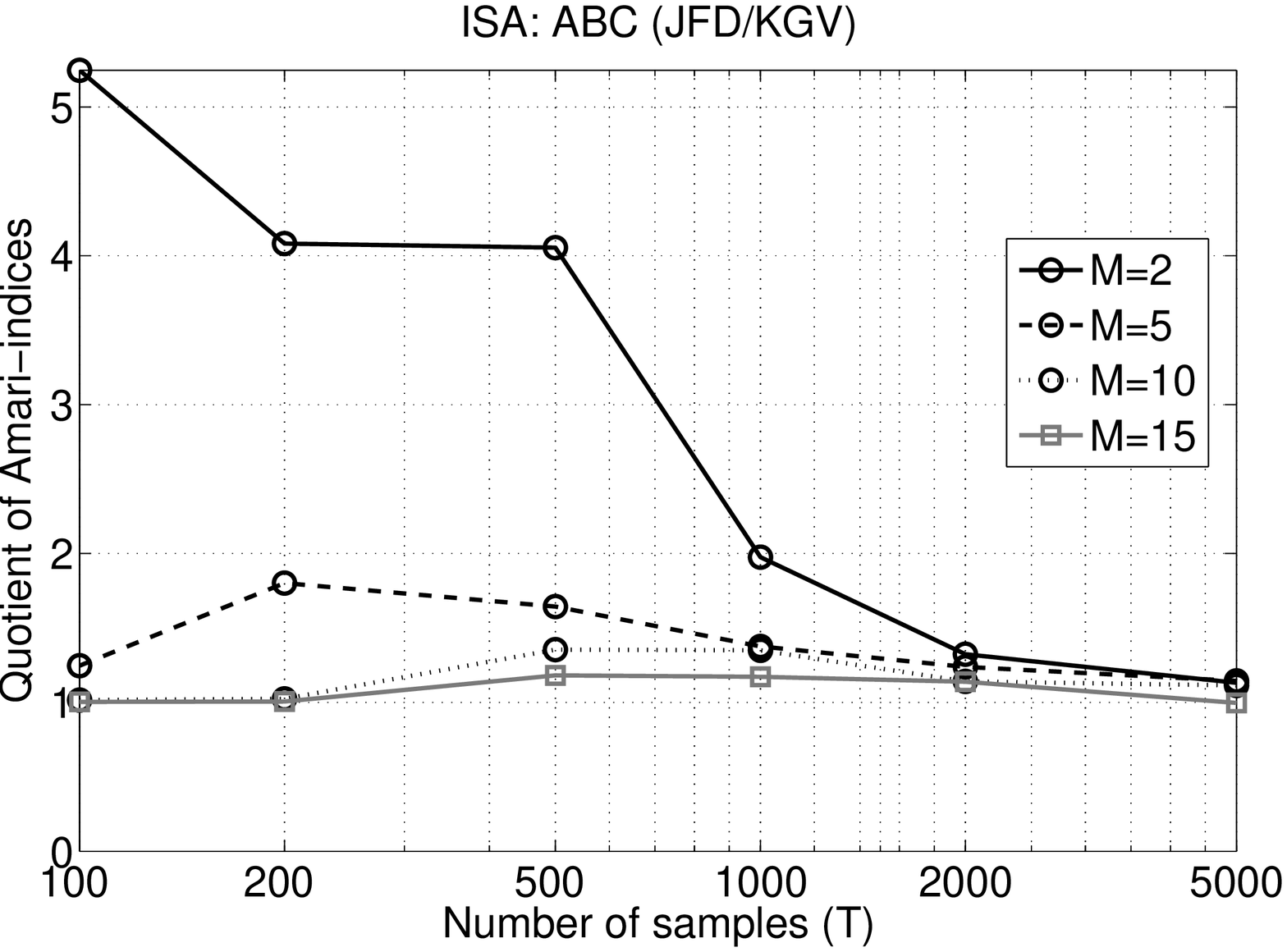}}\\%
\caption[]{(a) and (b): Amari-index of the KCCA  and KGV methods,
respectively, for the \emph{ABC} database as a function of sample
number and for different numbers of components $M$. (c) and (d):
Amari-index of the JFD method is divided by the Amari-index of the
KCCA method and
the KGV technique, respectively. For values larger (smaller) than 1 the kernel-ISA method is better (worse) than the JFD method.}%
\label{fig:KCCA-KGV-ABC-amari-vsJFD}%
\end{figure}

\begin{figure}%
\centering%
\subfloat[][]{\includegraphics[width=7.9cm]{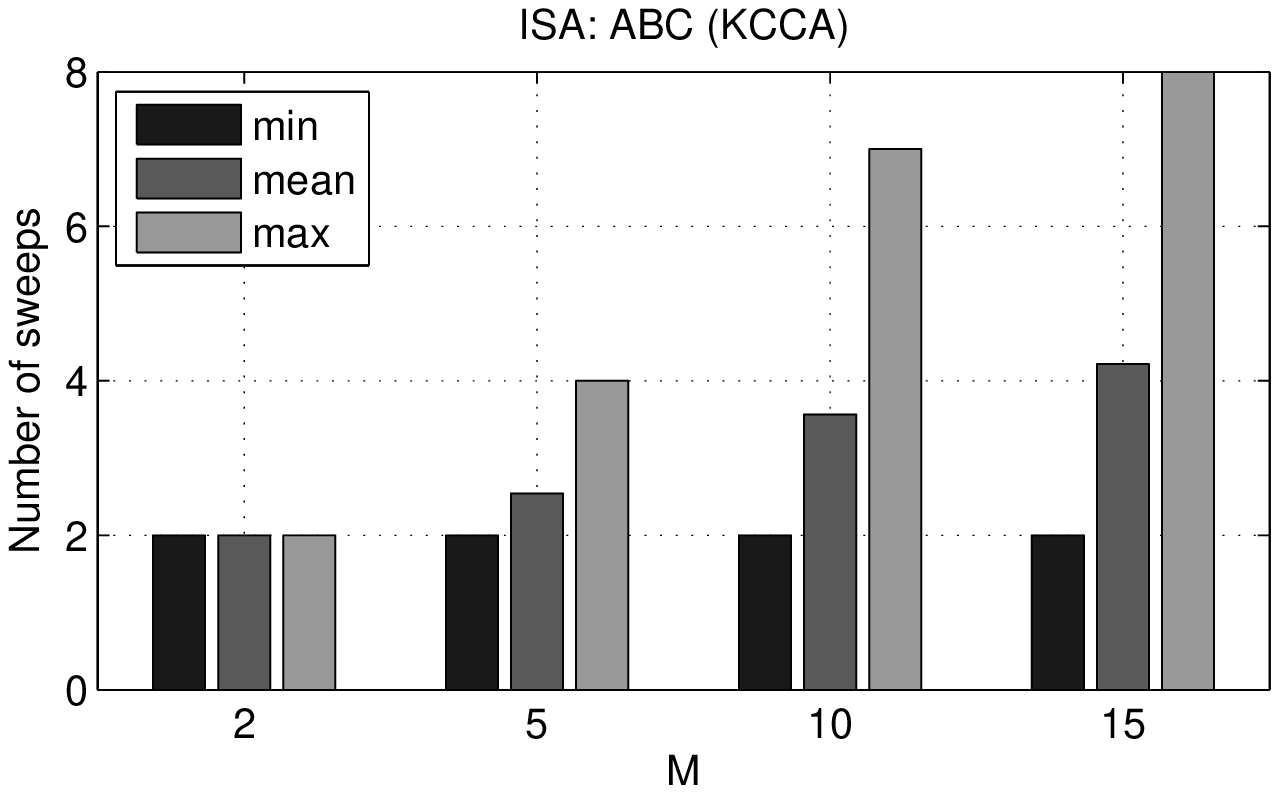}}%
\subfloat[][]{\includegraphics[width=7.9cm]{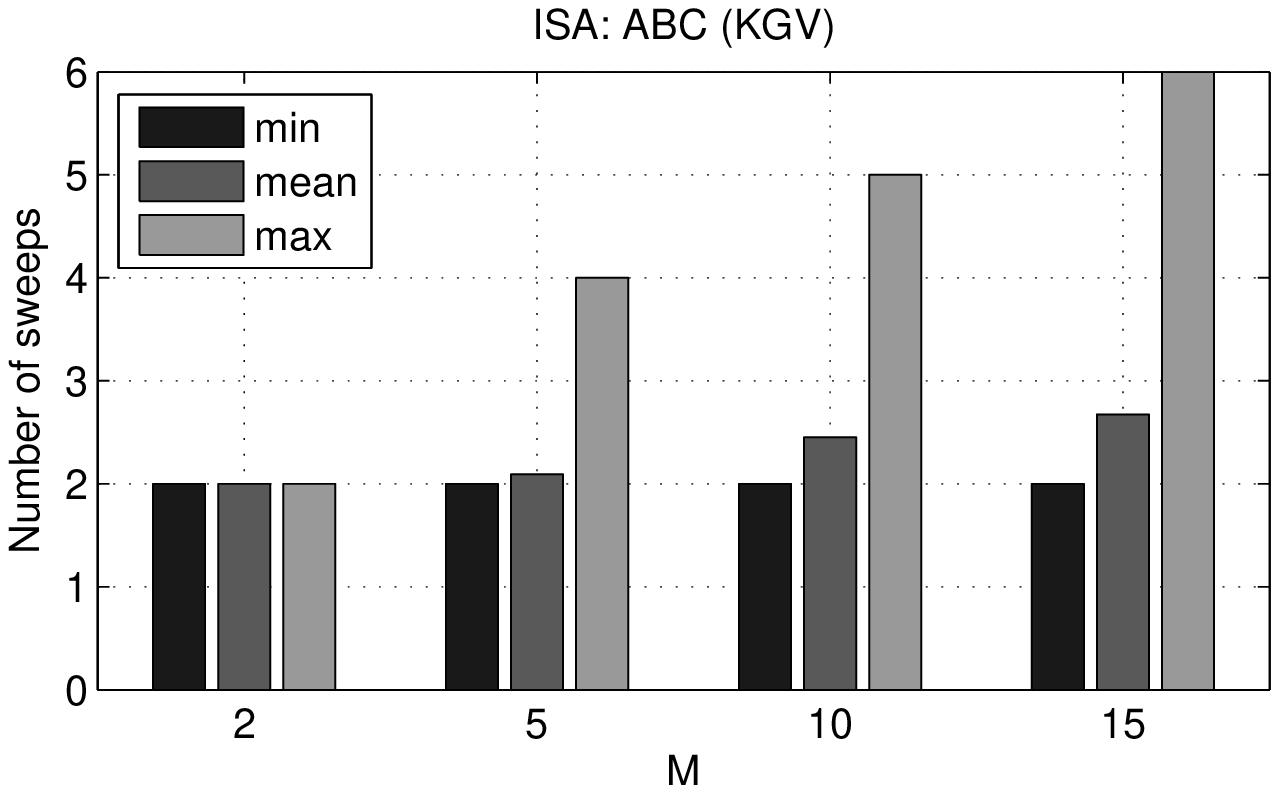}}%
\caption[]{Number of sweeps for the KCCA and the KGV methods needed for the optimization of permutations as a function
of component number $M$ on the \emph{ABC} database. (a): KCCA method, (b): KGV method. Black: minimum, gray: mean, light gray: maximum.}%
\label{fig:KCCA-KGV-sweep-ABC}%
\end{figure}

\begin{figure}%
\centering%
\hfill\subfloat[][]{\includegraphics[width=4cm]{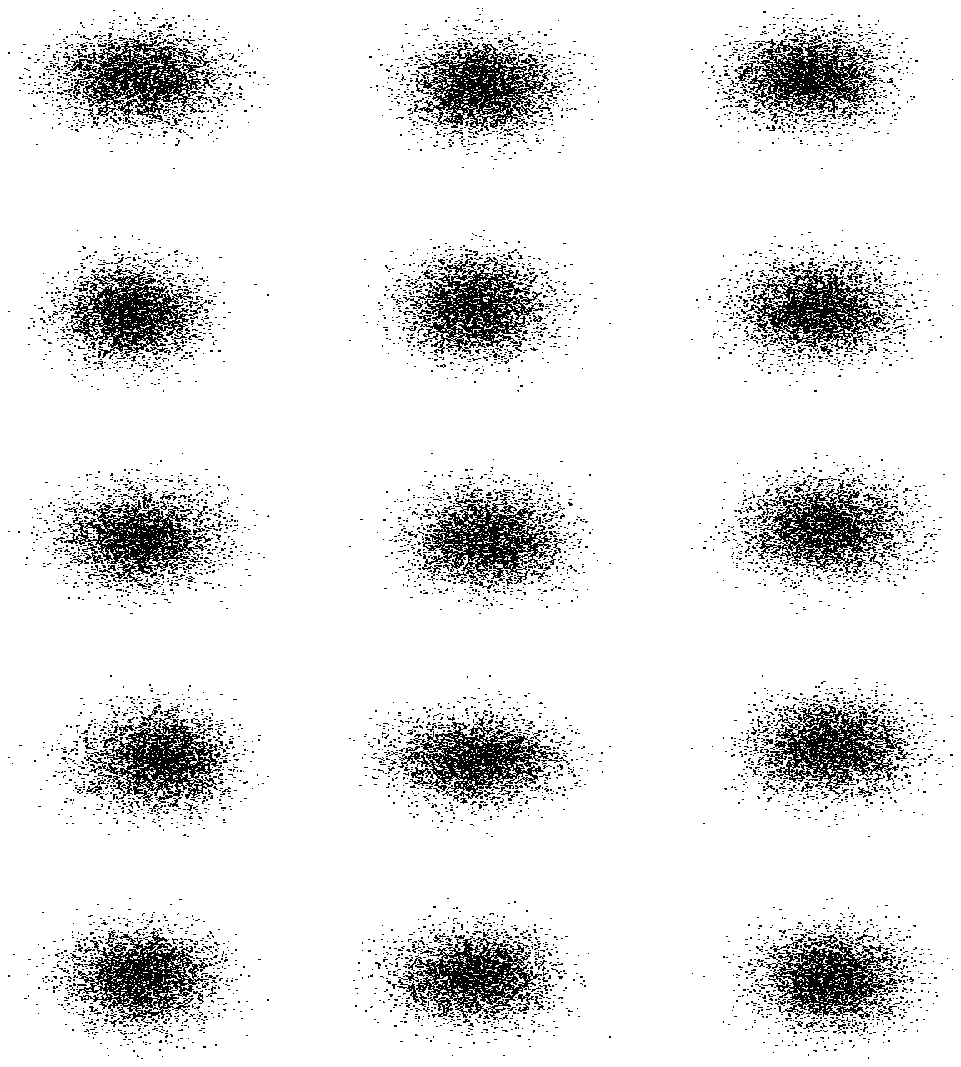}}\hfill%
\subfloat[][]{\includegraphics[width=4.15cm]{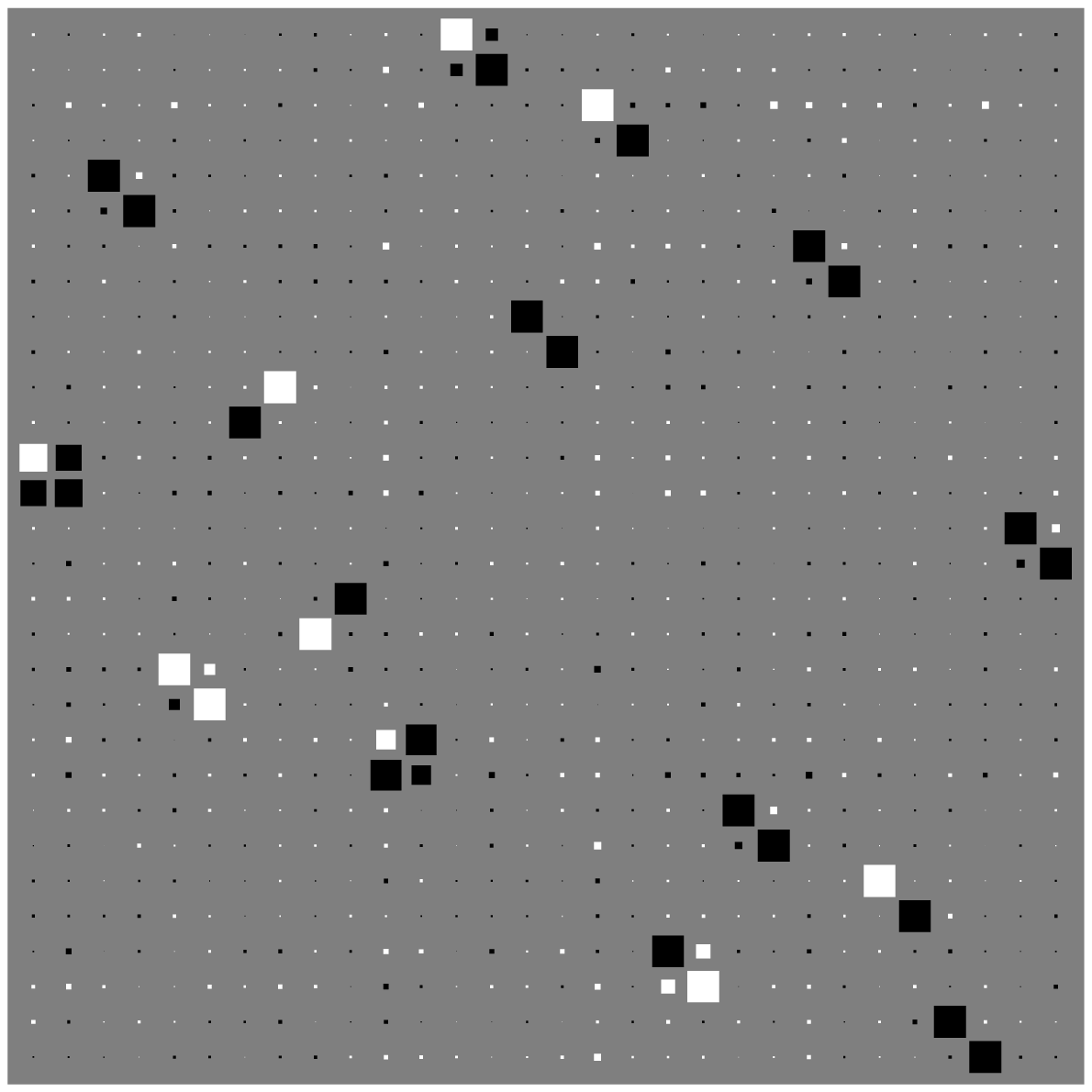}}\hfill%
\subfloat[][]{\includegraphics[width=4cm]{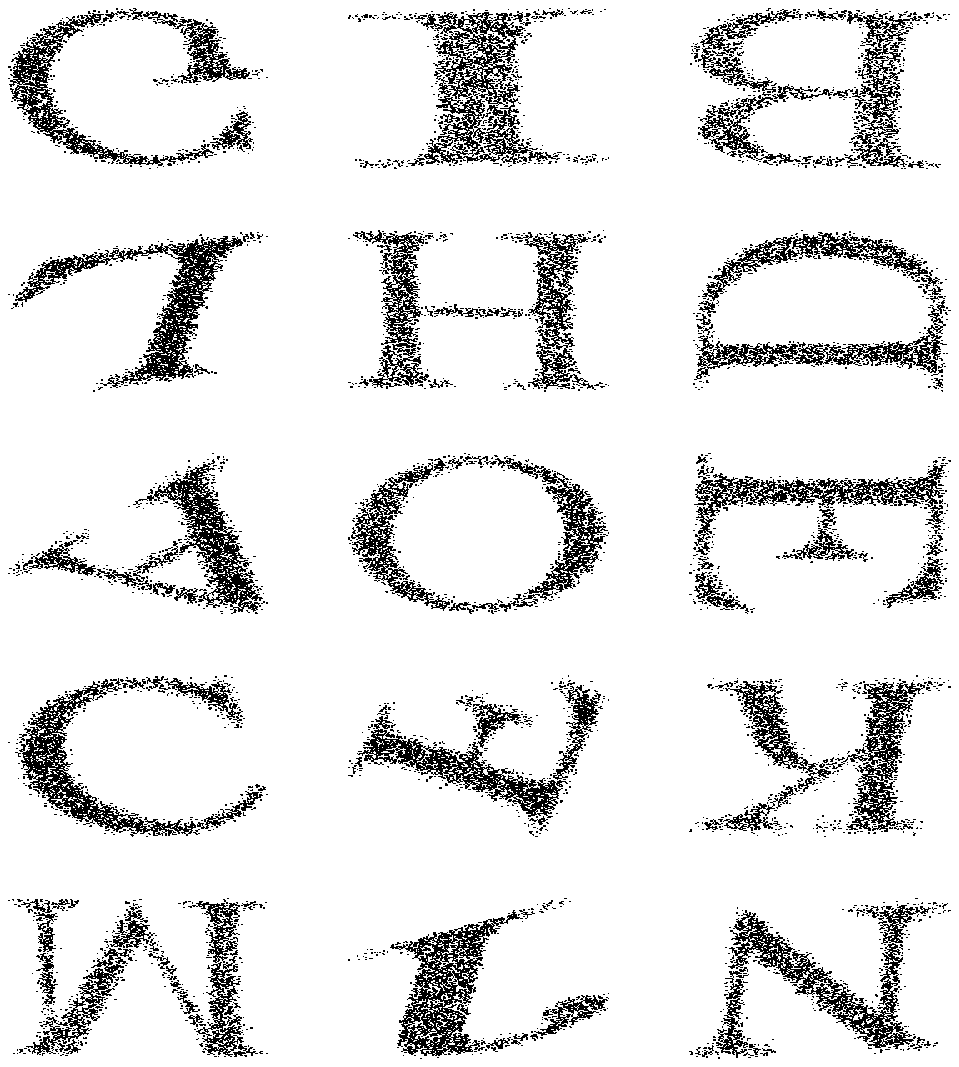}}\hfill\hfill%
\caption[]{Illustration of the KCCA method for the \emph{ABC}
database. Sample number: \mbox{$T=5,000$}. (a): observed mixed
signals $\b{x}(t)$. (b) Hinton-diagram: the product of the mixing
matrix of the derived ISA task and the estimated demixing matrix
(= approximately block-permutation matrix). (c): estimated
components.
Hidden components are recovered up to permutation and orthogonal transformation.}%
\label{fig:KCCA-ABC-demo}%
\end{figure}

Our first ISA example concerns the \emph{ABC} database. The
dimension of a component was $d=2$ and the number of the
components $M$ took different values ($M=2,5,10,15$). Precision of
the estimations is shown in
Figure~\ref{fig:KCCA-KGV-ABC-amari-vsJFD}, where the precision of
the JFD method on the same database is also depicted. The number
of sweeps required for the optimization of the permutations is
shown in Figure~\ref{fig:KCCA-KGV-sweep-ABC} for different sample
numbers and for different component numbers. The data are averaged
over $50$ random estimations. Figure~\ref{fig:KCCA-ABC-demo}
depicts the KCCA estimation for the \emph{ABC} database.

Figure~\ref{fig:KCCA-KGV-ABC-amari-vsJFD} shows that the KCCA and
KGV kernel-ISA methods give rise to high precision estimations on
the \emph{ABC} database even for small sample numbers. The KGV
method was more precise for all $M$ values studied than the JFD
method. The ratio of precisions could be as high as $4$ (see
sample number $500$). The KCCA method is somewhat weaker. For
smaller tasks ($M=2$ and $5$) and for small sample numbers it also
exceeds the precision of the JFD method. Precision of the method
become about the same for higher sample numbers and larger tasks.
Sweep numbers of the  KCCA (KGV) method were between $2$ and $8$
($2$ and $6$) (Figure~\ref{fig:KCCA-KGV-sweep-ABC}). Note that one
sweep is always necessary for our procedure
(Table~\ref{tab:ISA-pseudocode}). A single sweep may be
satisfactory if---by chance---the ICA provides the correct
permutation.
\begin{table}
    \centering
    \begin{tabular}{|c|c|c|c|c|}
    \hline
        &$M=2$ & $M=5$ & $M=10$ & $M=15$\\
    \hline\hline
        KCCA &$1.33\%$ $(\pm 0.48)$ & $1.20\%$ $(\pm 0.17)$ & $2.76\%$ $(\pm 2.86)$ & $3.00\%$ $(\pm 2.21)$\\
    \hline
        KGV &$1.26\%$ $(\pm 0.54)$ & $1.18\%$ $(\pm 0.17)$ & $1.51\%$ $(\pm 0.31)$ & $1.54\%$ $(\pm 0.34)$\\
    \hline
    \end{tabular}
    \caption{Amari-index for the KCCA and the KGV methods for database \emph{ABC}, for different component number $M$: average $\pm$ deviation.
    Number of samples: $T=5,000$. For other sample numbers between $100\le T < 5,000$, see
    Figure~\ref{fig:KCCA-KGV-ABC-amari-vsJFD}.}
    \label{tab:KCCA-KGV-amari-dists-ABC}
\end{table}

\begin{figure}%
\centering%
\subfloat[][]{\includegraphics[width=8cm]{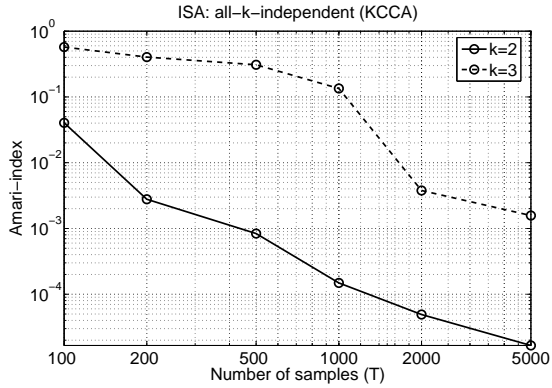}}%
\subfloat[][]{\includegraphics[width=8cm]{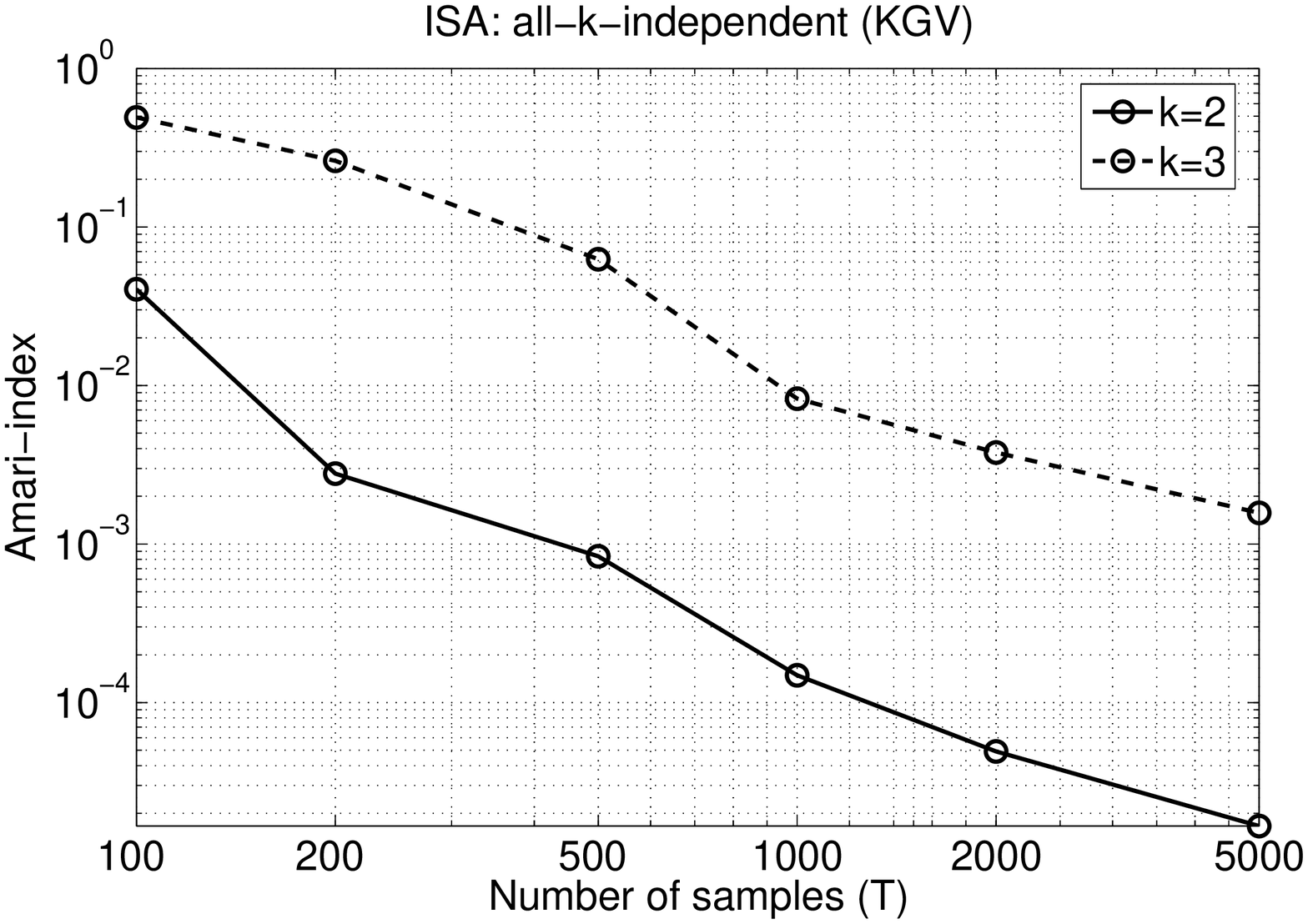}}\\%
\subfloat[][]{\includegraphics[width=8cm]{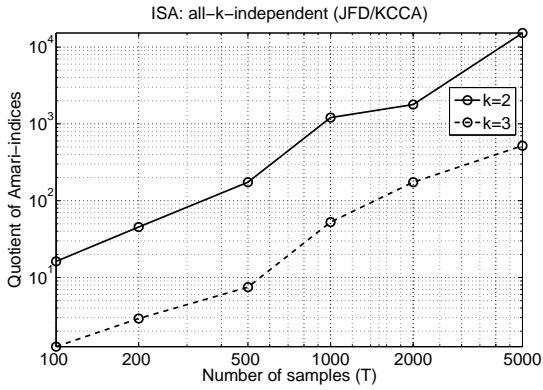}}%
\subfloat[][]{\includegraphics[width=8cm]{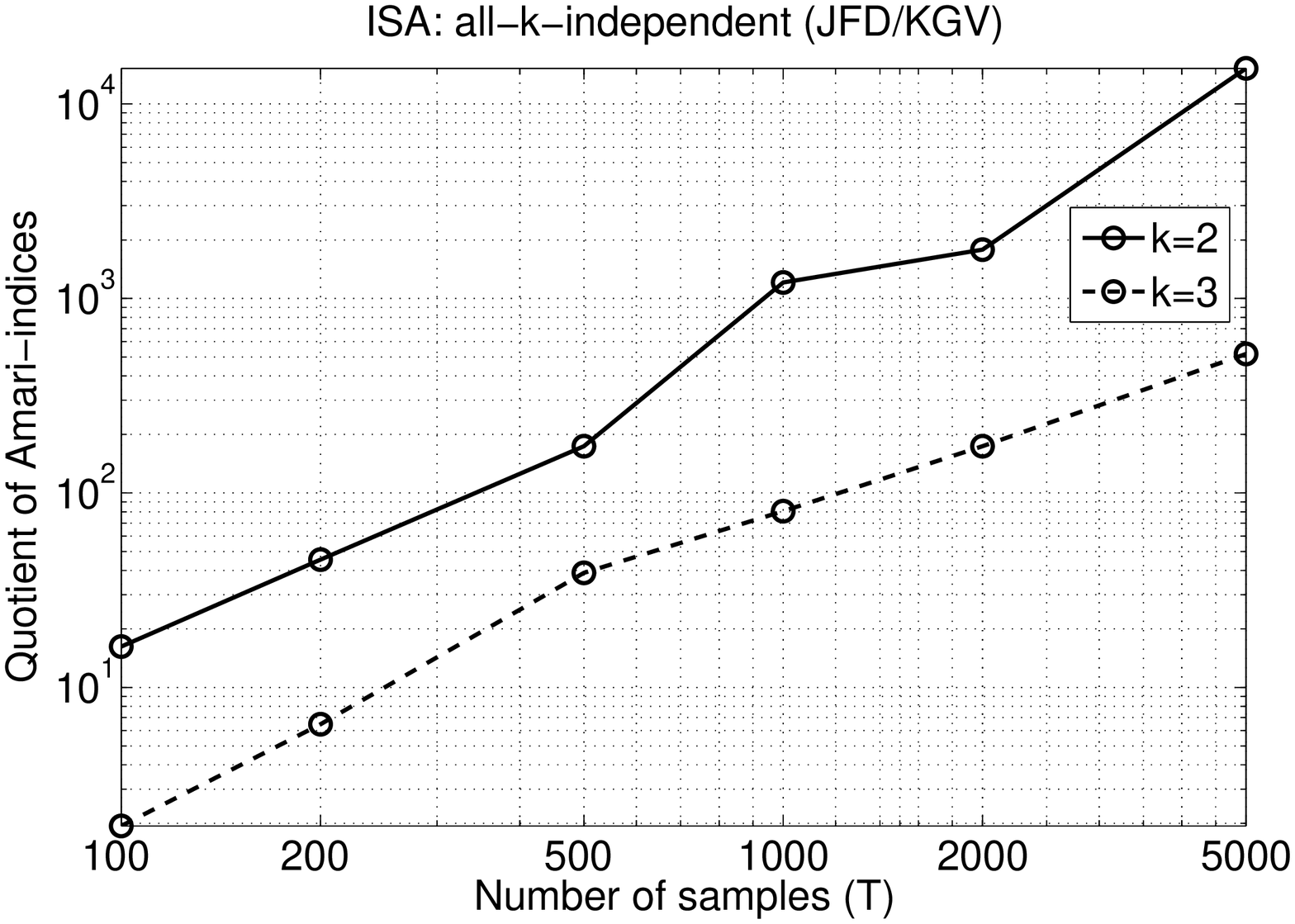}}\\%
\caption[]{(a) and (b): Amari-indices of the KCCA and KGV methods, respectively,  as a function of the sample number and for $k=2$ and $3$
on the \emph{all-k-independent} database. For more details, see Table~\ref{tab:KCCA-KGV-amari-dists-all-k-independent}. (c):
ratio of Amari-indices of JFD and KCCA methods, (d): ratio of Amari-indices of JFD and KGV methods.}%
\label{fig:KCCA-KGV-vsJFD-amari-all-k-independent}%
\end{figure}

\begin{figure}%
\centering%
\subfloat[][]{\includegraphics[width=8cm]{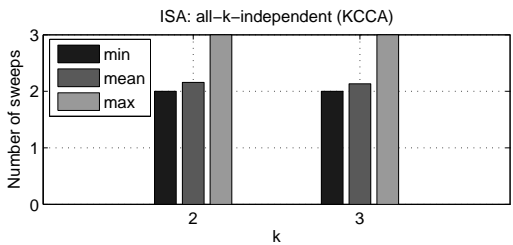}}%
\subfloat[][]{\includegraphics[width=8cm]{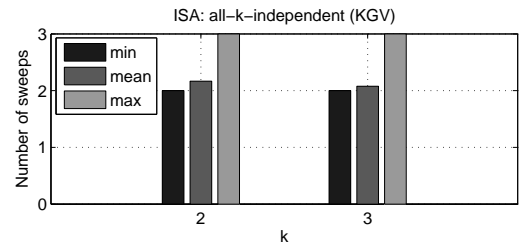}}%
\caption[]{Number of sweeps of permutation optimization for the KCCA (a) and KGV (b) methods for the
\emph{all-k-independent} database and for different $k$ values. Black: minimum, gray: mean, light gray: maximum.}%
\label{fig:KCCA-KGV-vsJFD-sweep-all-k-independent}%
\end{figure}

The other illustration concerns the \emph{all-k-independent}
database. This test can be difficult for ISA methods
\cite{szabo06cross}. Number of components $M$ was $2$. For
$k=2,3$, when the dimension of the components $d=3$ and $4$,
respectively, the KCCA and KGV kernel-ISA methods efficiently
estimated the hidden components. The precision of the KCCA and KGV
estimations as well as the comparison with the JFD method are
shown in Figure~\ref{fig:KCCA-KGV-vsJFD-amari-all-k-independent}.
The average number of sweeps required for the optimization of the
permutations for different $k$ values and for $50$ randomly
initialized computations is provided in
Figure~\ref{fig:KCCA-KGV-vsJFD-sweep-all-k-independent}. The
values of the Amari-indices are shown in
Table~\ref{tab:KCCA-KGV-amari-dists-all-k-independent}.

According to
Figure~\ref{fig:KCCA-KGV-vsJFD-amari-all-k-independent} the two
kernel-based methods exhibit similar precision. Both were superior
to the JFD technique. The ratio of the Amari-indices for sample
number $5,000$ and for $k=2$ is more than $15,000$, for $k=3$ it
is more than $500$. For details concerning the Amari-indices, see
Table~\ref{tab:KCCA-KGV-amari-dists-all-k-independent}. These
indices are close to each other for the KCCA and the KGV methods:
$0.0017\%$ for $k=2$, $0.16\%$ for $k=3$ on average. Both
kernel-ISA methods used $2-3$ sweeps for the optimization of the
permutations
(Figure~\ref{fig:KCCA-KGV-vsJFD-sweep-all-k-independent}).

\begin{table}
    \centering
    \begin{tabular}{|c|c|c|}
    \hline
        &$k=2$ & $k=3$\\
    \hline\hline
        KCCA/KGV &$0.0017\%$ $(\pm 0.0014)$ & $0.16\%$ $(\pm 0.04)$\\
    \hline
    \end{tabular}
    \caption{The Amari-index of the KCCA and KGV methods for database \emph{all-k-independent} for different $k$ values:
    average $\pm$ deviation. Number of samples: $T=5,000$. For other sample numbers between $100\le T < 5,000$, see
    Figure~\ref{fig:KCCA-KGV-vsJFD-amari-all-k-independent}.}
    \label{tab:KCCA-KGV-amari-dists-all-k-independent}
\end{table}

\section{Conclusions}\label{sec:conclusions}
We have introduced a new model, the blind subspace deconvolution
(BSSD) for data analysis. This model deals with the casual
convolutive mixture of multidimensional independent sources. The
undercomplete version (uBSSD) of the task has been presented, and
it has been shown how to derive an independent subspace analysis
(ISA) task from the uBSSD problem. Recent developments of the ISA
techniques enabled us to handle the emerging high dimensional
problems. Our earlier results, namely the ISA Separation Theorem
\cite{szabo06separation} motivated us to reduce the ISA task to
the search for the optimal permutation of the ICA components. The
components were grouped with a novel joint decorrelation
technique, the joint \mbox{f-decorrelation} (JFD) method
\cite{szabo06real}.

Also, we adapted other ICA techniques, such as the KCCA and KGV
methods to the ISA task and studied their efficiency. Simulations
indicated that although the KCCA and KGV methods give rise to
serious computational burden relative to the JFD method, they can
be advantageous for smaller ISA tasks and for ISA tasks when the
number of samples is small.

Finally, we note that we achieved small errors in these high dimensional computations. These small errors indicate that
the Separation Theorem is robust and might be extended to a larger class of noise sources.

\appendix

\section{Proof of the ISA Separation Theorem}\label{sec:RISA-sep-theorem}
We shall rely on entropy inequalities
(Section~\ref{sec:REPIineq-s}). In
Section~\ref{sec:connection2RICA} connection to the  ICA cost
function is derived (Lemma~\ref{R-prop}). The ISA Separation
Theorem then follows.

\subsection{EPI-type Relations}\label{sec:REPIineq-s} First, consider the so-called entropy
power inequality (EPI)
\begin{equation*}
e^{2H\left(\sum_{i=1}^Lu_i\right)}\ge \sum_{i=1}^L e^{2H(u_i)},
\end{equation*}
where $u_1,\ldots,u_L\,\in\R$ denote continuous random variables.
This inequality holds, for example, for independent continuous
variables \cite{cover91elements}.

If EPI is satisfied on the surface of the \mbox{$L$-dimensional}
unit sphere $\S^L$, then a further inequality holds:
\begin{lemma}\label{lem:suff}
Suppose that continuous random variables $u_1,\ldots,u_L \,\in\R$
satisfy the following inequality
\begin{equation}
    e^{2H\left(\sum_{i=1}^Lw_iu_i\right)}\ge \sum_{i=1}^L e^{2H(w_iu_i)}, \forall\b{w}\in \S^L.\label{eq:w-EPI}
\end{equation}
This inequality will be called the \emph{w-EPI} condition. Then Equation~\eqref{eq:suff} holds, too.
\end{lemma}

\begin{proof2} Assume that $\b{w}\in \S^L$. Applying
$\ln$ on condition \eqref{eq:w-EPI}, and using the monotonicity of
the $\ln$ function, we can see that the first inequality is valid
in the following inequality chain
\begin{align*}
2H\left(\sum_{i=1}^Lw_iu_i\right)&\ge \ln\left(\sum_{i=1}^L
e^{2H(w_iu_i)}\right)=\ln\left(\sum_{i=1}^Le^{2H(u_i)}
w_i^2\right)\\
&\ge\sum_{i=1}^Lw_i^2\ln\left(e^{2H(u_i)}\right)=2\sum_{i=1}^Lw_i^2H(u_i).
\end{align*}
Here,
\begin{enumerate}
    \item
        we used the relation
        \begin{equation*}
            H(w_iu_i)=H(u_i)+\ln\left(\left|w_i\right|\right)
        \end{equation*}
         for the entropy of the transformed variable \cite{cover91elements}. Hence
        \begin{equation*}
            e^{2H(w_iu_i)}=e^{2H(u_i)+2\ln\left(\left|w_i\right|\right)}=e^{2H(u_i)}
            e^{2\ln\left(\left|w_i\right|\right)}=e^{2H(u_i)}
            w_i^2.\label{eq:2entr-transf}
        \end{equation*}
    \item
        In the second inequality, the concavity of $\ln$ was exploited. $\square$
\end{enumerate}
\end{proof2}

\begin{note}
w-EPI holds, for example, for independent variables $u_i$, because independence is not affected by
multiplication with a constant.
\end{note}

\subsection{Connection to the Cost Function of the ICA
Task}\label{sec:connection2RICA} The ISA Separation Theorem will
be a corollary of the following claim:
\begin{lemma}\label{R-prop}
Let
$\b{y}=\left[\b{y}^1;\ldots;\b{y}^M\right]=\b{y}(\b{W})=\b{W}\b{s}\in
\R^D$, where $\b{W}\in \mathcal{O}^D$ ($D=\sum_{m=1}^Md_m$),
$\b{y}^m\in\R^{d_m}$ is the estimation of the $m^{th}$ component
of the ISA task. Let $y^m_i\in\R$ be the $i^{th}$ coordinate of
the $m^{th}$ component ($i=1,\ldots,d_m$). Similarly, let
$s^m_i\in\R$ stand for the $i^{th}$ coordinate of the $m^{th}$
source. Let us assume that the $\b{u}:=\b{s}^m\in\R^{d_m}$ sources
satisfy the condition~\eqref{eq:suff}. Then
\begin{equation}\label{eq:main-prop}
\sum_{m=1}^M\sum_{i=1}^{d_m}H\left(y^m_i\right)\ge
\sum_{m=1}^M\sum_{i=1}^{d_m}H\left(s^m_i\right).
\end{equation}
\end{lemma}

\begin{proof2}
Let us denote the $(i,j)^{th}$ element of matrix $\b{W}$ by
$W_{i,j}$. Coordinates of $\b{y}$ and $\b{s}$ will be denoted by
$y_i$ and $s_i$, respectively. Further, let $\G^m$ denote the
indices of the $m^{th}$ subspace ($m=1,\ldots,M$), that is,
$\G^m:=\{1+\sum_{i=1}^{m-1}d_i,\ldots,\sum_{i=1}^{m}d_i\}$
$(d_0:=0)$. Now, writing the elements of the $i^{th}$ row of
matrix multiplication $\b{y}=\b{W}\b{s}$, we have
\begin{equation}
y_i=\sum_{j\in \G^1} W_{i,j}s_j+\ldots+\sum_{j\in \G^M}
W_{i,j}s_j\label{eq:y=Ws}
\end{equation}
and thus,
\begin{align}
H\left(y_i\right)&=H\left(\sum_{m=1}^M\sum_{j\in \G^m} W_{i,j}s_j\right)\label{eq:H(y=Ws)}\\
&=H\left(\sum_{m=1}^M\left[\left(\sum_{l\in\G^m}W_{i,l}^2\right)^{\frac{1}{2}}\frac{\sum_{j\in\G^m}W_{i,j}s_j}{\left(\sum_{l\in\G^m}W_{i,l}^2\right)^{\frac{1}{2}}}\right]
\right)\label{eq:w^2-in}\\
&\ge\sum_{m=1}^M\left[\left(\sum_{l\in\G^m}W_{i,l}^2\right)H\left(\frac{\sum_{j\in\G^m}W_{i,j}s_j}{\left(\sum_{l\in\G^m}W_{i,l}^2\right)^{\frac{1}{2}}}\right)\right]
\label{eq:Lem2-applied}\\
&=\sum_{m=1}^M\left[\left(\sum_{l\in\G^m}W_{i,l}^2\right)
H\left(\sum_{j\in\G^m}\frac{W_{i,j}}{\left(\sum_{l\in\G^m}W_{i,l}^2\right)^{\frac{1}{2}}}s_j\right)\right]
\label{eq:Lem2-applied-again-pre}\\
&\ge\sum_{m=1}^M\left[\left(\sum_{l\in\G^m}W_{i,l}^2\right)
\sum_{j\in\G^m}\left(\frac{W_{i,j}}{\left(\sum_{l\in\G^m}W_{i,l}^2\right)^{\frac{1}{2}}}\right)^2H\left(s_j\right)\right]\label{eq:Lem2-applied-again}\\
&=\sum_{j\in\G^1}W_{i,j}^2H\left(s_j\right)+\ldots+\sum_{j\in\G^M}W_{i,j}^2H\left(s_j\right).\label{eq:H(yi-last)}
\end{align}
The above steps can be justified as follows:
\begin{enumerate}
    \item
        \eqref{eq:H(y=Ws)}: Equation~\eqref{eq:y=Ws} was inserted into the argument of $H$.
    \item
        \eqref{eq:w^2-in}: New terms were added for Lemma~\ref{lem:suff}.
    \item
        \eqref{eq:Lem2-applied}: Sources $\b{s}^m$ are independent of one another and this independence is preserved upon
        mixing \emph{within} the subspaces, and we could also use Lemma~\ref{lem:suff}, because $\b{W}$
        is an orthogonal matrix.
    \item
        \eqref{eq:Lem2-applied-again-pre}: Nominators were transferred into the $\sum_j$ terms.
    \item
        \eqref{eq:Lem2-applied-again}: Variables $\b{s}^m$ satisfy condition~\eqref{eq:suff} according
        to our assumptions.
    \item
        \eqref{eq:H(yi-last)}: We simplified the expression after squaring.
\end{enumerate}
Using this inequality, summing it for $i$, exchanging the order of
the sums, and making use of the orthogonality of matrix $\b{W}$,
we have
\begin{align*}
\sum_{i=1}^DH(y_i)&\ge\sum_{i=1}^D\left(\sum_{j\in\G^1}W_{i,j}^2H\left(s_j\right)+\ldots+\sum_{j\in\G^M}W_{i,j}^2H\left(s_j\right)\right)\\
&=\sum_{j\in\G^1}\left(\sum_{i=1}^DW^2_{i,j}\right)H\left(s_j\right)+\ldots+\sum_{j\in\G^M}\left(\sum_{i=1}^DW^2_{i,j}\right)H\left(s_j\right)\\
&=\sum_{j=1}^DH(s_j). \hspace*{11.2cm} \square
\end{align*}
\end{proof2}

\begin{cor}[ISA Separation Theorem]
ICA minimizes the l.h.s. of Equation~\eqref{eq:main-prop}, that is, it minimizes
\mbox{$\sum_{m=1}^M\sum_{i=1}^{d_m}H\left(y^m_i\right)$}. The set of minima is invariant to permutations and to
changes of the signs. Also, according to Proposition~\ref{R-prop}, $\{s^m_i\}$, that is, the coordinates of the
$\b{s}^m$ components of the ISA task belong to the set of the minima. $\square$
\end{cor}

\section{Kernel Covariance Technique for the ISA
Task}\label{sec:KC-deduction} For the sake of completeness, the
extension of the KC method \cite{gretton05kernel} for the ISA task
is detailed below. The extension is similar to the extensions
presented in Section~\ref{sec:KCCA} (see the KCCA method), and we
use the notations of that section.

First, we would like to measure the dependence of two 2 random variables $\b{u}\in\R^{d_1}$ and $\b{v}\in\R^{d_2}$. The
KC technique defines their dependence as their maximal covariance on the unit spheres $\S^{\b{u}}$, $\S^{\b{v}}$ of
function spaces $\F^{\b{u}}$, $\F^{\b{v}}$:

\[
J_{\text{KC}}(\b{u},\b{v},\F^{\b{u}},\F^{\b{v}}):=\sup_{g \in \S^{\b{u}},h \in \S^{\b{v}}} |E\{[g(\b{u})-E
g(\b{u})][h(\b{v})-Eh(\b{v})]\}|.
\]
This function $J_{\text{KC}}$ can be estimated empirically from $T$-element samples
$\b{u}_1,\ldots,\b{u}_T\in\R^{d_1}$, $\b{v}_1,\ldots,\b{v}_T\in\R^{d_2}$:

\[
J_{\text{KC}}^{\,emp}(\b{u},\b{v},\F^{\b{u}},\F^{\b{v}}):=\sup_{g \in \S^{\b{u}},h \in \S^{\b{v}}}
\left|\frac{1}{T}\sum_{t=1}^T[g(\b{u}_t)-\bar{g}][h(\b{v}_t)-\bar{h}]\right|.
\]
The estimation can be reduced to the following conditional maximization problem:
\begin{equation}
J_{\text{KC}}^{\,emp}(\b{u},\b{v},\F^{\b{u}},\F^{\b{v}})=\sup_{\b{c}_1^*\widetilde{\b{K}}^{\b{u}}\b{c}_1
\leq 1,\b{c}_2^*\widetilde{\b{K}}^{\b{v}}\b{c}_2 \leq 1}
\b{c}_1^*\widetilde{\b{K}}^{\b{u}}\widetilde{\b{K}}^{\b{v}}\b{c}_2.\label{eq:J-emp-KCCA}
\end{equation}
After the adaptation of the Lagrange multiplier technique and the
computation of the stationary points of \eqref{eq:J-emp-KCCA} it
can be realized that the values of $[\b{c}_1$;$\b{c}_2]$ and
$J^{\,emp}_{\text{KC}}$ can be computed as the solutions of the
generalized eigenvalue problem
\begin{equation}
  \begin{pmatrix}
    \widetilde{\b{K}}^{\b{u}} & \widetilde{\b{K}}^{\b{u}}\widetilde{\b{K}}^{\b{v}}\\
    \widetilde{\b{K}}^{\b{v}}\widetilde{\b{K}}^{\b{u}} & \widetilde{\b{K}}^{\b{v}}
  \end{pmatrix}
  \begin{pmatrix}
    \b{c}_1 \\
    \b{c}_2
  \end{pmatrix}
=\lambda
  \begin{pmatrix}
    \widetilde{\b{K}}^{\b{u}} & \b{0} \\
    \b{0} & \widetilde{\b{K}}^{\b{v}}
  \end{pmatrix}
\begin{pmatrix}
\b{c}_1\\
\b{c}_2
\end{pmatrix}.\label{eq:KC-gen-eig-problem-2}
\end{equation}

If the task is to measure the dependence between more than two
random variables $\b{y}^1 \in \R^{d_1} ,\ldots,$ \linebreak[4] $\b{y}^M \in
\R^{d_M}$ then \eqref{eq:KC-gen-eig-problem-2} is to be replaced
with the following generalized eigenvalue problem:

\[
\begin{pmatrix}
  \widetilde{\b{K}}^1 & \widetilde{\b{K}}^1\widetilde{\b{K}}^2 & \cdots & \widetilde{\b{K}}^1\widetilde{\b{K}}^M \\
  \widetilde{\b{K}}^2\widetilde{\b{K}}^1 &   \widetilde{\b{K}}^2 & \cdots & \widetilde{\b{K}}^2\widetilde{\b{K}}^M \\
  \vdots &           \vdots               &    & \vdots   \\
  \widetilde{\b{K}}^M\widetilde{\b{K}}^1 & \widetilde{\b{K}}^M\widetilde{\b{K}}^2 & \cdots & \widetilde{\b{K}}^M
\end{pmatrix}\begin{pmatrix}
  \b{c}_1 \\
  \b{c}_2 \\
  \vdots \\
  \b{c}_M
\end{pmatrix}=
\lambda\begin{pmatrix}
  \widetilde{\b{K}}^1 & \b{0} & \cdots & \b{0}\\
  \b{0} & \widetilde{\b{K}}^2 &  \cdots & \b{0}  \\
  \vdots &\vdots & &\vdots \\
  \b{0} & \b{0} & \cdots & \widetilde{\b{K}}^M\\
\end{pmatrix}\begin{pmatrix}
  \b{c}_1 \\
  \b{c}_2 \\
  \vdots\\
  \b{c}_M
\end{pmatrix}.
\]
Using the maximal eigenvalue of this problem, $J_{\text{KC}}$ can be estimated.

\bibliographystyle{splncs}

\begin{thebibliography}{10}

\bibitem{jutten91blind}
Jutten, C., Herault, J.:
\newblock Blind separation of sources: An adaptive algorithm based on
  neuromimetic architecture.
\newblock Signal Processing \textbf{24} (1991)  1--10

\bibitem{comon94independent}
Comon, P.:
\newblock Independent component analysis, a new concept?
\newblock Signal Processing \textbf{36} (1994)  287--314

\bibitem{hyvarinen01independent}
Hyv{\"a}rinen, A., Karhunen, J., Oja, E.:
\newblock Independent Component Analysis.
\newblock John Wiley \& Sons (2001)

\bibitem{cichocki02adaptive}
Cichocki, A., Amari, S.:
\newblock Adaptive Blind Signal and Image Processing.
\newblock John Wiley \& Sons (2002)

\bibitem{hyvarinen00emergence}
Hyv{\"a}rinen, A., Hoyer, P.O.:
\newblock Emergence of phase and shift invariant features by decomposition of
  natural images into independent feature subspaces.
\newblock Neural Computation \textbf{12} (2000)  1705--1720

\bibitem{cardoso98multidimensional}
Cardoso, J.:
\newblock Multidimensional independent component analysis.
\newblock In: Proceedings of International Conference on Acoustics, Speech, and
  Signal Processing (ICASSP '98). Volume~4., Seattle, {WA}, {USA} (1998)
  1941--1944

\bibitem{theis05blind}
Theis, F.J.:
\newblock Blind signal separation into groups of dependent signals using joint
  block diagonalization.
\newblock In: Proceedings of International Society for Computer Aided Surgery
  (ISCAS 2005), Kobe, Japan (2005)  5878--5881

\bibitem{akaho99MICA}
Akaho, S., Kiuchi, Y., Umeyama, S.:
\newblock {MICA}: Multimodal independent component analysis.
\newblock In: Proceedings of International Joint Conference on Neural Networks
  (IJCNN '99). (1999)  927--932

\bibitem{hyvarinnen06fastisa}
Hyv{\"a}rinen, A., K{\"o}ster, U.:
\newblock Fast{ISA}: A fast fixed-point algorithm for independent subspace
  analysis.
\newblock In: Proceedings of European Symposium on Artificial Neural Networks
  (ESANN 2006), Bruges, Belgium (2006)

\bibitem{vollgraf01multi}
Vollgraf, R., Obermayer, K.:
\newblock Multi-dimensional {ICA} to separate correlated sources.
\newblock In: Proceedings of Neural Information Processing Systems (NIPS 2001).
  Volume~14. (2001)  993--1000

\bibitem{bach03beyond}
Bach, F.R., Jordan, M.I.:
\newblock Beyond independent components: Trees and clusters.
\newblock Journal of Machine Learning Research \textbf{4} (2003)  1205--1233

\bibitem{strogbauer04least}
St{\"o}gbauer, H., Kraskov, A., Astakhov, S.A., Grassberger, P.:
\newblock Least dependent component analysis based on mutual information.
\newblock Physical Review E - Statistical, Nonlinear, and Soft Matter Physics
  \textbf{70} (2004)

\bibitem{poczos05independent1}
P{\'o}czos, B., L{\H{o}}rincz, A.:
\newblock Independent subspace analysis using k-nearest neighborhood distances.
\newblock Artificial Neural Networks: Formal Models and their Applications -
  ICANN 2005, pt 2, Proceedings \textbf{3697} (2005)  163--168

\bibitem{poczos05independent2}
P{\'o}czos, B., L{\H{o}}rincz, A.:
\newblock Independent subspace analysis using geodesic spanning trees.
\newblock In: Proceedings of International Conference on Machine Learning (ICML
  2005), Bonn, Germany (2005)  673--680

\bibitem{theis05multidimensional}
Theis, F.J.:
\newblock Multidimensional independent component analysis using characteristic
  functions.
\newblock In: Proceedings of European Signal Processing Conference (EUSIPCO
  2005). (2005)

\bibitem{theis06towards}
Theis, F.J.:
\newblock Towards a general independent subspace analysis.
\newblock In: Proceedings of Neural Information Processing Systems (NIPS 2006).
  (2006)

\bibitem{szabo06real}
Szab{\'o}, Z., L{\H{o}}rincz, A.:
\newblock Real and complex independent subspace analysis by generalized
  variance.
\newblock In: Proceedings of {ICA} Research Network International Workshop
  (ICARN 2006), Liverpool, {U.K.} (2006)  85--88
  {http://arxiv.org/abs/math.ST/0610438}.

\bibitem{nolte06identifying}
Nolte, G., Meinecke, F.C., Ziehe, A., M{\"u}ller, K.R.:
\newblock Identifying interactions in mixed and noisy complex systems.
\newblock Physical Review E \textbf{73} (2006)

\bibitem{pedersen07survey}
Pedersen, M.S., Larsen, J., Kjems, U., Parra, L.C.:
\newblock A survey of convolutive blind source separation methods.
\newblock In: Springer Handbook of Speech (to appear).
\newblock Springer Press (2007)

\bibitem{macdonald05derivation}
MacDonald, A., Cain, S.:
\newblock Derivation and application of an anisoplanatic optical transfer
  function for blind deconvolution of laser radar imagery.
\newblock Unconventional Imaging \textbf{5896} (2005)  9--20

\bibitem{hedgepeth99expectation}
Hedgepeth, J.B., Gallucci, V.F., O'Sullivan, F., Thorne, R.E.:
\newblock An expectation maximization and smoothing approach for indirect
  acoustic estimation of fish size and density.
\newblock {ICES} Journal of Marine Science \textbf{56} (1999)  36--50

\bibitem{vural06blind}
Vural, C., Sethares, W.A.:
\newblock Blind image deconvolution via dispersion minimization.
\newblock Digital Signal Processing \textbf{16} (2006)  137--148

\bibitem{douglas05natural}
Douglas, S.C., Sawada, H., Makino, S.:
\newblock Natural gradient multichannel blind deconvolution and speech
  separation using causal {FIR} filters.
\newblock IEEE Transactions on Speech and Audio Processing \textbf{13} (2005)
  92--104

\bibitem{mitianoudis03audio}
Mitianoudis, N., Davies, M.E.:
\newblock Audio source separation of convolutive mixtures.
\newblock {IEEE} Transactions on Speech and Audio Processing \textbf{11} (2003)
   489--497

\bibitem{roan03blind}
Roan, M.J., Gramann, M.R., Erling, J.G., Sibul, L.H.:
\newblock Blind deconvolution applied to acoustical systems identification with
  supporting experimental results.
\newblock The Journal of the Acoustical Society of America \textbf{114} (2003)
  1988--1996

\bibitem{araki03fundamental}
Araki, S., Makino, S., Mukai, R., Nishikawa, T., Saruwatari, H.:
\newblock Fundamental limitation of frequency domain blind source separation
  for convolved mixture of speech.
\newblock {IEEE} Transactions on Speech and Audio Processing \textbf{11} (2003)
   109--116

\bibitem{akyldiz02wireless}
Akyildiz, I.F., Su, W., Sankarasubramaniam, Y., Cayirci, E.:
\newblock Wireless sensor networks: a survey.
\newblock Computer Networks \textbf{38} (2002)  393--422

\bibitem{deligianni06source}
Deligianni, F., Lo, B., Yang, G.:
\newblock Source recovery for body sensor network.
\newblock In: Proceedings of International Workshop on Wearable and Implantable
  Body Sensor Networks 2006 (BSN 2006). (2006)  199--202

\bibitem{jung00independent}
Jung, T., Makeig, S., Lee, T., McKeown, M.J., Brown, G., Bell,
A.J., Sejnowski,
  T.J.:
\newblock Independent component analysis of biomedical signals.
\newblock In: Proceedings of International Workshop on Independent Component
  Analysis and Signal Separation (ICA 2000), Helsinki (2000)  633--644

\bibitem{glover99deconvolution}
Glover, G.H.:
\newblock Deconvolution of impulse response in event-related {BOLD fMRI}.
\newblock NeuroImage \textbf{9} (1999)  416--429

\bibitem{dyrholm06model}
Dyrholm, M., Makeig, S., Hansen, L.K.:
\newblock Model selection for convolutive {ICA} with an application to
  spatio-temporal analysis of {EEG}.
\newblock Neural Computation (2007)

\bibitem{kotzer98generalized}
Kotzer, T., Cohen, N., Shamir, J.:
\newblock Generalized projection algorithms with applications to optics and
  signal restoration.
\newblock Optics Communications \textbf{156} (1998)  77--91

\bibitem{karsli06further}
Karsl\i, H.:
\newblock Further improvement of temporal resolution of seismic data by
  autoregressive ({AR}) spectral extrapolation.
\newblock Journal of Applied Geophysics \textbf{59} (2006)  324--336

\bibitem{szabo06separation}
Szab{\'o}, Z., P{\'o}czos, B., L{\H{o}}rincz, A.:
\newblock Separation theorem for $\mathbb{K}$-independent subspace analysis
  with sufficient conditions.
\newblock Technical report, E\"otv\"os Lor\'and University, Budapest (2006)
  http://arxiv.org/abs/math.ST/0608100.

\bibitem{szabo06cross}
Szab{\'o}, Z., P{\'o}czos, B., L{\H{o}}rincz, A.:
\newblock Cross-entropy optimization for independent process analysis.
\newblock In: Proceedings of Independent Component Analysis and Blind Signal
  Separation (ICA 2006). Volume 3889 of LNCS., Springer (2006)  909--916

\bibitem{bach02kernel}
Bach, F.R., Jordan, M.I.:
\newblock Kernel independent component analysis.
\newblock Journal of Machine Learning Research \textbf{3} (2002)  1--48

\bibitem{rajagopal03multivariate}
Rajagopal, R., Potter, L.C.:
\newblock Multivariate {MIMO FIR} inverses.
\newblock IEEE Transactions on Image Processing \textbf{12} (2003)  458 -- 465

\bibitem{theis04uniqueness1}
Theis, F.J.:
\newblock Uniqueness of complex and multidimensional independent component
  analysis.
\newblock Signal Processing \textbf{84} (2004)  951--956

\bibitem{cover91elements}
Cover, T.M., Thomas, J.A.:
\newblock Elements of Information Theory.
\newblock John Wiley and Sons, New York, USA (1991)

\bibitem{fevotte03unified}
F\'{e}votte, C., Doncarli, C.:
\newblock A unified presentation of blind source separation for convolutive
  mixtures using block-diagonalization.
\newblock In: Proceedings of Independent Component Analysis and Blind Signal
  Separation (ICA 2003), Nara, Japan (2003)  349--354

\bibitem{meraim04algorithms}
Abed-Meraim, K., Belouchrani, A.:
\newblock Algorithms for joint block diagonalization.
\newblock In: Proceedings of European Signal Processing Conference (EUSIPCO
  2004). (2004)  209--212

\bibitem{fang90symmetric}
Fang, K.T., Kotz, S., Ng, K.W.:
\newblock Symmetric Multivariate and Related Distributions.
\newblock Chapman and Hall (1990)

\bibitem{takano95inequalities}
Takano, S.:
\newblock The inequalities of {F}isher information and entropy power for
  dependent variables.
\newblock In: Symposium on Probability Theory and Mathematical Statistics.
  (1995)

\bibitem{aronszajn50theory}
Aronszajn, N.:
\newblock Theory of reproducing kernels.
\newblock Transactions of the American Mathematical Society \textbf{68} (1950)
  337--404

\bibitem{wahba99support}
Wahba, G.:
\newblock Support vector machines, reproducing kernel hilbert spaces, and
  randomized {GACV}.
\newblock In: Advances in Kernel Methods, MIT Press (1999)  69--88

\bibitem{scholkopf99advances}
Sch{\"o}lkopf, B., Burges, C.J.C., Smola, A.J.:
\newblock Advances in Kernel Methods - Support Vector Learning.
\newblock MIT Press, Cambridge, {MA} (1999)

\bibitem{fukumizu07statistical}
Fukumizu, K., Bach, F.R., Gretton, A.:
\newblock Statistical consistency of kernel canonical correlation analysis.
\newblock Journal of Machine Learning Research \textbf{8} (2007)  361--383

\bibitem{gretton05kernel}
Gretton, A., Smola, A., Bousquet, O., Sch{\"o}lkopf, B.:
\newblock Kernel methods for measuring independence.
\newblock Journal of Machine Learning Research \textbf{6} (2005)  2075--2129

\bibitem{edelman98geometry}
Edelman, A., Arias, T., Smith, S.T.:
\newblock The geometry of algorithms with orthogonality constraints.
\newblock SIAM Journal on Matrix Analysis and Applications \textbf{20} (1998)
  303--353

\bibitem{lippert98nonlinear}
Lippert, R.A.:
\newblock Nonlinear Eigenvalue Problems.
\newblock PhD thesis, Massachusetts Institute of Technology (1998)

\bibitem{plumbley04lie}
Plumbley, M.D.:
\newblock Lie group methods for optimization with orthogonality constraints.
\newblock In: Proceedings of Independent Component Analysis and Blind Signal
  Separation (ICA 2004). Volume 3195 of LNCS., Springer (2004)  1245--1252

\bibitem{quinquis06efficient}
Quinquis, N., Yamada, I., Sakaniwa, K.:
\newblock Efficient dual {Cayley} parametrization technique for {ICA} with
  orthogonality constraints.
\newblock In: Proceedings of {ICA} Research Network International Workshop
  (ICARN 2006), Liverpool, {U.K.} (2006)  123--126

\bibitem{nishimori06riemannian}
Nishimori, Y., Akaho, S., Plumbley, M.D.:
\newblock Riemannian optimization method on the flag manifold for independent
  subspace analysis.
\newblock In: Proceedings of Independent Component Analysis and Blind Signal
  Separation (ICA 2006). Volume 3889 of LNCS., Springer (2006)  295--302

\bibitem{poczos06noncombinatorial}
P{\'o}czos, B., L{\H{o}}rincz, A.:
\newblock Non-combinatorial estimation of independent autoregressive sources.
\newblock Neurocomputing Letters \textbf{69} (2006)  2416--2419

\bibitem{amari96new}
Amari, S., Cichocki, A., Yang, H.H.:
\newblock A new learning algorithm for blind signal separation.
\newblock Advances in Neural Information Processing Systems \textbf{8} (1996)
  757--763

\bibitem{hyvarinen97fast}
Hyv{\"a}rinen, A., Oja, E.:
\newblock A fast fixed-point algorithm for independent component analysis.
\newblock Neural Computation \textbf{9} (1997)  1483--1492

\end{thebibliography}

\end{document}